\documentclass[letterpaper,11pt]{article}

\usepackage[utf8]{inputenc}
\usepackage[T1]{fontenc}
\usepackage{amsmath}
\usepackage{cfr-lm}

\usepackage{relsize}
\usepackage{exscale}

\usepackage{afterpage}

\usepackage{stmaryrd}
\usepackage{cite}
\usepackage[dvips]{graphicx}
\usepackage{graphpap}
\usepackage{ifthen}
\usepackage{enumerate}
\usepackage{amsmath}
\usepackage{amssymb}
\usepackage{mathrsfs}
\usepackage[errorshow]{tracefnt}
\usepackage{fancybox}
\usepackage{amsfonts}
\usepackage{amssymb}
\usepackage{multirow}
\usepackage{array}
\usepackage{mathtools}
\usepackage{soul}
\usepackage{pict2e}
\usepackage{bbold}
\DeclareMathAlphabet{\mathpzc}{OT1}{pzc}{m}{it}

\usepackage{cancel}
\usepackage{youngtab}

\Yboxdim{4pt}

\usepackage{lscape}
\usepackage{xypic}
\usepackage[dvipsnames]{xcolor}
\usepackage[
color,matrix,arrow]{xy}
\usepackage{setspace}

\usepackage{dynkin-diagrams}
\usepackage{caption}

\usepackage{epsf}

\usepackage{enumitem}

\usepackage{graphicx}

\renewcommand{\leq}{\leqslant}
\renewcommand{\geq}{\geqslant}

\makeatletter
\newcommand{\thickhline}{%
    \noalign {\ifnum 0=`}\fi \hrule height 1pt
    \futurelet \reserved@a \@xhline
}
\newcolumntype{'}{@{\hskip\tabcolsep\vrule width 1pt\hskip\tabcolsep}}
\makeatother
\newcolumntype{"}{@{\hskip\tabcolsep\vrule width 1.5pt\hskip\tabcolsep}}
\makeatother

\newcommand{\scr}{\mathscr}

\def\ie{{\it i.e.}}

\def\eg{{\it e.g.}}

\newcommand\BB{{\mathscr B}}
\newcommand\ZZ{{\mathbb Z}}

\def\small#1{{\hbox{$#1$}}}
\def\fraction#1{\small{1\over#1}}
\def\fr{\fraction}

\def\boxit#1{\vbox{\hrule\hbox{\vrule\kern3pt
             \vbox{\kern3pt#1\kern3pt}\kern3pt\vrule}\hrule}}

\newcommand{\Red}[1]{{\color{red} #1}}

\newcommand{\beq}{\begin{equation}}
\newcommand{\beqn}{\begin{equation*}}
\newcommand{\eeq}{\end{equation}}
\newcommand{\eeqn}{\end{equation*}}
\newcommand{\beqa}{\begin{eqnarray}}
\newcommand{\beqan}{\begin{eqnarray*}}
\newcommand{\eeqa}{\end{eqnarray}}
\newcommand{\eeqan}{\end{eqnarray*}}
\newcommand{\bdm}{\begin{displaymath}}
\newcommand{\edm}{\end{displaymath}}

\newcommand{\ba}{\begin{align}}
\newcommand{\ea}{\end{align}}

\newcommand\nn{\nonumber}

\newcommand\benu{\begin{enumerate}}
\newcommand\eenu{\end{enumerate}}
\newcommand\bit{\begin{itemize}}
\newcommand\eit{\end{itemize}}

\def\der'{\mathfrak{der}'\,}
\def\der{\mathfrak{der}\,}
\def\str'{\mathfrak{str}'\,}
\def\str{\mathfrak{str}\,}

\def\sl{\mathfrak{sl}}
\def\gl{\mathfrak{gl}}

\newcommand{\be}{\beta}

\newcommand{\dd}{{\mathsf{d}}}
\newcommand{\KK}{{\mathsf{K}}}

\newcommand{\Etendiagram}{
\begin{picture}(330,70)(45,-10)
\put(42,-10){$-1$}
\put(88,-10){$0$}
\put(128,-10){$1$}
\put(168,-10){$2$}
\put(208,-10){$3$}
\put(248,-10){$4$}
\put(288,-10){$5$}
\put(328,-10){$6$}
\put(368,-10){$7$}
\put(303,48){$8$}
\thicklines
\multiput(50,10)(40,0){9}{\circle{10}}
\multiput(55,10)(40,0){8}{\line(1,0){30}}
\put(290,50){\circle{10}}
\put(290,15){\line(0,1){30}}
\end{picture}
}

\newcommand{\Ediagram}{
\begin{picture}(330,70)(45,-10)
\put(44,-10){${\alpha_{-n}}$}
\put(78,-10){${\alpha_{-(n-1)}}$}
\put(118,-10){${\alpha_{-(n-2)}}$}
\put(246,-10){${\alpha_{4}}$}
\put(286,-10){${\alpha_{5}}$}
\put(326,-10){${\alpha_{6}}$}
\put(366,-10){${\alpha_{7}}$}
\put(303,48){${\alpha_{8}}$}
\thicklines
\put(50,10){\line(1,1){3.5}}
\put(50,10){\line(-1,1){3.5}}
\put(50,10){\line(1,-1){3.5}}
\put(50,10){\line(-1,-1){3.5}}
\multiput(50,10)(40,0){3}{\circle{10}}
\multiput(250,10)(40,0){4}{\circle{10}}
\multiput(55,10)(40,0){2}{\line(1,0){30}}
\put(135,10){\line(1,0){20}}
\put(165,10){\line(1,0){10}}
\put(185,10){\line(1,0){10}}
\put(205,10){\line(1,0){10}}
\put(225,10){\line(1,0){20}}
\put(255,10){\line(1,0){30}}
\put(295,10){\line(1,0){30}}
\put(335,10){\line(1,0){30}}
\put(290,50){\circle{10}}
\put(290,15){\line(0,1){30}}
\end{picture}
}

\newcommand{\EdiagramTwoGreyNodes}{
\begin{picture}(330,70)(45,-10)
\put(44,-10){${\beta_{-n}}$}
\put(78,-10){${\beta_{-(n-1)}}$}
\put(118,-10){${\beta_{-(n-2)}}$}
\put(246,-10){${\beta_4}$}
\put(286,-10){${\beta_5}$}
\put(326,-10){${\beta_6}$}
\put(367,-10){${\beta_7}$}
\put(303,48){${\beta_8}$}
\thicklines
\put(50,10){\line(1,1){3.5}}
\put(50,10){\line(-1,1){3.5}}
\put(50,10){\line(1,-1){3.5}}
\put(50,10){\line(-1,-1){3.5}}
\put(90,10){\line(1,1){3.5}}
\put(90,10){\line(-1,1){3.5}}
\put(90,10){\line(1,-1){3.5}}
\put(90,10){\line(-1,-1){3.5}}
\multiput(50,10)(40,0){3}{\circle{10}}
\multiput(250,10)(40,0){4}{\circle{10}}
\multiput(55,10)(40,0){2}{\line(1,0){30}}
\put(135,10){\line(1,0){20}}
\put(165,10){\line(1,0){10}}
\put(185,10){\line(1,0){10}}
\put(205,10){\line(1,0){10}}
\put(225,10){\line(1,0){20}}
\put(255,10){\line(1,0){30}}
\put(295,10){\line(1,0){30}}
\put(335,10){\line(1,0){30}}
\put(290,50){\circle{10}}
\put(290,15){\line(0,1){30}}
\end{picture}
}

\newcommand{\Arplusplusdiagram}{
\begin{picture}(210,70)(45,-10)
\put(42,-10){$-1$}
\put(88,-10){$0$}
\put(128,-10){$1$}
\put(168,-10){$2$}
\put(240,-10){$r-1$}
\put(182,50){$r$}
\thicklines
\multiput(50,10)(40,0){4}{\circle{10}}
\multiput(55,10)(40,0){3}{\line(1,0){30}}
\put(175,10){\line(1,0){25}}
\multiput(204,10)(8,0){2}{\line(1,0){4}}
\put(220,10){\line(1,0){25}}
\put(250,10){\circle{10}}
\put(170,50){\circle{10}}
\put(94,12){\line(2,1){71}}
\put(246,12){\line(-2,1){71}}
\end{picture}
}

\newcommand{\AOneplusplusdiagram}{
\begin{picture}(90,70)(45,-10)
\put(42,-10){$-1$}
\put(88,-10){$0$}
\put(128,-10){$1$}
\thicklines
\multiput(50,10)(40,0){3}{\circle{10}}
\put(55,10){\line(1,0){30}}
\put(95,12){\line(1,0){30}}
\put(95,8){\line(1,0){30}}
\end{picture}
}

\newcommand{\ad}{\mathrm{ad}\,}

\def\fg{{\mathfrak g}}

\def\fa{{\mathfrak a}}

\def\sh{\sharp}

\def\*{\partial}

\def\id
{\mathbb{1}}

\def\shift#1#2{\underset
  {\scriptscriptstyle\Red{#1}}{#2{}^{\mathstrut}_{\mathstrut}}}

\def\adj{\hbox{\bf adj}}
\def\rank{\hbox{rank}}

\def\EWeight#1#2#3#4{\bigl({}^{\mathstrut}_{#1\mathstrut}{}_{#2\mathstrut}^{#4\mathstrut}{}_{#3\mathstrut}^{\mathstrut}\bigr)}

\def\EEWeight#1#2#3#4{\bigl({}^{\mathstrut}_{#1\mathstrut}{}_{#2\mathstrut}^{\hspace{-6pt}#4\hspace{-1pt}\mathstrut}{}_{#3\mathstrut}^{\mathstrut}\bigr)}

\def\PP{{\mathbb P}}

\def\be{\begin{equation}}
\def\ee{\end{equation}}

\def\II{{\mathscr I}}

\def\dd{{\sf d}}
\def\KK{{\sf K}}

\def\AA{{\mathscr A}}
\def\BB{{\mathscr B}}

\def\Cn#1{\hskip1pt{\buildrel{\scriptscriptstyle #1}\over c}\hskip1pt}

\def\compl{{\scriptscriptstyle\complement}}

\def\deltaeq{\overset{\scriptscriptstyle\Delta}=}

\def\curlybrack#1{{\{\mskip-.8\thinmuskip#1\mskip-.6\thinmuskip\}}}
\def\MP{\curlybrack M}
\def\NP{\curlybrack N}
\def\PP{\curlybrack P}
\def\QP{\curlybrack Q}

\def\shift#1#2{\underset
  {\scriptscriptstyle\Red{#1}}{#2{}^{\mathstrut}_{\mathstrut}}}

\makeatletter
\DeclareRobustCommand{\loplus}{\mathbin{\mathpalette\dog@lsemi{+}}}
\DeclareRobustCommand{\lotimes}{\mathbin{\mathpalette\dog@lsemi{\times}}}
\DeclareRobustCommand{\roplus}{\mathbin{\mathpalette\dog@rsemi{+}}}
\DeclareRobustCommand{\rotimes}{\mathbin{\mathpalette\dog@rsemi{\times}}}

\newcommand{\dog@rsemi}[2]{\dog@semi{#1}{#2}{-90,90}}
\newcommand{\dog@lsemi}[2]{\dog@semi{#1}{#2}{270,90}}
\newcommand{\dog@semi}[3]{%
  \begingroup
  \sbox\z@{$\m@th#1#2$}%
  \setlength{\unitlength}{\dimexpr\ht\z@+\dp\z@\relax}%
  \makebox[\wd\z@]{\raisebox{-\dp\z@}{%
    \begin{picture}(1,1)
    \linethickness{\variable@rule{#1}}
    \roundcap
    \put(0.5,0.5){\makebox(0,0){\raisebox{\dp\z@}{$\m@th#1#2$}}}
    \put(0.5,0.5){\arc[#3]{0.5}}
    \end{picture}%
  }}%
  \endgroup
}
\newcommand{\variable@rule}[1]{%
  \fontdimen8  
  \ifx#1\displaystyle\textfont3\else
    \ifx#1\textstyle\textfont3\else
      \ifx#1\scriptstyle\scriptfont3\else
        \scriptscriptfont3\relax
  \fi\fi\fi
}
\makeatother

\RequirePackage{hyperref}

\numberwithin{equation}{section}

\begin{document}


\pgfkeys{/Dynkin diagram, edge length=1cm,
fold radius=.6cm, 
root-radius=.11cm,
indefinite edge/.style={
draw=black, fill=white, dotted,
thin}}

\frenchspacing

\null\vspace{-18mm}

\includegraphics[height=2cm]{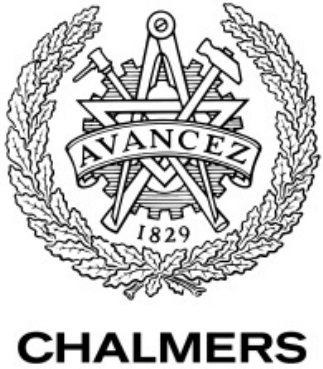}
\hspace{2mm}
\includegraphics[height=1.85cm]{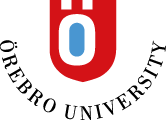}

\vspace{-17mm}
{\flushright Gothenburg preprint \\ 
\today\\}

\vspace{4mm}

\hrule

\vspace{16mm}


\pagestyle{empty}

\begin{center}
  {\Large \bf \sc Tensor hierarchy algebra extensions}
  \\[3mm]
  {\Large \bf \sc of over-extended Kac--Moody algebras}
    \\[10mm]
    
{\large
Martin Cederwall${}^1$ and Jakob Palmkvist${}^{2}$}

\vspace{10mm}
       {\footnotesize ${}^1${\it Department of Physics,
         Chalmers Univ. of Technology,\\
 SE-412 96 Gothenburg, Sweden}}

\vspace{2mm}
       {\footnotesize ${}^2${\it Department of Mathematics,
         \"Orebro Univ.,\\
 SE-701 82 \"Orebro, Sweden}}

\end{center}

\vfill

\begin{quote}
  
\textbf{Abstract:} 
Tensor hierarchy algebras are infinite-dimensional generalisations of Cartan-type Lie superalgebras.
They are not contragredient, exhibiting
an asymmetry between positive and negative levels.
These superalgebras have been a focus of attention due to the fundamental
r\^ole they play for extended geometry. 
In the present paper, we examine tensor hierarchy
algebras which are super-extensions of over-extended (often, hyperbolic) Kac--Moody algebras.
They contain novel algebraic structures. Of particular interest is the extension
of a over-extended algebra by its fundamental module, an extension that
contains and generalises the extension of an affine Kac--Moody algebra by a
Virasoro derivation $L_1$.
A conjecture about the complete superalgebra is formulated, relating
it to the corresponding Borcherds superalgebra.

\end{quote} 

\vfill

\hrule

\noindent{\tiny email:
  martin.cederwall@chalmers.se, jakob.palmkvist@oru.se}

\newpage

\tableofcontents

\pagestyle{plain}

\section{Introduction}

Tensor hierarchy algebras
\cite{Palmkvist:2013vya,Carbone:2018xqq}
constitute a class of simple non-contragredient Lie superalgebras, whose finite-dimensional members are 
those of Cartan type
in the classification of Kac \cite{Kac77B}.
They were originally invented due to the need to accommodate the
embedding tensor of gauged supergravities in the algebra \cite{Palmkvist:2013vya,Greitz:2013pua}, and
have been the subject of attention recently due to their importance in extended geometry
\cite{Cederwall:2017fjm,Cederwall:2018aab,Cederwall:2019qnw,Cederwall:2019bai}.
It had already been clear that 
certain classes of Lie superalgebras were relevant for the gauge structure and dynamics of extended geometry
\cite{Palmkvist:2015dea,Cederwall:2015oua,Cederwall:2017fjm,Cederwall:2018aab}, but it turns out that the tensor hierarchy algebras provide all information necessary, so that they may indeed be viewed as basic building blocks of extended geometry.
This becomes increasingly clear when situations with infinite-dimensional structure groups are encountered
\cite{Bossard:2017wxl,Bossard:2017aae,Bossard:2018utw,Bossard:2019ksx}.

Extended geometry contains double geometry
\cite{Hull:2004in,Hull:2009mi,Berman:2014jba,Cederwall:2014kxa,Cederwall:2014opa,Cederwall:2016ukd}
and exceptional geometry
\cite{Hull:2007zu,Berman:2010is,Berman:2012vc,Cederwall:2013naa,Cederwall:2013oaa,Aldazabal:2013mya,Blair:2013gqa,Hohm:2013vpa,Hohm:2013uia,Hohm:2014fxa,Cederwall:2015ica,Butter:2018bkl,Bossard:2017aae,Bossard:2018utw,Bossard:2019ksx}
as special cases.
The hope of understanding extended geometry for infinite-dimensional, in particular hyperbolic, structure groups is the main motivation of the present work,
and we hope that it will lead to a reformulation of gravity or supergravity
where the hyperbolic Belinskii--Khalatnikov--Lifshitz group 
\cite{BKL,Damour:2001sa,Damour:2002cu}
not only emerges in extreme situations, but is an integral part of the formulation of the theory.
This may eventually put the $E_{10}$
\cite{Damour:2002cu}
and $E_{11}$ \cite{West:2001as,West:2017vhm,Bossard:2017wxl} proposals on a firm ground, and provide a mechanism for the emergence of space(-time).

The subject of the present paper is an attempt to take a step in this direction. The results are purely mathematical, but have obvious bearing to extended geometry.
An extensive summary is provided in Section \ref{SummarySection}.

\section{Summary and discussion\label{SummarySection}}

The aim of this Section is to give an extensive summary of the lines of thought, methods and results of this paper, with some minimum of technical detail, but without all the calculations of the following Sections. 
Hopefully, this organisation of the paper will make it more accessible to readers, especially those not inclined to indulge in technical issues concerning generators and Serre relations of non-standard type.
We also comment on physical applications, and on the possible continued investigation of tensor hierarchy extensions of ``very extended'' Kac--Moody algebras, such as $E_{11}$ or $A_1^{+++}$. 

\begin{figure}
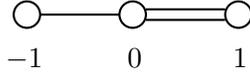

\begin{center}
\AOneplusplusdiagram
\caption{\it Dynkin diagram of $A_1^{++}$, with the convention for numbering of nodes.}  \label{A1-figur}  
\end{center}
\end{figure}

\begin{figure}
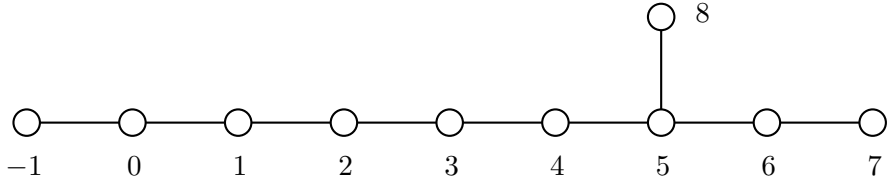

\begin{center}
\Etendiagram
\caption{\it Dynkin diagram for $E_{10}$, with the convention for numbering of nodes.} \label{E-figur}  
\end{center}
\end{figure}

The goal of the paper is one and clear: to provide a working (tensor) formalism to deal with 
tensor hierarchy algebras $S(\fg^{++})$, where $\fg^{++}$ is an over-extended Kac--Moody algebra. 
This means
that it is the further extension of an affine Kac--Moody algebra $\fg^+$.
Often it is a hyperbolic, but not always
(the standard counterexamples being $A_r^{++}$, with Dynkin diagrams shown in Figures \ref{A1-figur} and \ref{A-figur}, which are hyperbolic only for $r\leq8$).
We are particularly interested
in the case where $\fg=E_8$, which means that $\fg^+=E_9$ and $\fg^{++}=E_{10}$ (indeed hyperbolic), with Dynkin diagram shown in Figure \ref{E-figur}.

The original rationale for introducing tensor hierarchy algebras \cite{Palmkvist:2013vya}
was the need for an algebraic foundation for gauged supergravities---the embedding tensor is a set of elements in this superalgebra, and this is the basic reason for the asymmetry between positive and negative levels. 
The asymmetry means that the tensor hierarchy algebras are not contragredient and is the main difference compared to the related Borcherds superalgebras
that had been considered previously \cite{Henry-Labordere:2002xau,Henneaux:2010ys,Palmkvist:2011vz}.
It was then realised that extended geometry \cite{Cederwall:2017fjm,Cederwall:2018aab} in general can be understood as founded on the algebraic structure of tensor hierarchy algebras
\cite{Cederwall:2019qnw,Cederwall:2019bai}. Tensor hierarchy algebras are the unique objects that capture all fields and local symmetries in extended geometry in a correct way, and this becomes increasingly obvious when the structure group becomes infinite-dimensional.

One main goal for this line of research is to formulate extended geometry with hyperbolic structure groups. This would amount to a reformulation of gravity, where the hyperbolic BKL symmetry \cite{BKL}
becomes a built-in feature, in a well-defined sense.
However, in order to achieve this, one needs the tensor hierarchy algebra extension of a ``very extended'' Kac--Moody algebra $\fg^{+++}$; the superalgebra underlying extended geometry is a super-extension of the already once extended Kac--Moody algebra. 
So far, we only have partial results for such algebras \cite{CederwallPalmkvistForth}.
The present paper is a step on the way, investigating the presumably much better controlled class of superalgebras $S(\fg^{++})$.
We find it interesting and instructive to understand the novel mathematical structures arising already here, and they will certainly form a basis for the continued investigations.

The definition of a tensor hierarchy algebra $S(\fa)$ (where $\fa$ is a generic Kac--Moody algebra) in terms of generators and relations
leads, when the algebra $\fa$ is infinite-dimensional, to the appearance of ``extra'' generators.
By ``extra'' generators, we understand generators that appear at non-negative levels in the tensor hierarchy algebra, but not in the corresponding Borcherds superalgebra $\BB(\fa)$.
The phenomenon is quite analogous to how the generators and relations for an affine Kac--Moody algebra allow the continued brackets with raising operators to produce the loop generators.
In fact, as will be displayed explicitly in Section \ref{SAffineSection}, the stage where this happens is more or less the same as for the affine algebra. The simplest instance is indeed for the super-extension of an affine algebra, then producing a ``scalar'' generator, which can be identified with a Virasoro generator $L_1$, providing a derivation on the algebra.
This operator has a direct interpretation in the physical application to extended geometry
\cite{Cederwall:2019bai}.

For the case at hand, when an over-extended Lie algebra $\fg^{++}$ is extended to 
a tensor hierarchy algebra
$S(\fg^{++})$, an even more interesting structure appears. Not only do extra modules appear at all non-negative levels, the extra module appearing at level $0$ then is a fundamental of $\fg^{++}$,
which is the lowest weight module $R(-\lambda)$ with lowest weight $-\lambda$, $\lambda$ being the fundamental weight associated with the leftmost node in $\fg^{++}$.
 We are led to a structure which is
at the heart of the superalgebra, namely the (bosonic) extension of an over-extended Kac--Moody algebra 
$\fg^{++}$ by generators spanning its fundamental module $F=R(-\lambda)$.
The Lie algebra structure on $\AA=\fg^{++}\oplus F$ is investigated in detail in Section
\ref{ExtendByFund}.
Let $\{T_\alpha\}$ be a basis of generators of $\fg^{++}$ and $\{J_M\}$ a basis for its fundamental module. We then find bracket relations
\begin{align}
  [T_\alpha,T_\beta]&=f_{\alpha\beta}{}^\gamma T_\gamma\;,\nn\\
  [T_\alpha,J_M]&=-t_{\alpha M}{}^NJ_N+u_{M\alpha}{}^\beta T_\beta\;,\label{TJBracketsIntro}\\
  [J_M,J_N]&=g_{MN}{}^PJ_P\;,\nn
\end{align}
where $f$ and $t$ are structure constants and representation matrices, and $u$ and $g$ ``new'' sets of coefficients. The Jacobi identities of course restrict these structure constants through bilinear relations.
Had the $J$'s been a collection of ``scalars'', so that the second of these equations only had read
``$[J,T]=-uT$'', the $J$'s would have been (outer) derivations of the algebra $\fg^{++}$ and $g$ the structure constants of a Lie algebra of derivations. Now the $J$'s also transform in the fundamental representation. We call this structure ``transforming derivations''.
Infinitely many of the generators in the fundamental---the ones with light-like weights which can be brought to the lowest weight by a Weyl transformation---are $L_1$'s corresponding to some embedding of an affine Lie
algebra $\fg^+\subset\fg^{++}$, but there are also the ones with time-like weights, together filling up the whole fundamental module. Among simple Lie algebras, the over-extended ones are the ``smallest'' ones allowing a non-trivial extension of this kind.

The Lie algebra $\AA$ turns out to have interesting representations, among them the fundamental module $F$ itself. There are representation matrices $j_{MN}{}^P$,
such that
\begin{align}
  [t_\alpha,t_\beta]&=f_{\alpha\beta}{}^\gamma t_\gamma\;,\nn\\
  [t_\alpha,j_M]&=-t_{\alpha M}{}^Nj_N
      +u_{M\alpha}{}^\beta t_\beta\;,\label{FundIds}\\
  [j_M,j_N]&=g_{MN}{}^Pj_P\;.\nn
 \end{align}
The structure constants and fundamental representation matrices are related as $g_{MN}{}^P=-2j_{[MN]}{}^P$.

The structure constants in eq. \eqref{TJBracketsIntro} and the identities among them expressing the Jacobi identities, together with the representation matrices of the fundamental module, turn out to provide all information we need in order to start building the superalgebra $S(\fg^{++})$ in Section
\ref{SgplusplusSection}. The structure constants of $S(\fg^{++})$ are constructed from this set of ``tensors''.

An important experimental arena is provided by certain gradings. 
Of special interest is the double grading with respect to the fermionic node and one more node, chosen so that the remaining Dynkin diagram is of type $A$. In the $E$ series, for example, the second node is the exceptional node (node 8 in Figure \ref{E-figur}), and in the $A$ series one of the nodes in the circle attached to the affine node
(node $r$ in Figure \ref{A-figur}. The important observation is that each grade with respect to this second node will consist of a module of a finite-dimensional Cartan-type superalgebra (tensor hierarchy algebra) 
$W(A_{d-1})\simeq W(d)$. This is described in Sections \ref{GLGradingSection} and 
\ref{ExtraModulesSubSection}.
Such $W(d)$ modules have the property that they, as vector spaces, consist of sums of elements which are formed as the tensor product of some $\gl(d)$ module with the sum of all form (antisymmetric) modules.
Even without knowing the precise content of the $W(d)$ modules, this constraint alone, together with the observation that all elements at non-negative levels in the corresponding Borcherds superalgebra also appear in the tensor hierarchy algebra, allows us to ``discover'' many extra modules experimentally, which are practically unreachable using the basic definitions.
Pictures of this double grading for $S(E_9)$,  $S(E_{10})$ and $S(A_1^{++})$ are given in Figures
\ref{SE9Figure},  \ref{SE10Figure} and \ref{SA1++Figure}.

The conjecture of Section \ref{SgplusplusSection}, that level $\ell$ in $S(\fg^{++})$ consists of 
\be \label{Rell}
R_\ell=\BB_\ell\oplus\BB_{\ell+1}\;,
\ee
where $\BB_\ell$ is level $\ell$ in the Borcherds superalgebra $\BB(\fg^{++})$, is corroborated both by our tensorial formalism and by the predictions of the 
$\gl$ gradings. In the case $\ell=0$, eq.~(\ref{Rell}) gives back the decomposition $\scr A= \fg^{++} \oplus F$ described above.

The continued work will focus on the tensor hierarchy extensions $S(\fg^{+++})$ of very extended 
Kac--Moody algebras, \eg\ $S(E_{11})$. We expect these superalgebras to display a much more unkempt structure, and it will probably no longer be possible to make statements to all levels.
An investigation of the local superalgebra in some grading should still be possible. In particular, the tools developed in the present paper will be useful in a $\fg^{++}$-covariant double grading, which is also precisely the grading relevant for application to extended geometry.
Some information about extended geometry with hyperbolic structure group is already revealed by the results of the present paper. For example, there will most likely be fields not only in a coset
of the hyperbolic group, but also in the fundamental module. Parameters of local symmetries will appear not only in the fundamental module, but also in $\BB_2$. 

\begin{figure}
\begin{center}
\Arplusplusdiagram
\caption{\it Dynkin diagram of $A_r^{++}$, $r\geq2$, with the convention for numbering of nodes.}\label{A-figur}
\end{center}
\end{figure}

\section{Some notation}
We let $\fg$ be a finite-dimensional simple Lie algebra, $\fg^+$ its (untwisted) affine extension, $\fg^{++}$ the further extension to what in most interesting cases is a hyperbolic Kac--Moody algebra, etc.
By $\fg^+$ we understand the loop algebra centrally extended by $\KK$ (but not $\dd\sim L_0$).
A generic Kac--Moody algebra is occasionally denoted $\fa$.

The nodes in $\fg$ are numbered $1,\ldots,r$, with $r=\rank(\fg)$.
Node $0$ is the affine node, node $-1$ the ``hyperbolic'' node etc. This means that the fermionic node for 
$\BB(\fg^{(n)})$ or $S(\fg^{(n)})$ (superscript $(n)$ meaning $n$ $+$ signs) will have number $-n$.
Often, we also distinguish fermionic raising and lowering operators from bosonic ones by writing  $\epsilon_{-n}$ instead of $e_{-n}$ and $\phi_{-n,\alpha}$ instead of $f_{-n,\alpha}$.

When using tensor notation for modules,
\begin{itemize}
\item $A,B,\ldots$ are indices for the adjoint of $\fg$,
\item $\mu,\nu,\ldots$ are indices for the fundamental module $f$ of the affine algebra $\fg^+$,
\item $\alpha,\beta,\ldots$ are adjoint indices for the over-extended Kac--Moody algebra $\fg^{++}$,
\item $M,N,\ldots$ are indices for the fundamental module $F$ of $\fg^{++}$.
\end{itemize}
By fundamental, we understand lowest weight modules, so that
$F=R(-\Lambda_{-1})$ (the module with lowest weight $-\Lambda_{-1}$), $\bar F=R(\Lambda_{-1})$.

The leading irreducible $\fg^{++}$ module (the one with lowest lowest weight) in the symmetric product 
$\vee^2F$ is denoted $S_2^\compl$, the remaining part (the one appearing at level $2$ in $\BB(\fg^{++})$) $S_2$. The leading module in  $\wedge^2F$ is denoted $A_2^\compl$, the remaining part (the one appearing at level $2$ in $\fg{^{+++}}$) $A_2$. 
Analogously, for products of the affine fundamentals (lowest weight modules at $k=-1$) we use lowercase letters:
$f$, $s_2^\compl$, $s_2$, $a_2^\compl$, $a_2$.

Level $\ell$ in the grading with respect to the fermionic node of a Borcherds superalgebra $\BB$ is denoted $\BB_\ell$.
Level $\ell$ in a tensor hierarchy algebra $S$ is denoted $R_\ell$.

Brackets with more than two arguments are built from binary brackets from right to left, for example,
\begin{align}
[x,y,z,w]\deltaeq[x,[y,|z,w]]]\;.
\end{align} 
The bracket $[x,y]$ of two elements may be symmetric or antisymmetric depending on the $\mathbb{Z}_2$-degree of the elements. 

The algebras and superalgebras in this paper can be viewed as algebras over the complex numbers, but since the definitions we employ use root and weight decompositions, one may also view them as algebras over the real numbers, the (split) real form being defined by these bases. The latter is relevant for applications to extended geometry.

\section{Tensor hierarchy algebras: definitions and generalities
\label{THAGeneralSection}}

Tensor hierarchy algebras were
 originally introduced in ref.~\cite{Palmkvist:2013vya} as extensions of simple finite-dimensional complex Lie algebras $\fg$
(with particular focus on the cases $\fg=E_n$ for $4\leq n \leq 8$).
In ref.~\cite{Carbone:2018xqq}, a construction from generators and relations was given, and the notation $S(\fg)$ was introduced for the tensor hierarchy algebras
considered in ref.~\cite{Palmkvist:2013vya}.
In addition to $\fg$ itself, the data defining $S(\fg)$ also consists of 
a simple root $\alpha_1$ of $\fg$.
The data can be visualised by a Dynkin diagram where a ``grey'' node $0$ is connected to node $1$, corresponding to $\alpha_1$,
in the Dynkin diagram of $\fg$.
Thus, for any choice of 
a simple finite-dimensional complex Lie algebras $\fg$ and a simple root of it, there is an associated tensor hierarchy algebra $S(\fg)$.
There are also two other related Lie superalgebras associated to this data, denoted by $\scr B(\fg)$ and $W(\fg)$ in ref.~\cite{Carbone:2018xqq}.
The Lie superalgebra $\scr B(\fg)$ is a contragredient Lie superalgebra (more specifically, a Borcherds superalgebra) 
which means that it can be 
constructed from
the Dynkin diagram by generators and relations
in a way which is symmetric (up to signs) under the exchange of raising and lowering operators \cite{Kac77B}. 
The Lie superalgebra $W(\fg)$ (also called tensor hierarchy algebra) can be considered as 
a unification of $S(\fg)$ and $\scr B(\fg)$. Both $S(\fg)$ and $W(\fg)$ are constructed by a modification of the usual generators and relations
for $\scr B(\fg)$, from the same Dynkin diagram. This modification breaks the symmetry between raising and lowering operators and thus makes the tensor hierarchy algebras non-contragredient. In ref.~\cite{Carbone:2018xqq} it was shown that when $\fg=A_{n-1}$ and node $1$ is one of the outermost nodes, the tensor hierarchy algebras $W(\fg)$ and $S(\fg)$ are finite-dimensional and coincide with the well known
Lie superalgebras $W(n)$ and $S(n)$ of Cartan type \cite{Kac77B}.

An advantage with the construction of $S(\fg)$ from generators and relations is that $\fg$ can be generalised to any Kac--Moody algebra,
not necessarily a finite-dimensional one. In \cite{Cederwall:2019qnw,Cederwall:2019bai}, we considered $S(\fg^+)$, where $\fg^+$ is a Kac--Moody algebra
obtained by adding a node to the Dynkin diagram of $\fg$
(generalising ``node $1$'' in the description of $S(\fg)$ above, but ``node $0$'' with the numbering of nodes that we will use below).
It can be attached to one or more nodes in the Dynkin diagram of $\fg$, with one or more lines.
The information to which nodes 
it is attached, and with how many lines,
is given by an integral dominant
weight $\lambda$ of $\fg$.
For a particular choice of $\lambda$, the Kac--Moody algebra $\fg^+$ is the affine extension of $\fg$.
In the present paper, we focus on this choice, but we also go one step further, and consider $S(\fg^{++})$, where $\fg^{++}$
is the ``over-extension'' of $\fg$, usually hyperbolic.
Particularly, if $\fg=E_8$ and $\lambda$ is the fundamental weight associated to the ``outermost'' node, then $\fg^+=E_9$ and $\fg^{++}=E_{10}$. 

In the present paper, the nodes in the Dynkin diagram of $\fg$ are numbered $1,2,\ldots,r$, where $r$ is the rank of $\fg$.
When we extend $\fg$ to the affine Lie algebra $\fg^+$, we add a node $0$ (``the affine node'') and when we 
extend $\fg^+$ to $\fg^{++}$, we add a node $-1$ (``the hyperbolic node''). In this way, we can define a Kac--Moody algebra $\fg^{(n)}$
of rank $r+n$, where
the superscript $(n)$ means $n$ plus signs, for any integer $n\geq1$, starting from the simple finite-dimensional Lie algebra $\fg$
and its affine extension $\fg^+$.
Node $-n+1$ is connected with a single line to node $-n+2$, but disconnected from nodes $-n+3,\ldots,r$.
If $A^{(n)}_{ij}$ is the Cartan matrix of $\fg^{(n)}$, where $i,j=-n+1,\ldots,r$, this means that 
\begin{align}
A^{(n)}_{-n+1,j}=A^{(n)}_{j,-n+1}=
\begin{cases}
\phantom{-}2 \qquad \text{if} \qquad j=-n+1 \\
-1 \qquad \text{if} \qquad j=-n+2 \\
\phantom{-}0 \qquad \text{if} \qquad j\geq-n+3 
\end{cases}
\end{align}
and $A^{(n)}_{ij}=A^{(n-1)}_{ij}$ if $i,j\neq -n+1$. Within the Dynkin diagram of $\fg$, we allow for multiple lines and arrows; it does not need to be simply laced.

The Dynkin diagram of $\fg^{(n)}$ can be extended to the Dynkin diagram of the Lie superalgebras $S(\fg^{(n)})$, $W(\fg^{(n)})$ and $\scr B(\fg^{(n)})$
by replacing node $-n$ in the Dynkin diagram of $\fg^{(n+1)}$ with a ``grey node'' (drawn as $\otimes$). This means two things: first (as we will see below), that
there are associated generators that are odd elements in the Lie superalgebra (and thus we may call it a ``fermionic node'' as opposed to the ordinary 
white or ``bosonic'' nodes), second, that the associated diagonal entry in the Cartan matrix is not $2$, but $0$.
More precisely, if $B^{(n+1)}_{ij}$ is the Cartan matrix of $S(\fg^{(n)})$, $W(\fg^{(n)})$ or $\scr B(\fg^{(n)})$,
where $i,j=-n,\ldots,r$, then 
\begin{align}
B^{(n+1)}_{-n,j}=B^{(n+1)}_{j,-n}=
\begin{cases}
\phantom{-}0 \qquad \text{if} \qquad j=-n \\
-1 \qquad \text{if} \qquad j=-n+1 \\
\phantom{-}0 \qquad \text{if} \qquad j\geq-n+2 
\end{cases}
\end{align}
and $B^{(n+1)}_{ij}=A^{(n)}_{ij}$ if $i,j\neq -n$.

Note that the numbering of nodes here, where node $0$ always is the affine node,
is different from the ones in refs. \cite{Palmkvist:2013vya,Carbone:2018xqq,Cederwall:2019qnw,Cederwall:2019bai}, where node $0$ always is the grey node, as reviewed above.
We will mostly in this paper not consider the situation with a general $n$, but focus on the cases $n=1$ and $n=2$. Note that our construction involves a specific position for the grey node, \ie, a choice of $\lambda$. In principle, other tensor hierarchy algebras corresponding to other choices of $\lambda$ can be constructed, but are not expected to share many of the attractive features of the superalgebras presently considered.

We note that the determinant of the Cartan matrix is
$|B|=-|A^-|$,
where $A^-$ is the Cartan matrix for an algebra where the node connecting to the fermionic one is removed, \ie,
\be
|B^{(n+1)}|=-|A^{(n-1)}|\;.
\ee
The matrix $B^{(n+1)}$ is non-singular for $n=0$, but singular for $n=1$.
The Cartan matrix for the super-extension of an affine Kac--Moody algebra is non-singular without further bosonic extension, while the Cartan matrix for the super-extension of a hyperbolic Kac--Moody algebra is singular. This is associated to a difficulty with the definition of $W(\fg^{++})$; we will return to this issue in Section \ref{LocalSubalgebraSection}.

\subsection{Generators and relations\label{GenRelSection}} 

We will now explain how to construct the Lie superalgebras $S(\fg^{(n)})$, $W(\fg^{(n)})$ or $\scr B(\fg^{(n)})$
from its Dynkin diagram described above, or from the associated Cartan matrix $B_{ij}$,
where $i,j=-n,\ldots,r$ (we drop the superscript $(n+1)$ on $B_{ij}$).
 
We start with $\scr B(\fg^{(n)})$ and the well known construction of it as a contragredient Lie superalgebra.
To each node $i=-n,\ldots,r$ in the Dynkin diagram we associate three generators $h_i,e_i,f_i$. 
Among these $3r$ generators, $e_{-n}$ and $f_{-n}$ are odd elements, whereas all the other generators
(including $h_{-n}$) are even. 

Then 
$\scr B(\fg^{(n)})$ is defined as the Lie superalgebra generated by all $h_i,e_i,f_i$ modulo the 
Chevalley--Serre relations
\begin{gather}
[h_i,e_j]=B_{ij}e_j\;,  \qquad   
[h_i,f_j]=-B_{ij}f_j\;, \qquad   
[e_i,f_j]=\delta_{ij}h_j\;,  \label{eigen1}\\
({\rm ad}\,e_i)^{1-B_{ij}}(e_j)=({\rm ad}\,f_i)^{1-B_{ij}}(f_j)=0\;\quad (\text{for } B_{ij}\leq0)\;.
\label{serre1}
\end{gather}
 
When we go from $\scr B(\fg^{(n)})$ to $S(\fg^{(n)})$ we modify the set of generators. We
keep all the generators associated to the white nodes $-n+1,\ldots,r$, as well as $e_{-n}$, 
but we remove $h_{-n}$ and we replace
the single odd generator $f_{-n}$ with $n+r-1$ odd generators $f_{-n,i}$, where $i=-n+2,\ldots,r$.
For the relations, we keep (\ref{eigen1}) and (\ref{serre1}) above (but whenever $f$ appears
with a single index as subscript,
it is understood that it does not take the value $-n$ anymore).
We then add the following relations to those above,
\begin{gather}
[e_{-n},f_{-n,i}]=h_i\;, \quad\!\!
[h_i,f_{-n,j}]=-B_{i,-n}f_{-n,j}\;, \quad\!\!
[e_{i'},[f_{j'},f_{-n,k}]]= \delta_{i'j'}B_{kj'}f_{-n,j'}\;, \label{eifjf0a}\\
[e_{-n+1},f_{-n,i}]=
[f_{-n+1},[f_{-n+1},f_{-n,i}]]=0\;,\label{ffphiis0}\\
[f_{-n,i},f_{-n,j}]=[f_{-n,i''},[f_{-n,j''},f_{-n+1}]]=0\;,
\label{IdealJ2'}
\end{gather}
where $i',j'=-n+2,\ldots,r$ and $i'',j''=-n+3,\ldots,r$.
We let $\tilde S(\fg^{(n)})$ be the Lie superalgebra defined on these sets of generators and relations.

For each node $i$ we get a $\mathbb{Z}$-grading of $\tilde S(\fg^{(n)})$ by letting the generators
$e_i$ and $f_i$ (or $f_{-n,j}$, if $i=-n$) have degree $1$ and $-1$, respectively, and all other generators degree $0$,
since the relations above respect such a grading.
All these $\mathbb{Z}$-gradings can be put together to a $\mathbb{Z}^{n+r}$-grading,
where the subalgebra at degree $(0,0,\ldots,0)$ contains all the $h$ generators, but no $e$ or $f$ generators.
We let $J$ be the maximal ideal of $\tilde S(\fg^{(n)})$ that intersects the degree $(0,0,\ldots,0)$ subalgebra trivially
(obtained as the sum of all ideals with this property)
and set $S(\fg^{(n)})=\tilde S(\fg^{(n)})/J$,
in analogy with the definition of a contragredient Lie superalgebra.
(For a Kac--Moody algebra, or more generally a Borcherds superalgebra,
this step in the construction is equivalent to imposing the Serre relations. For tensor hierarchy algebras, we do no know the full set of relations generating $J$.
It is possible that some of the relations that we impose above are redundant in the sense that
they are contained in this ideal.)

The tensor hierarchy algebra $W(\fg^{(n)})$ is constructed analogously to $S(\fg^{(n)})$, but when we go from $\scr B(\fg^{(n)})$ to $W(\fg^{(n)})$
we add the
$n+r-1$ generators $f_{-n,i}$ without removing the generators $h_{-n}$ and $f_{-n}$.
More precisely, we replace $f_{-n}$ with $n+r$ generators $f_{-n,i}$, where $i=-n,-n+2,-n+3,\ldots,r$. The relations that we impose take 
the same form as for $S(\fg^{(n)})$, and those involving $f_{-n}$ for $\scr B(\fg^{(n)})$ then follows by setting $f_{-n}=f_{-n,-n}$. Thus, in practice,
we rename the generator $f_{-n}$ rather than replacing it. Since $W(\fg^{(n)})$ 
contains both $f_{-n}=f_{-n,-n}$ and $f_{-n,i}$ for $i=-n,-n+2,-n+3,\ldots,r$, it
can be seen as a unification of $S(\fg^{(n)})$ and $\scr B(\fg^{(n)})$,
``usually'' containing them as subalgebras. (We are aware of one case, see Section \ref{grading-fermionic}, where $S(\fg^{(n)})$ is obviously not a subalgebra
of $W(\fg^{(n)})$.
However, in many cases they certainly are, and for all practical purposes that we have encountered so far, both $S(\fg^{(n)})$ and $\scr B(\fg^{(n)})$ can safely be considered as subalgebras of $W(\fg^{(n)})$.)
As we will see, when the Cartan matrix $B_{ij}$ is singular, $\scr B(\fg^{(n)})$ is contained already in $S(\fg^{(n)})$, and a different definition of
$W(\fg^{(n)})$ might be more useful.

\begin{figure}
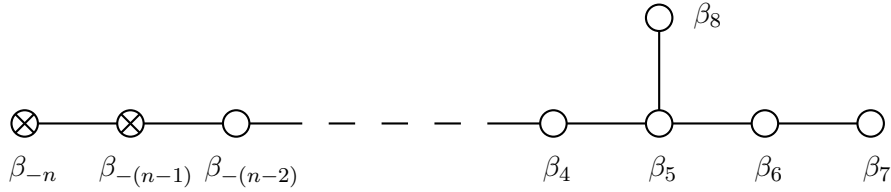

\begin{center}
\EdiagramTwoGreyNodes
\caption{\it The Dynkin diagram for $\BB(\fg^{(n)})$ and $S(\fg^{(n)})$, when $\fg=E_8$, after an odd Weyl reflection, displaying the $\BB(\fg^{(n-1)})$ and $S(\fg^{(n-1)})$ subalgebras.} \label{beta-diagram}
\label{DynkinSETwoGreyNodes}
\end{center}
\end{figure}

From now on, we will often use Greek letters for the fermionic generators and write $\epsilon_{-n}$ instead of $e_{-n}$ and $\phi_{-n,i}$ instead of $f_{-n,i}$.
This allows us to furthermore set
\begin{align}
\epsilon_{-n+1} &= [\epsilon_{-n},e_{-n+1}]\,, &  \phi_{-n+1,i} &= -[\phi_{-n,i},f_{-n+1}]\,, \label{oddweyl}
\end{align}
for $i=-n+3,\ldots,r$.
It then follows that $\tilde S(\fg^{(n-1)})$ is embedded in $\tilde S(\fg^{(n)})$. We will use this embedding occasionally to ``lift''
calculations done for $\tilde S(\fg^{+})$
to $\tilde S(\fg^{++})$.
The relations (\ref{oddweyl}) corresponds to a so called odd Weyl reflection for $\scr B(\fg^{(n)})$, where the set of simple roots
$\alpha_{-n},\ldots, \alpha_r$ are mapped to another set of roots $\beta_{-n},\ldots, \beta_r$, which can also be taken as simple roots, but corresponding to
a different Dynkin diagram \cite{Dobrev:1985qz}. We then have
\begin{align}
\beta_{-n}&=-\alpha_{-n}\,, & \beta_{-n+1}&=\alpha_{-n}+\alpha_{-n+1}\,, & \beta_i &=\alpha_i\,,
\end{align}
for $i=-n+2,\ldots,r$. The two different, but equivalent, Dynkin diagrams for $\BB(\fg^{(n)})$ and $S(\fg^{(n)})$ when $\fg=E_8$ are shown in 
Figures~\ref{alpha-diagram} and \ref{beta-diagram}. The embedding of $S(E_{7+n})$ in $S(E_{8+n})$ corresponds to removal of the node associated to the simple root 
$\beta_{-n}$ in Figure~\ref{beta-diagram}.

As usual, we introduce an inner product in the vector space dual to the Cartan subalgebra of $\fg^{(n-1)}$, such that the Cartan matrix 
$A_{ij}^{(n-1)}$ of $\fg^{(n-1)}$ is given by
\begin{align}
A_{ij}^{(n-1)}=2\frac{(\alpha_i,\alpha_j)}{(\alpha_i,\alpha_i)}\,,
\end{align}
where 
$\alpha_{-n+2},\ldots,\alpha_r$ are the simple roots.
Let $\alpha$ be a root of $\fg^{(n-1)}$, expanded in the basis of simple roots as
\begin{align}
\alpha = a_{-n+2} \alpha_{-n+2} +\cdots+ a_r\alpha_r\,.
\end{align}
We then let
\begin{align}
h_{\alpha^\vee}=\tilde a_{-n+2} h_{-n+2} +\cdots+ \tilde a_rh_r\,,
\end{align}
where $\tilde a_i=a_i\frac{(\alpha_i,\alpha_i)}2$ 
($i=-n+2,\ldots,r$)
be a corresponding element in the Cartan subalgebra of $\fg^{(n-1)}$, 
so that $[h_{\alpha^\vee},e_\beta]=(\alpha,\beta)e_\beta$,
for any  
generator $e_\beta$ corresponding to a root $\beta$.
Set
\begin{align}
\phi_{-n,{\alpha^\vee}}=\tilde a_{-n+2} \phi_{-n,-n+2} +\cdots+ \tilde a_r\phi_{-n,r}\,.
\end{align}
As noted in ref. \cite{Cederwall:2019qnw} (in the simply laced case), and shown in ref. \cite{Carbone:2018xqq} 
(in the case where $\alpha,\beta,\gamma$ are simple roots), we now have
\begin{align}
(\alpha,\beta)[e_\alpha,\phi_{-n,{\gamma^\vee}}]=(\alpha,\gamma)[e_\alpha,\phi_{-n,{\beta^\vee}}]
\label{PhiCartanEq}
\end{align} 
for any roots $\alpha,\beta,\gamma$ of $\fg^{(n-1)}$.
This identity, reflecting that $\langle \phi_{-n,i}\rangle$ provides the Cartan subspace of an adjoint 
$\fg^{(n-1)}$-module, will prove very useful, as well as the following ones:
\begin{align}
[e_{-n+1},[e_{-n+2},\phi_{-n,i}]]=0\,,\label{e-identitet}\\
[f_{-n+1},[e_{-n+2},\phi_{-n,i}]]=0\,\label{f-identitet}.
\end{align}
In the case $i=-n+2$, the first identity (\ref{e-identitet}) follows by
\begin{align}
2[e_{-n+1},e_{-n+2},\phi_{-n,-n+2}]&=[e_{-n+1},e_{-n+2},e_{-n+2},f_{-n+2},\phi_{-n,-n+2}]\nn\\
&=2[e_{-n+2},e_{-n+1},e_{-n+2},f_{-n+2},\phi_{-n,-n+2}]\nn\\
&\quad\,-[e_{-n+2},e_{-n+2},e_{-n+1},f_{-n+2},\phi_{-n,-n+2}]\\
&=4[e_{-n+2},e_{-n+1},\phi_{-n,-n+2}]\nn\\
&\quad\,-[e_{-n+2},e_{-n+2},f_{-n+2},e_{-n+1},\phi_{-n,-n+2}]=0\,,\nn
\end{align} 
and we can then go to the case of a general $i$ by (\ref{PhiCartanEq}). The second identity (\ref{f-identitet}) now follows by acting on (\ref{e-identitet})
with $f_{-n+1}$ twice,
\begin{align}
0&=[f_{-n+1},f_{-n+1},e_{-n+1},e_{-n+2},\phi_{-n,i}]\nn\\
&=-[f_{-n+1},h_{-n+1},e_{-n+2},\phi_{-n,i}]+[f_{-n+1},e_{-n+1},f_{-n+1},e_{-n+2},\phi_{-n,i}]\nn\\
&=-[h_{-n+1},f_{-n+1},e_{-n+2},\phi_{-n,i}]+[e_{-n+1},f_{-n+1},f_{-n+1},e_{-n+2},\phi_{-n,i}]\\
&=2[f_{-n+1},e_{-n+2},\phi_{-n,i}]\,.\nn
\end{align}

A few words on the generators at level $-1$. This is where tensor hierarchy algebras differ from contragredient superalgebras, and it is also the ultimate source of much of the interesting behaviour we will encounter later. 
When $n=0$, \ie, when we consider a tensor hierarchy algebra extension $S(\fg)$ of the finite-dimensional Lie algebra $\fg$, we have the generators 
$\phi_{0,\alpha}$, where $\alpha$ is a root of $\fg^-$, the Lie algebra obtained by deleting
node $1$.
Given the relation (from \eqref{eifjf0a}) $[\epsilon_0,\phi_{0,\alpha}]=h_\alpha$, and the observation that $\epsilon_0$ carries $\fg^-$ weight $0$, it follows that $\phi_{0,\alpha}$ forms the Cartan part of a $\fg^-$ adjoint
 \cite{Carbone:2018xqq}.
 Then, eq. \eqref{PhiCartanEq} becomes natural, it simply states how raising operators act on the Cartan part of the adjoint module.
 The highest state in this adjoint module of course is $[e_{\theta^-},\phi_{0,2}]$, where $\theta^-$ is the highest root of $\fg^-$.
 It is also the highest state in the $\fg$-module at level $-1$, with highest weight
 $\lambda+\theta^-$.
 As we will see, it is this last statement that becomes modified when one considers extensions of infinite-dimensional Kac--Moody algebras, \eg\ $S(\fg^+)$, and leads to the appearance of generators 
where $\epsilon$ and $\phi$ appear simultaneously without cancelling out.

\begin{figure}
\begin{center}
\Ediagram
\caption{\it The Dynkin diagram for $\BB(\fg^{(n)})$ and $S(\fg^{(n)})$ when $\fg=E_8$.} \label{alpha-diagram}
\label{DynkinSE}
\end{center}
\end{figure}

\subsection{Gradings}

As described above, there is a $\mathbb{Z}$-grading of $\tilde S(\fg^{(n)})$ with respect to each node in the Dynkin diagram. These
$\mathbb{Z}$-gradings 
are inherited to $S(\fg^{(n)})$, and can be used to derive the content of the algebra, as we will now describe in some important cases.

\subsubsection{Grading with respect to the fermionic node} \label{grading-fermionic}

The $\mathbb{Z}$-grading with respect to the fermionic node plays a distinguished r\^ole. We call the corresponding degree ``level'' and denote it by $\ell$.
This is a consistent grading, which means that even (bosonic) elements appear at even levels, and odd (fermionic) elements at odd levels.
Since removal of the grey node yields the Dynkin diagram of $\fg^{(n)}$, we have a $\fg^{(n)}$-module $R_\ell$ at each level $\ell$.
The corresponding $\fg^{(n)}$-module $\scr B_\ell$ that appears at level $\ell$ of $\scr B(\fg^{(n)})$ is usually a submodule of $R_\ell$
for $\ell \geq0$. The reason for this is that all elements in $\scr B_\ell$ for $\ell \geq0$ are precisely those that can be formed without 
using the generator $f_{-n}$. Since the remaining generators (except for $h_{-n}$) are included also in $S(\fg^{(n)})$, and the relations among them are the same,
these elements can be formed in $S(\fg^{(n)})$ as well. (However, it might happen that such an element is factored out when going from $\tilde S(\fg^{(n)})$
to $S(\fg^{(n)})$. We are only aware of one case where this happens, which is for a singlet at level $6$ in $S(E_6)$. We have not seen it in the cases
that we focus on in this paper, and therefore we say that $\scr B_\ell$ is ``usually'' a submodule of $R_\ell$).
But there are in general more elements that can be formed at non-negative levels in 
$S(\fg^{(n)})$ than in $\BB(\fg^{(n)})$, since the $\epsilon$ and $\phi$ generators do not need to cancel out in multiple brackets, as they do in the 
contragredient Lie superalgebra $\scr B(\fg^{(n)})$. A first example of such an ``extra'' element is (\ref{J-generator}) that
appears at level $0$ in $S(\fg^+)$
as we will describe below.

\subsubsection{$\gl$-covariant gradings\label{GLGradingSection}}

The $\mathbb{Z}$-grading with respect to the fermionic node can be refined to a $(\mathbb{Z}\times\mathbb{Z})$-grading (or ``double grading'')
with respect to both the fermionic node and one other node. We still call the degree with respect to the fermionic node $\ell$ ``level''.
In the $\gl$-covariant gradings,
the other grading is with respect to a node, which when annihilated leaves the Dynkin diagram of
$A_{d-1}$, extended by a fermionic node at one end, \ie, the Dynkin diagram of
$\BB(A_{d-1})\simeq\sl(d|1)$, $S(A_{d-1})\simeq S(d)$ or $W(A_{d-1})\simeq W(d)$, where $d=r+n$.
This degree $m$ is simply called ``degree''. The second node is node $8$ for $S(E_{8+n})$ (``the exceptional node'') and node $r$ for $S(A_r^{(n)})$ according to the numbering in Figures \ref{E-figur} and \ref{A-figur} (we will call it node $r$ also in the general case).
This results in a double grading, where each element is a $\gl(d)$ module.
In addition, each degree constitutes a module of $W(d)$, and the subalgebra at $m=0$ is $W(d)$ itself.

This can be shown by constructing, as in ref. \cite{Bossard:2017wxl}, a simple local Lie superalgebra $T_{-1} \oplus T_{0} \oplus T_1$,
where $T_0$ is the Lie superalgebra $W(d)$ and $T_\pm$ are suitable $W(d)$-modules, at degree $m=\pm1$,
with a bracket $[T_1,T_{-1}]\subseteq T_0$ defined such that the Jacobi identity
for one element in each of three subspaces is satisfied \cite{Kac77B}. (The brackets $[T_0,T_{\pm1}]\subseteq T_{\pm1}$ are given by the representations that these modules come with.)
With a suitable identification of the generators $\{e_i,f_i,f_{-n,i},h_i\}$ as elements in the three subspaces $T_0,T_{\pm1}$ (corresponding to degree $m=0,\pm1$)
it can then be shown that the bracket relations in this local Lie superalgebra follow from the defining relations for the generators, and conversely.
This means that the minimal $\mathbb{Z}$-graded Lie superalgebra with the local part $T_{-1} \oplus T_{0} \oplus T_1$ is isomorphic to $S(\fg^{(n)})$.

The identification of the generators is such that $e_r$ and $f_r$ belong to $T_1$ and $T_{-1}$, respectively,
whereas all the other generators belong to $T_0$.
These other generators are the generators of $W(A_{d-1})\simeq W(d)$,
up to a redefinition of the $h$ generators (needed since there is an $h_r$ but no $h_{-n}$ at degree $m=0$ in
$S(\fg^{(n)})$, whereas there is an $h_{-n}$ but no $h_r$ in $W(A_{d-1})$). In ref. \cite{Carbone:2018xqq} the identification of the generators was given
in the case $\fg^{(n)}=E_{8+n}$. The $W(d)$-module $T_1$ is then (as in all simply laced cases) the Grassmann algebra $\Lambda(d)$ on $d$ anticommuting generators that
$W(d)$ naturally acts on as its derivation algebra \cite{Kac77B}, but the bracket $[T_0,T_1]$ also contains a trace term (related to the redefinition of the $h$ generators).
The module $T_{-1}$ has a more complicated structure, given in ref.
\cite{Bossard:2017wxl}, but in order to construct the algebra, it is sufficient to know the bracket relation $[f_r, T_1]\subseteq T_0$, given in ref. \cite{Carbone:2018xqq}.
The module $T_{-1}$  will then simply be the $W(d)$-module generated by the single element $f_r$, and the general brackets $[T_{-1},T_1]\subseteq T_0$
can be defined so that the Jacobi identity is automatically satisfied.
Besides giving information about the content of tensor hierarchy algebra, this construction also shows that it exists (non-trivially), which is not obvious from the definition by generators and relations, since there is {\it a priori} the risk that the ideal generated by the relations will be the whole free Lie superalgebra. 

This type of grading was used for $S(E_{11})$ in ref.
\cite{Bossard:2017wxl}, and constructed in detail up to degree $4$.
We will not construct the $W(d)$ representations for higher degree $m$. Instead we will make use of their property, that they consist of sums of ``columns'', where each column is a tensor product
of some $\gl(d)$ module $r_{\ell,m}$ with all forms:
\be
C_{\ell,m}=\bigoplus_{p=0}^dc_{\ell-p,m}=r_{\ell,m}\otimes\bigoplus_{p=0}^d\Lambda^p\;,
\ee
where $\Lambda^p$ is the $p$-form module.
The term $c_{\ell-p,m}=r_{\ell,m}\otimes\Lambda^p$ occurs at level $\ell-p$.
Some further restrictions on $r_{\ell,m}$ are derived in Section
\ref{InvariantFormSection}.

Examples of this type of double grading are given for $S(E_9)$ in Figure \ref{SE9Figure}, for $S(A_1^{++})$ in Figure \ref{SA1++Figure}, and for $S(E_{10})$ in Figure \ref{SE10Figure}. A list of $r_{\ell,m}$ is given for $S(E_{10})$ up to $m=7$ in Table \ref{SE10ColumnsTable} and for $S(A_1^{++})$ up to $m=7$ in Table \ref{SA1++ColumnsTable}.
They are found using the technique described in Section \ref{ExtraModulesSubSection}.

\subsubsection{ ($\fg^-\oplus\gl(n+2))$-covariant grading} \label{gl(n+2)covgrad}

Besides the $\gl$-grading, which is not always present,
there are other interesting and useful double gradings
given by the fermionic node $-n$ and one other node. As the other node, we can take node $1$. 
Removal of this node from the Dynkin diagram of $\fg$ gives the Dynkin diagram of the subalgebra $\fg^-$, and removal of both nodes
from the Dynkin diagram of $S(\fg^{(n)})$ then gives the Dynkin diagram of $A_{n+1} \oplus \fg^-$.
If we now set $d=n+2$ and let $m$ be the degree with respect to node $1$, we have a module of both $\gl(d)$ and $\fg^-$ at each $(\ell,m)$,
and the subalgebra at degree $m=0$ is the direct sum of $W(d)$ and the tensor product $\Lambda(d) \otimes \fg^-$
of the Grassmann algebra $\Lambda(d)$ and $\fg$. 
At degree $m=1$ we have the tensor product of $\Lambda(d)$ and the fundamental module of $\fg^-$ (the $\bf 56$ of $E_7$ in the case of $S(E_{n+8})$).
In the general there will be a column structure like in the $\gl$-covariant grading at any $m$, but now each column is a tensor product
of some $\fg^-$ module with all forms.

Similarly to the $\gl$-covariant construction, we can now construct a local Lie superalgebra $T_{-1}\oplus T_0\oplus T_1$
and show that the minimal $\mathbb{Z}$-graded Lie superalgebra with this local part is isomorphic to $S(\fg^{(n)})$,
and thereby that the tensor hierarchy algebra is non-trivial, also in cases that do not admit a $\gl$-grading.

It might appear more natural to consider a double grading of $S(\fg^{(n)})$ with respect to the
node $0$ rather than node $1$ (together with the fermionic node $-n$). Such a grading would be covariant under $\fg\oplus\gl(n+1))$
instead of $\fg^-\oplus\gl(n+2))$. However, the ``extra'' elements would then appear already in the local part, which would make it more complicated. In particular
the $W(d)$-module at degree $m=1$ would not be irreducible anymore.

\subsubsection{Using gradings to find ``extra'' modules
\label{ExtraModulesSubSection}}

The superalgebra $S(\fg^{(n)})_{\geq0}$ , \ie, the non-negative level subalgebra of $S(\fg^{(n)})$, 
contains a subalgebra $\fg^{(n)}\loplus\BB(\fg^{(n)})_{>0}$ (up to the possibility that elements in $\BB(\fg^{(n)})_{>0}$ are factored out when 
going from $\tilde S(\fg^{(n)})$ to $S(\fg^{(n)})$, described in Section \ref{grading-fermionic}). Everything else we call ``extra modules''.
Using the generators and relations, it is straightforward to show that extra modules are absent for $n=0$ and that only a ``singlet'', $L_1$, arises when $n=1$ (see Section \ref{SAffineSection}).
For $n\geq2$, the analysis in terms of generators and relations becomes complicated. To some extent (to reasonably low levels) we can make use of the tensorial methods of Sections \ref{ExtendByFund} and \ref{SgplusplusSection}, but find it very useful to verify the findings by experimental evidence.

The $\gl$-covariant gradings of Section \ref{GLGradingSection} provide such a laboratory. We start from the $\gl$-covariant decomposition of the Borcherds superalgebras, which we calculate to some level and degree $(\ell,m)$ using the methods of Appendix \ref{BorcherdsAppendix} (to relate the Borcherds modules to certain lowest weight modules) and Appendix \ref{GradingAppendix} (to calculate branchings to $\gl$ modules).
Then, this decomposition can be tested against the prediction that $\gl$-modules at each degree form ``columns'' of the type described in Section \ref{GLGradingSection}, demanded by modules of the $W(d)$ superalgebra at degree $0$, with the other restrictions stated in Section \ref{InvariantFormSection}.
We thus use the interplay between $\fg^{(n)}$ modules at each $\ell$ and $W(r+n)$ modules at each $m$, but without knowing the details of the latter.
Finding a $\gl$ module predicted by a column, but not appearing in the decomposition of the Borcherds superalgebra, signals the appearance of a new $\fg^{(n)}$ module, starting with this $\gl$ module.
This can be done sequentially for increasing $m$.

By computing the $\gl(9)$ decomposition of $S(E_9)$ up to $\ell=22$, $m=7$
(Figure \ref{SE9Figure}), and the $\gl(2)$ decomposition of $S(A_1^+)$ up to $\ell=12$, $m=10$, we test this method against the known result
\cite{Cederwall:2019qnw}
that $L_1$ represents the only extra module (see Section \ref{SAffineSection}).

For the tensor hierarchy algebra extensions of hyperbolic algebras,
we have computed the $\gl(10)$ decomposition of $S(E_{10})$ up to $\ell=22$, $m=7$
(Figure \ref{SE10Figure}), and the $\gl(3)$ decomposition of $S(A_1^{++})$ up to $\ell=12$, $m=10$.
(Figure \ref{SA1++Figure}).
In both cases, we observe perfect agreement with the conjecture that
$R_\ell=\BB_\ell\oplus\BB_{\ell+1}$.
The $\gl$ modules $r_{\ell,m}$ leading the columns for $S(E_{10})$ are listed in Table
\ref{SE10ColumnsTable}, and those for $S(A_1^{++})$ are found in Table~\ref{SA1++ColumnsTable}.

This method will prove valuable when continuing to very extended algebras, \eg\ $E_{11}$.
So far, we only have partial results \cite{CederwallPalmkvistForth}.

\subsection{The invariant bilinear form\label{InvariantFormSection}}

The tensor hierarchy algebra $S(\fg^{(n)})$ is not contragredient, meaning that there is no symmetry, or ``duality'', between positive and negative levels. Nevertheless, it has been shown that it allows for the existence of an invariant non-degenerate bilinear form $\omega$, analogous to the Killing form of Kac--Moody algebras with non-degenerate Cartan matrices. There is thus a ``duality'', as strong as contragredience.
The existence of $\omega$ was one of the features that were more or less built into the original construction of the tensor hierarchy algebras as extensions of exceptional Lie algebras
\cite{Palmkvist:2013vya}.
Since $\omega$ is not (generically) centered around the origin for any grading, it is difficult to derive using the defining generators and relations \cite{Carbone:2018xqq}, and it indeed remains somewhat of a mystery why this very basic property results from these definitions. 
The bilinear form is not present in $W(\fg^{(n)})$. Indeed, it has non-vanishing eigenvalue under
$h_{-n}$. The absence of this generator in $S(\fg^{(n)})$ allows for the existence of an invariant $\omega$ not centered around the origin.
A formalism combining the universality of the generators and relations with the manifestation of the duality would be highly desirable.

The existence of $\omega$ was shown in ref. \cite{Bossard:2017wxl}, using the type of grading described in Section \ref{GLGradingSection} for the specific case of the $E$-series. The proof is straightforward to extend to the more general setting of $S(\fg^{(n)})$, as it only relies of the fact that one member of the series ($n=1$) is the extension of an affine Kac--Moody algebra.
This is also a necessary condition; in a general tensor hierarchy $S(\mathfrak{a})$, there is no such bilinear form $\omega$ (as is evident from the finite-dimensional
$S(A_{n-1})\simeq S(n)$).

The bilinear form $\omega$ pairs elements whose levels add to $1-n$. They thus belong to dual modules of the level $0$ subalgebra.

The bilinear form is invariant in the sense that
\be
\omega([a,b],c)+(-1)^{|a||b|}\omega(b,[a,c])=0\label{OmegaInvariance}
\ee
for all $a,b,c\in S(\fg^{(n)})$.
Interpreting $\omega([a,b],c)$ as (shorthand for) structure constants with the last index lowered and using the graded symmetry of $\omega$, $\omega(a,b)=(-1)^{|a||b|}\omega(b,a)$, eq. 
\eqref{OmegaInvariance} immediately leads to
\be
\omega([a,b],c)=-(-1)^{|b||c|}\omega([a,c],b)\;,
\ee
stating the total graded antisymmetry of the structure constants.
For $n$ even, $\omega$ is fermionic, in the sense that it pairs bosonic with fermionic generators, while for $n$ odd it pairs generators of the same statistics.

In the $E$ series, the $\gl$ grading for $S(E_{8+n})$ is dual under reflection in 
$(\ell,m)=(\frac{1-n}2,\frac32)$ (also for non-positive $n$). 
From this, it easily follows, since the bottom of the lowest column will be dual to an element at $m=0$, that the top of the lowest columns for given $m\geq4$ is at $\ell=9$, and the corresponding $r_{\ell,m}$ are $r_{9,m}=(E_{8+n})_{m-3}=S(E_{8+n})_{0,m-3}$, degree $m-3$ in the $\gl(8+n)$-grading of $E_{8+n}$.

In the $A$ series, labelled by $r$ and $n$, the center of duality for $S(A_r^{(n)})$ in the $\gl$ grading is at 
$(\ell,m)=(\frac{1-n}2,\frac12)$. The top of the lowest columns for $m\geq2$ appear at 
$\ell=r+1$, and the corresponding $r_{\ell,m}$ are 
$r_{r+1,m}=(A_r^{(n)})_{m-1}=S(A_r^{(n)})_{0,m-1}$.
For all $r$ and $n$, the subalgebra (superalgebra for odd $r$) at the line $\ell=rm$ is freely generated.

In the $(\fg^- \oplus \gl(n+2))$-covariant grading, the center of duality is always at $(\ell,m)=(\frac{1-n}2,1)$. Thus the subspace at $m=1$,
which can be considered as the tensor product of $\Lambda(d)$ with the fundamental module of $\fg^-$ (se Section \ref{gl(n+2)covgrad}),
is paired to itself, where it is given by a combination of the volume form on $\Lambda(d)$ and the symplectic form on
the fundamental module of $\fg^-$.

\section{The tensor hierarchy extension of an affine algebra
\label{SAffineSection}}

We now restrict to the case $n=1$ and consider $S(\fg^+)$.

As mentioned at the end of Section \ref{GenRelSection}, the generators $\phi_{-1,\alpha}$ reside as the Cartan part of an adjoint $\fg$ module at level $-1$. This statement still holds true. However,
$[e_\theta,\phi_{-1,1}]$ is not the highest state in this module. One is now allowed to act with $e_0$ to obtain $[e_0,e_\theta,\phi_{-1,1}]$. This is easily shown to be a highest weight state in a $\fg^+$ module, an anti-fundamental.

Then, one is free to act with $\epsilon_{-1}$ and obtain
$[\epsilon_{-1},e_0,e_\theta,\phi_{-1,1}]$ (at level $0$).
There is no way to use the identities of Section \ref{GenRelSection} to relate this expression to one without $\epsilon$ and $\phi$.
For the first time, there exists a generator, formed from the simple raising and lowering operators, where a raising and a lowering operator associated with the fermionic node, \ie, an $\epsilon$ and a $\phi$, do not cancel out. This generator is
\begin{align} \label{J-generator}
J=[\epsilon_{-1},e_0,e_\theta,\phi_{-1,1}]=[\phi_{-1,1},e_\theta,e_0,\epsilon_{-1}]\;.
\end{align}
For tensor hierarchy algebras $S(\fg)$, with $\fg$ finite-dimensional, this can not happen
\cite{Carbone:2018xqq}.
There are no relations that can be used to simplify $J$. It arises due to the light-likeness of the $\fg^+$ root $\alpha_0+\theta$.

We will now investigate the properties of the generator $J$.
For $i\geq1$, \ie, for brackets with raising and lowering operators within the finite-dimensional algebra $\fg$,  we have
\begin{align}
  [e_i,J]&=[\epsilon_{-1},e_i,e_0,e_\theta,\phi_{-1,1}]
  =[\epsilon_{-1},[e_i,e_0],e_\theta,\phi_{-1,1}]=0\;,\nn\\
  [f_i,J]&=[\epsilon_{-1},e_0,e_\theta,f_i,\phi_{-1,1}]=0\;.\label{Jeifi}
\end{align}
The second of these equations follows directly from $[e_0,e_\theta,f_i,\phi_{-1,1}]=0$, which
holds (see Section \ref{GenRelSection})
since $[e_\theta,f_i,\phi_{-1,1}]$ is lower than $[e_\theta,\phi_{-1,1}]$ in
the $\fg$ adjoint. The first then follows since
$[h_i,e_0,e_\theta,f_{-1,1}]=0$, and thus
$[f_i,e_i,e_0,e_\theta,\phi_{-1,1}]=[e_i,f_i,e_0,e_\theta,\phi_{-1,1}]=0$ and
\begin{align}
  0&=[e_i,f_i,e_i,e_0,e_\theta,\phi_{-1,1}]\nn\\ 
  &=[h_i,e_i,e_0,e_\theta,\phi_{-1,1}]+[f_i,[e_i,e_i,e_0],e_\theta,\phi_{-1,1}]\\
  &=2[e_i,e_0,e_\theta,\phi_{-1,1}]\;.\nn
\end{align}

We then consider the action of $e_0$ and $f_0$.
Lowering gives
\begin{align}
  [f_0,J]&=[\epsilon_{-1},f_0,e_0,e_\theta,\phi_{-1,1}]\nn\\
  &=-[\epsilon_{-1},h_0,e_\theta,\phi_{-1,1}]+[\epsilon_{-1},e_0,e_\theta,f_0,\phi_{-1,1}]
  \label{Jf0}\\
  &=[e_\theta,\epsilon_{-1},\phi_{-1,1}]=[e_\theta,h_1]=-e_\theta\;.\nn
\end{align}
Using $[e_0,e_0,J]=0$ and $[h_0,J]=0$ we then get
\be
[e_0,J]=-{1\over2}[e_0,e_0,e_\theta]\;.\label{Je0}
\ee
The generator $J$ can thus be considered a 
{\it weakly lowest weight state},
by which we understand that acting with lowering operators only gives generators in another module (in the present case, the adjoint).
$J$ is indeed also a weakly highest weight state, a singlet.

In a mode expansion, we can identify
\begin{align}
e_0&=e_{-\theta,1}\;,\nn\\
f_0&=e_{\theta,-1}\;,
\end{align}
where the last index is the mode number, since
\be
[e_0,f_0]=h_0=h_{\delta-\theta}=\KK-h_\theta\;.
\ee
We then have $[e_0,e_\theta]=[e_{-\theta,1},e_\theta]=-h_{\theta,1}$,
and $[e_0,[e_0,e_\theta]]=-[e_{-\theta,1},h_{\theta,1}]=-2e_{-\theta,2}$.
Eqs. \eqref{Jeifi}, \eqref{Jf0} and \eqref{Je0} then state that
\be
[J,T^A_m]=-mT^A_{m+1}\;,\quad m=0,\pm1\;,
\ee
implying that $J$ acts as $L_1$ on all loop generators.

At level $1$, the lowest weight state of the adjoint fundamental is
$\epsilon_{-1}$. Its bracket with $J$ is
\begin{align}
[\epsilon_{-1},J]=[\epsilon_{-1},\epsilon_{-1},e_0,e_\theta,\phi_{-1,1}]=0\;.
\end{align}
Note that $L_1$ (at $k=-1$) annihilates this state. Together with the action on the generators, this shows that $J$ acts as $L_1$ at $k=-1$ on the fundamental module $f$ at level $1$ in the 
tensor hierarchy algebra, which coincides with level $1$ in the Borcherds superalgebra.

We know that there is a non-degenerate quadratic form, pairing level $\ell$ with level $-\ell$. 
The unique restriction to level $0$ is given by
\be
\omega_0=-\KK\otimes L_1-L_1\otimes \KK
+\sum_{m\in\ZZ}\eta_{AB}T^A_m\otimes T^B_{1-m}\;,\label{Omega0}
\ee
where $\eta_{AB}$ is the Killing metric of $\fg$.
At level $-1$ we must necessarily have an anti-fundamental $\bar f$. Given the ``shift'' in 
eq. \eqref{Omega0}, this anti-fundamental should be considered as shifted by $1$ in mode number, \ie, by $-1$ in the eigenvalue of $\dd$ (when it is included, see below).

\begin{figure}
  \begin{center}
\includegraphics[scale=.7]{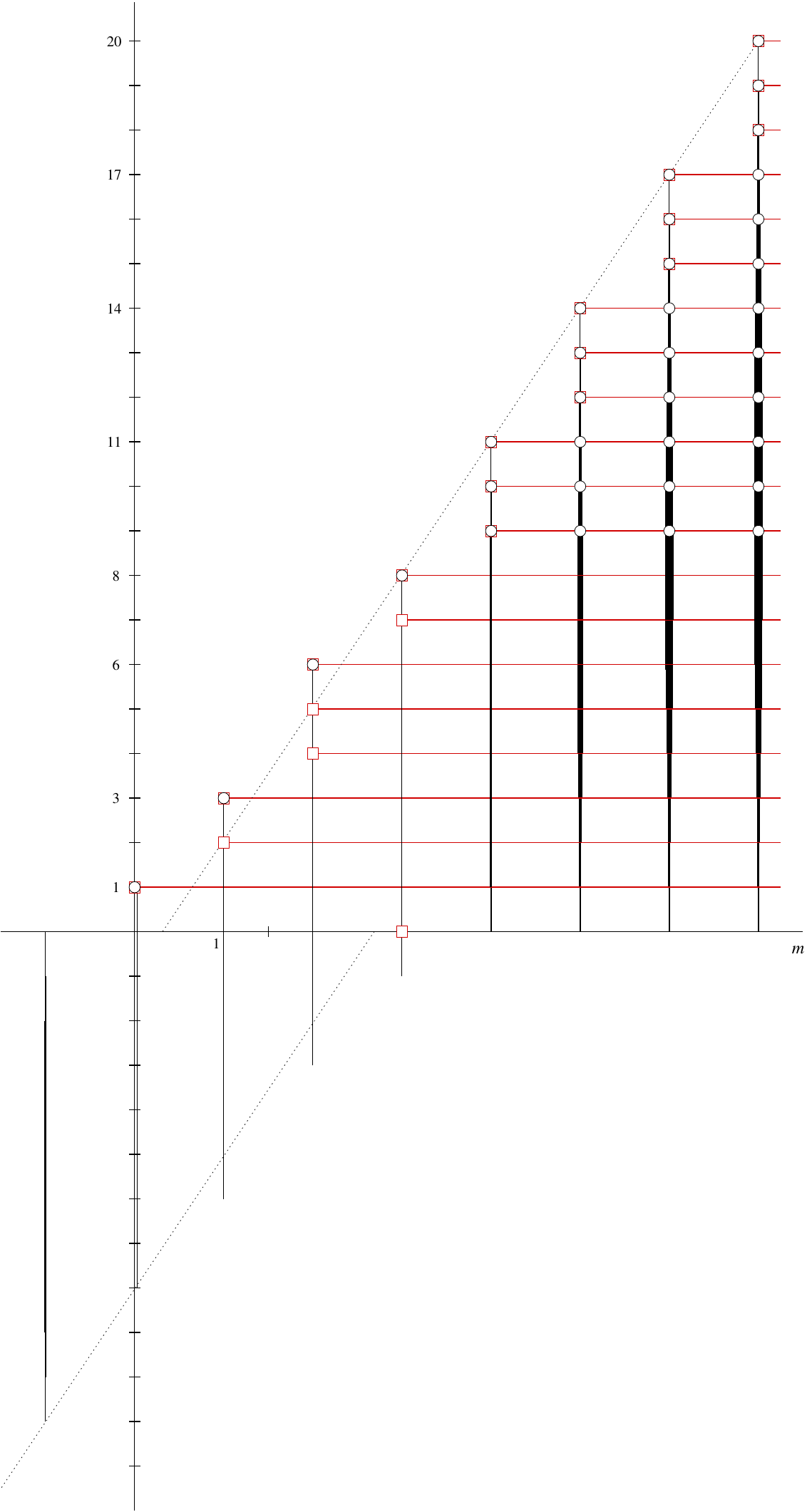}
\caption{\it The double
grading of $S(E_{9})$, with 
respect to the fermionic and the exceptional node. Each grade contains a
$\gl(9)$ module. These modules organise into ``columns'', containing
tensor products of some $\gl(9)$ module with all forms. The collection
of such columns at given degree $m$ with respect to the exceptional node forms a
module of the tensor hierarchy algebra $W(9)$ at degree $0$. Tops of
columns are marked with a black circle. Horizontally, the $\gl(9)$
modules are organised in $E_{9}$ modules. The lowest states in
$R_\ell$  at level $\ell$, and $L_1$, are denoted with red squares. The only extra module is the singlet $L_1$ at $\ell=0$, $m=3$.}    
\label{SE9Figure}
\end{center}
\end{figure}

Using a $\gl$-covariant double grading as described in Section \ref{GLGradingSection}, and the calculational methods of Appendices \ref{BorcherdsAppendix} and \ref{GradingAppendix}, it is verified that no extra modules except $J$ arise at positive levels in $S(A_1^+)$ and $S(E_9)$, as compared to the corresponding Borcherds superalgebras. This is depicted in Figure \ref{SE9Figure}.

The tensor hierarchy algebra $S(\fg^+)$ is very similar to $\BB(\fg^+)$. Instead of the unshifted quadratic form, containing $\dd$ (effectively, $L_0$), the shifted bilinear form with $L_1$ is used.  
Note that $\dd$ is not present in $S(\fg^+)$. The Cartan matrix of the super-extension is non-degenerate without its inclusion.

The superalgebra $W(\fg^+)$ contains both $S(\fg^+)$ and $\BB(\fg^+)$ as subalgebras.
Positive levels coincide in all three superalgebras.
In addition to the content in $S(\fg^+)$ we find $\dd$ at level $0$, an unshifted anti-fundamental at level $-1$, and so on. $W(\fg^+)$ is included as the subalgebra at $q=0$ of $S(\fg^{++})$ in Table
\ref{DoubleGradingTable}.

\section{Extending an over-extended algebra with its fundamental module
\label{ExtendByFund}}

\subsection{Affine subalgebras and Virasoro algebras}
\begin{figure}
\begin{center}
\includegraphics[scale=.8]{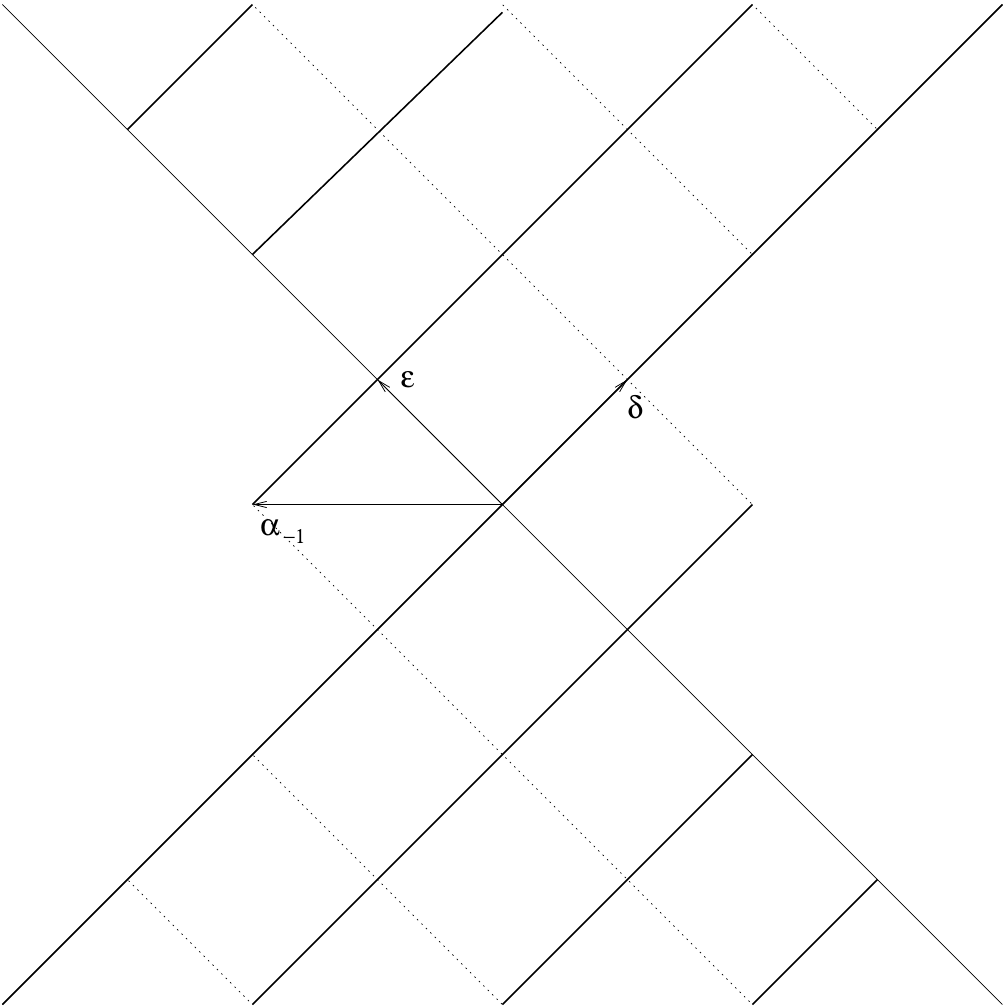}
\caption{\it The branching of the over-extended algebra into
modules of two affine algebras, indicated by the solid and dotted
diagonal lines. Each 
point in the diagram contains a finite-dimensional, finitely reducible
$\fg$-module. 
\label{HyperbolicFigure}}    
\end{center}
\end{figure}

Consider
the extension to a over-extended Kac--Moody (KM) algebra $\fg^{++}$ of an
(untwisted) affine KM algebra $\fg^+$, which in turn is the central extension
of a loop algebra over a finite-dimensional Lie algebra $\fg$.
Number the simple roots of $\fg$ as $\alpha_i$‚ $i=1,\ldots,r$.
The affine root is $\alpha_0=\delta-\theta$, where $(\delta,\delta)=0$ and
$\theta$ is the highest root for $\fg$. The hyperbolic root is
$\alpha_{-1}=\epsilon-\delta$, where $(\epsilon,\epsilon)=0$,
$(\epsilon,\delta)=-1$.
Thus, $\epsilon$ and $\delta$ are vectors on the forward light-c\^one.
We visualise it with the projection on the $\delta\epsilon$-plane as in Figure \ref{HyperbolicFigure}.

The affine algebra is
\be
[T^A_m,T^B_n]=f^{AB}{}_CT^C_{m+n}+\eta^{AB}m\delta_{m+n,0}\KK\;,
\ee
where the central element $\KK$ is identified as the Cartan element $\KK=h_\delta$. The
eigenvalues of $\KK$, the affine level, are denoted $k$.

To the affine algebra is associated a Virasoro algebra, obtained by
the Sugawara construction \cite{Sugawara:1967rw}, which at affine level $k$ is
\be
L_m={1\over2(k+g^\vee)}\sum\limits_{n\in\ZZ}\eta_{AB}:T^A_nT^B_{m-n}:\;,
\ee
$g^\vee$ being the dual Coxeter number.

The embedding of the affine algebra in the over-extended one is not
unique. Any Weyl transformation of $\delta$ maps it to a primitive
light-like root, and the affine algebra to another affine algebra. The
simplest such transformation is a Weyl reflection in the hyperplane
orthogonal to $\alpha_{-1}$, which is a reflection in the vertical axis
through the origin in Figure \ref{HyperbolicFigure}, and exchanges affine mode number and
negative affine level. 
The new affine generators are, for $m\geq0$:
\begin{align}
\tilde T^A_m&=\frac1{m!}(\ad e_{-1})^mT^A_m\;,\nn\\
\tilde T^A_{-m}&=\frac{(-1)^m}{m!}(\ad f_{-1})^mT^A_{-m}\;,
\end{align}
and they fulfil
\be
[\tilde T^A_m,\tilde T^B_n]=f^{AB}{}_C\tilde T^C_{m+n}+\eta^{AB}m\delta_{m+n,0}h_\epsilon\;.
\ee
$T^A_m$ and $\tilde T^A_m$ are, for all values of $A$, lowest and
highest states in an $(m+1)$-dimensional module of the $\sl(2)$ associated to the
hyperbolic root.

The question about where and how to find the infinite number of Virasoro algebras associated to affine subalgebras of a hyperbolic algebra has been asked before
\cite{Kleinschmidt:2005bq}.
We will take a modest approach and limit the question to how the extension of an affine algebra by $L_1$, which we saw arises in the tensor hierarchy extension of an affine algebra, generalises to the over-extended case.  

The level $0$ subalgebra of $S(\fg^{++})$ will contain a subalgebra $\langle L_1\rangle\loplus\fg^+$.
The $\fg^{++}$ weight of $L_1$ is $-\lambda=\delta$. Clearly, it can not be a ''scalar'' under $\fg^{++}$. In particular, it is not consistent to set $[e_{-1},L_1]$ to $0$ if $[f_{-1},L_1]=0$. Instead, we are led to a situation where $L_1$ is the lowest weight state of a fundamental module $F=R(-\lambda)$.
As will be shown in Section \ref{SgplusplusSection} from the generators and relations of Section \ref{THAGeneralSection}, the action of $L_1$ in the affine (lowest weight) fundamental
$f$ at $k=-1$ in $\fg^{++}$ and in the shifted fundamental $f[1]$ at $k=-1$ in $F$ are
\begin{align}
  [L_1,T_\mu]&=-(\ell_1)_\mu{}^\nu T_\nu-J_\mu\;,\nn\\
  [L_1,J_\mu]&=-(\ell_1)_\mu{}^\nu J_\nu\;,\label{LOneAction}
\end{align}
where $\mu,\nu,\ldots$ are fundamental affine indices and $\ell_1$ is the $k=1$ representation matrix for $L_1$.
This Jordan cell structure is a starting point for the full Lie algebra structure on
$\AA=\fg^{++}\oplus F$.

\begin{figure}
\begin{center}
\includegraphics[scale=.9]{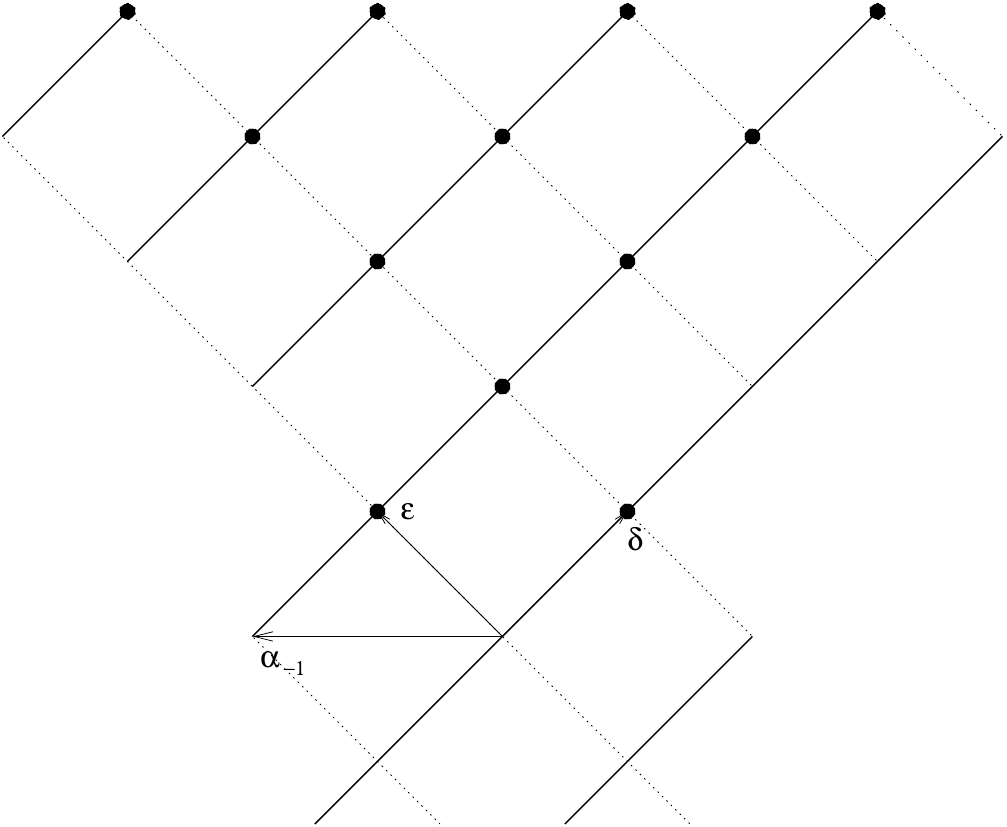}
\caption{\it The projection of weights (black dots) in the over-extended
  lowest weight fundamental module to the $\delta\epsilon$-plane.}    
\end{center}
\end{figure}

\subsection{Transforming derivations\label{TranfDerSection}}
In the case of extending the affine algebra $\fg^+$ by $L_1$, the latter acts as an outer derivation, and one forms the semidirect sum $\langle L_1\rangle\loplus\fg^+$. If we want to extend $\fg^{++}$ in a similar way, it is seen already in eq. \ref{LOneAction} that the fundamental not only transforms the algebra into itself, but also transforms under the algebra. The structure is the following. Let $\{T_\alpha\}$ be generators of a Lie algebra $A$ (in our case $A=\fg^{++}$), and $\{J_M\}$ a basis for some module (the fundamental $F$). 
Then we extend $A$ by $F$ with the brackets
\begin{align}
  [T_\alpha,T_\beta]&=f_{\alpha\beta}{}^\gamma T_\gamma\;,\label{TTBracket}\\
  [T_\alpha,J_M]&=-t_{\alpha M}{}^NJ_N+u_{M\alpha}{}^\beta T_\beta\;,\label{TJBracket}\\
  [J_M,J_N]&=g_{MN}{}^PJ_P\;.\label{JJBracket}\
\end{align}
The Jacobi identities demand that the structure constants $f,t,u,g$ satisfy
\begin{align}
0&=  f_{[\alpha\beta}{}^\epsilon
  f_{\gamma]\epsilon}{}^\delta;,\label{ffid}\\
0&=  [t_\alpha,t_\beta]_M{}^N-f_{\alpha\beta}{}^\gamma t_{\gamma
  M}{}^N\;,\label{ttid}\\
0&=  f_{\alpha\beta}{}^\delta u_{M\delta}{}^\gamma
    +2u_{M[\alpha}{}^\delta f_{\beta]\delta}{}^\gamma
    +2u_{N[\alpha}{}^\gamma t_{\beta]M}{}^N\;,\label{fuid}\\
0&=  g_{MN}{}^Q t_{\alpha Q}{}^P
    +2t_{\alpha[M}{}^Q g_{N]Q}{}^P
    +2t_{\beta[M}{}^P u_{N]\alpha}{}^\beta\;,\label{gtid}\\
0&=  [u_M,u_N]_\alpha{}^\beta-g_{MN}{}^Pu_{P\alpha}{}^\beta
    \;,\label{uuid}\\
0&=  g_{[MN}{}^Rg_{P]R}{}^Q\;.\label{ggid}
\end{align}
We have already made the assumption that $[J,J]$ does not contain $T$. This, and the absence of other modules, \ie, the consistency of eqs.  \eqref{TTBracket}--\eqref{JJBracket} and \eqref{ffid}--\eqref{ggid} for the case at hand, will be shown shortly in Section \ref{ClosureSection}.
We can call $J_M$ ''transforming derivations'', since removing the first term in eq. \eqref{TJBracket} gives the structure of a Lie algebra extended by an algebra of derivations.

The extension of an affine algebra---extended by $\dd$, so the Killing
metric becomes non-singular---by $L_1$ is already an example of this
structure, although somewhat degenerate---there, $L_1$ is a ``shifted
scalar'', which transforms under the Cartan generator $\dd$,
and transforms $T^A_m$, 

The brackets and identities look symmetric, in the sense that $f$ and $g$ are structure constants for two Lie algebras (eqs. \eqref{ffid},\eqref{ggid}), and the algebras form modules with respect to one another with representation matrices $t$ and $u$ (eqs. \eqref{ttid},\eqref{uuid}). Eq. \eqref{fuid} is a modification of the derivation property, and eq. 
\eqref{gtid} its ''mirror''.
The actual situation is however far from symmetric. $\fg^{++}$, the Lie algebra we extend, is simple, while the algebra of its fundamental, as we will see, is far from simple.

In order to claim that one has found a non-trivial structure of this kind, it must of course be non-trivial, in the sense that it can not be brought to the trivial form of a semi-direct sum $A\loplus F$ by redefinition of the generators. This is already clear from the affine subalgebra and $L_1$.

\subsection{Affine decomposition}
The brackets and identities are those given in Section \ref{TranfDerSection}.
We will use the decomposition into affine modules as a tool both for exhibiting the structure of the algebra and to derive some of its properties.

Expansion in negative affine level $-k$ is the same as level expansion with respect to node $-1$ of the over-extended Kac--Moody algebra $\fg^{++}$, since $h_\lambda=-\KK$ counts the number of $e_{-1}$'s.
The local Lie algebra between generators at $k=0,\pm1$ follows straightforwardly from the generators and relations.
We want to continue to $k=-2$ and check that the null states there are consistent 
without any generators beyond the adjoint and $F$.
Let $\mu,\nu\ldots$ be fundamental $\fg^+$ indices.
At $k=0$ the generators in $A$ are $T^A_m$, $\KK$ and $\dd$, and in
$F$ $L_1$. At $k=1$ there is $\bar T^\mu$ in $A$, and at $k=-1$
$T_\mu$ in $A$ and $J_\mu$ in $F$.
Note that $J_\mu$ is a shifted fundamental compared to $T_\mu$, which is reflected in their transformations under $\dd$ in eq. \eqref{dTdJEq}.
We then have
\begin{align}
[T^A_m,T^B_n]&=f^{AB}{}_CT^C_{m+n}+m\delta_{m+n,0}\eta^{AB}\KK\;,\nn\\
[\KK,T^A_m]&=0\;,
\nn\\
[\dd,T^A_m]&=-mT^A_m\;,\\
[L_1,T^A_m]&=-mT^A_{m+1}\;,\nn
\end{align}
\begin{align}
[\dd,L_1]&=-L_1\;,& [\KK,\dd]=[\KK,L_1]&=0
\end{align}
\begin{align}
[T^A_m,T_\mu]&=-(t^A_m)_\mu{}^\nu T_\nu\;,
            &[T^A_m,J_\mu]&=-(t^A_m)_\mu{}^\nu J_\nu\;,\nn\\
[\KK,T_\mu]&=-T_\mu\;,
&[\KK,J_\mu]&=-J_\mu\;,\nn\\
[\dd,T_\mu]&=-(\ell_0-1)_\mu{}^\nu T_\nu\;,
            &[\dd,J_\mu]&=-(\ell_0)_\mu{}^\nu J_\nu\;,\label{dTdJEq}\\
[L_1,T_\mu]&=-(\ell_1)_\mu{}^\nu T_\nu-J_\mu\;,
&[L_1,J_\mu]&=-(\ell_1)_\mu{}^\nu J_\nu\;,\nn
\end{align}
\begin{align}
[T^A_m,\bar T^\mu]&=(t^A_m)_\nu{}^\mu \bar T^\nu\;,\nn\\
[\KK,\bar T^\mu]&=\bar T^\mu\;,\nn\\
[\dd,\bar T^\mu]&=(\ell_0-1)_\nu{}^\mu\bar T^\nu\;,\\
[L_1,\bar T^\mu]&=(\ell_1)_\nu{}^\mu \bar T^\nu\;,\nn
\end{align}
\begin{align}
[T_\mu,\bar T^\nu]&=\delta_\mu^\nu\dd+(\ell_0-1)_\mu{}^\nu\KK
      -\sum\limits_{m\in\ZZ}\eta_{AB}(t^A_m)_\mu{}^\nu T^B_{-m}\;,\nn\\
[J_\mu,\bar T^\nu]&=-\delta_\mu^\nu L_1-(\ell_1)_\mu{}^\nu\KK
      +\sum\limits_{m\in\ZZ}\eta_{AB}(t^A_m)_\mu{}^\nu T^B_{1-m}\;.
\end{align}

These brackets give
\begin{align}
[[T_\mu,T_\nu],\bar T^\kappa]
   &=-2(\Cn0_{[\mu\nu]}{}^{\kappa\lambda}-2\delta_{\mu\nu}^{\kappa\lambda})T_\lambda
   \;,\nn\\
[[T_\mu,J_\nu],\bar T^\kappa]
   &=\Cn1_{\nu\mu}{}^{\kappa\lambda}T_\lambda
   -(\Cn0_{\mu\nu}{}^{\kappa\lambda}-2\delta_{\mu\nu}^{\kappa\lambda})J_\lambda
   \;,\label{TJbarTEqn}\\
[[J_\mu,J_\nu],\bar T^\kappa]
   &=-2\Cn1_{[\mu\nu]}{}^{\kappa\lambda}J_\lambda\;.\nn
\end{align}
(The term $2\delta_{\mu\nu}^{\kappa\lambda}J_\lambda$ is a combination
of terms from $[T_\mu,\bar T^\kappa]$ and $[J_\nu,\bar T^\kappa]$.)
Here,
\be
\Cn m=1\otimes \ell_m+\ell_m\otimes 1
-\sum\limits_{n\in\ZZ}\eta_{AB}t^A_nt^B_{m-n}
\ee
is the matrix form (representation matrix) of $g^\vee+2$
times the coset Virasoro generator
$1\otimes L^{(1)}_m+L_m^{(1)}\otimes1-L^{(2)}_m$, the superscript indicating $k$ in the Sugawara construction. 

In the tensor product of two affine fundamentals, we let
$\vee^2f=s_2^\compl\oplus s_2$ and $\wedge^2f=a_2^\compl\oplus a_2$, where 
$s_2^\compl$ and $a_2^\compl$ are the lowest (leading) symmetric and antisymmetric modules.
Recall that $s_2^\compl$ is annihilated
by $\Cn0$ and $\Cn1$, while
$a_2^\compl$ is annihilated by $\Cn0-2$.
It follows from the first equation in \eqref{TJbarTEqn}
that the leading module in $[T,T]$ is an ideal.
In the second equation, 
both the $T$ and $J$ terms annihilate
$s_2^\compl$, but only the $J$ term annihilates $a_2^\compl$. The
leading symmetric module is an ideal. In the third equation, there is {\it a priori}
no ideal.

The modules in the three brackets are thus
\begin{align}
[T,T]\;&:\quad \wedge^2 f\ominus a_2^\compl=a_2\;,\nn\\
[T,J]\;&:\quad f\otimes f[1]\ominus s_2^\compl[1]
    =s_2[1]\oplus a_2^\compl[1]\oplus a_2[1]\;,\\
[J,J]\;&:\quad \wedge^2 f[1]=a_2^\compl[2]\oplus a_2[2]\;,\nn
\end{align}
where the shifts in brackets are relative to $T$ and its tensor
products.

On the other hand, the modules in $T$ and $J$ at $k=-2$
are $a_2$ and $a_2[1]\oplus s_2[1]$, respectively.
If $T_{\mu\nu}\in a_2$, then $(\Cn1 T)_{\mu\nu}\in a_2^\compl[1]\oplus
a_2[1]$, and similarly $J_{[\mu\nu]}\in a_2[1]$ gives
$(\Cn1 J)_{[\mu\nu]}\in a_2^\compl[2]\oplus a_2[2]$.
Now, the module $a_2^\compl[2]\oplus a_2[2]$ is smaller than $a_2[1]$.
The coset Virasoro character contained in the antisymmetric tensor product of two $k=-1$ modules has non-decreasing coefficients (we assume this, but suspect it is a known fact). For $E_9$ \cite{Bossard:2017aae}, it is
\begin{align}
\chi(q)=\frac{\phi(q^2)}{\phi(q)}
&=1+q+q^2+2q^3+2q^4+3q^5+4q^6+5q^7+6q^8\nn\\
&\quad+8q^9+9q^{10}+12q^{11} +15q^{12}+\ldots\;,
\end{align}
where $\phi(q)=\prod_{n=1}^\infty(1-q^n)$.
The statement that $(a_2^\compl\oplus a_2)[1]\subset a_2$ amounts to the statement that
$\chi(q)-1-q\chi(q)$ has no negative coefficients. We assume this to hold in general. It ensures that anything produced by $[J_\mu,J_\nu]$ can be interpreted as part of $F$ at $k=-2$. 
This is a heuristic argument helping us to understand what to expect at $k=-2$, and does not play a r\^ole in  the arguments of Section \ref{ClosureSection}.

The brackets \eqref{TJbarTEqn} are reproduced by
\begin{align}
[T_\mu,T_\nu]&=2T_{\mu\nu}\;,\nn\\
[T_\mu,J_\nu]&=
      -(\Cn1(\Cn0-2)^{-1})_{\mu\nu}{}^{\kappa\lambda}T_{\kappa\lambda}
      +J_{\mu\nu}\;,\label{OneOneTwo}\\
[J_\mu,J_\nu]&=-2(\Cn1(\Cn0-2)^{-1})_{[\mu\nu]}{}^{\kappa\lambda}J_{\kappa\lambda}\nn
\end{align}
(also defining $T_{\mu\nu}$ and $J_{\mu\nu}$),
and
\begin{align}
[T_{\mu\nu},\bar T^\kappa]
      &=-(\Cn0-2)_{[\mu\nu]}{}^{\kappa\lambda}T_\lambda\;,\nn\\
[J_{\mu\nu},\bar T^\kappa]
      &=\Cn1_{(\mu\nu)}{}^{\kappa\lambda}T_\lambda
      -(\Cn0-1+\sigma)_{\mu\nu}{}^{\kappa\lambda}J_\lambda
      \;,
\label{JTwoTTwoTbar}
\end{align}
where
$\sigma_{\mu\nu}{}^{\kappa\lambda}=\delta_\mu^\lambda\delta_\nu^\kappa$
is the permutation operator, taking the eigenvalue $1$ on the symmetric
part  and $-1$ on the antisymmetric part of a tensor.
Notice that the right hand sides of eq. \eqref{JTwoTTwoTbar} annihilate the
ideals, $a_2^\compl$ and $a_2^\compl\oplus s_2^\compl$, respectively,
and that the inverse of $\Cn0-2$ is well defined on the antisymmetric
modules it acts on (when needed, we define the inverse to be $0$ on the kernel).
It is also interesting to note that there is no $T$ term in $[J,J]$.
This holds to all affine levels (values of $k$), which is demonstrated in Section
\ref{ClosureSection}.
The $T$ term in $[T,J]$ and the $J$ term in $[J,J]$ in eq. \eqref{JTwoTTwoTbar} are shifted by one unit of mode number by $\Cn1$. In view of the discussion on the Virasoro character, they therefore do not contain all of $T_{\mu\nu}$ and $J_{[\mu\nu]}$, respectively. 

The reason for choosing this normalisation of the generators at $k=-2$
is to minimise the occurrence of the coset Virasoro generators.
In this normalisation, the (inverse) Killing form for the over-extended algebra is encoded in
\begin{align}
K&=-\dd\otimes\KK-\KK\otimes\dd
+\sum_{m\in\ZZ}\eta_{AB}T^A_m\otimes T^B_{-m}\nn\\
&+T_\mu\otimes\bar T^\mu+\bar T^\mu\otimes T_\mu
+((\Cn0-2)^{-1})_{\mu\nu}{}^{\kappa\lambda}
\left(T_{\kappa\lambda}\otimes\bar T^{\mu\nu}
       +\bar T^{\mu\nu}\otimes T_{\kappa\lambda}\right)
+\ldots\;.
\end{align}

The brackets between $L_1$ and generators at $k=-2$ are deduced
from eq. \eqref{OneOneTwo}. One obtains
\begin{align}
[L_1,T_{\mu\nu}]&=-((\Cn0-2)(1\otimes\ell_1+\ell_1\otimes1)
(\Cn0-2)^{-1})_{\mu\nu}{}^{\kappa\lambda}T_{\kappa\lambda}
-J_{[\mu\nu]}\;,\nn\\
[L_1,J_{\mu\nu}]
&=-(1\otimes\ell_1+\ell_1\otimes1)_{(\mu\nu)}{}^{\kappa\lambda}J_{\kappa\lambda}\label{LDeltaTwo}\nn\\
&\qquad-((\Cn0-2)(1\otimes\ell_1+\ell_1\otimes1)
      (\Cn0-2)^{-1})_{[\mu\nu]}{}^{\kappa\lambda}J_{\kappa\lambda}\\
&=-((\Cn0-1+\sigma)^{{1-\sigma\over2}}(1\otimes\ell_1+\ell_1\otimes1)
      (\Cn0-1+\sigma)^{-{1-\sigma\over2}})_{\mu\nu}{}^{\kappa\lambda}J_{\kappa\lambda}
\nn\;
\end{align}
(the expressions with exponents containing the permutation operator $\sigma$ should simply be read by inserting its $\pm1$ eigenvalues when acting on the symmetric/antisymmetric parts).
Verification of the consistency of eq. \eqref{LDeltaTwo} with
eq. \eqref{OneOneTwo} and the brackets
\be
[L_1,\left(
\begin{matrix}
T_\mu\\ J_\mu
\end{matrix}
\right)]
=\left(
\begin{matrix}
-(\ell_1)_\mu{}^\nu&-\delta_\mu^\nu\\
0&-(\ell_1)_\mu{}^\nu
\end{matrix}
\right)
\left(
\begin{matrix}
T_\nu\\ J_\nu
\end{matrix}
\right)
\ee
relies on the commutators
\begin{align}
[1\otimes\ell_1+\ell_1\otimes1,\Cn0\,]&=\Cn1\;,\nn\\
[1\otimes\ell_1+\ell_1\otimes1,\Cn1\,]&=0\;.
\label{CCCommutators}
\end{align}
Note that adjoint terms are absent also in $[L_1,J_{\mu\nu}]$.
$L_1$ acts at $k=-2$ as $1\otimes\ell_1+\ell_1\otimes1$
in some basis,
and, in addition, there is a Jordan cell structure, so eq. \eqref{LDeltaTwo}
can be rewritten as
\be
[L_1,\left(
\begin{matrix}
T_{\mu\nu}\\ J_{\mu\nu}
\end{matrix}
\right)]
=\left(
\begin{matrix}
-r_{\mu\nu}{}^{\kappa\lambda}&-\delta_{\mu\nu}^{\kappa\lambda}\\
0&-r_{\mu\nu}{}^{\kappa\lambda}
\end{matrix}
\right)
\left(
\begin{matrix}
T_{\kappa\lambda}\\ J_{\kappa\lambda}
\end{matrix}
\right)\;.\label{LDeltaLevelTwo}
\ee
where
\be
r=(\Cn0-1+\sigma)^{{1-\sigma\over2}}(1\otimes\ell_1+\ell_1\otimes1)
(\Cn0-1+\sigma)^{-{1-\sigma\over2}}\;.
\ee
In a suitable basis, this representation matrix is simply given by $1\otimes\ell_1+\ell_1\otimes1$.
The individual irreducible affine modules in $s_2$ and $a_2$ are not
preserved by the action of $L_1$.

Similarly, the Jacobi identities $[L_1,T_{\mu\nu},\bar T^\kappa]$
and $[L_1,J_{\mu\nu},\bar T^\kappa]$ are shown to hold, using
eq. \eqref{CCCommutators}.

\subsection{The closure of the algebra\label{ClosureSection}}
We would
like to show by induction in $-k$ that nothing
beyond $T$'s and $J$'s is generated. The considerations above show
this to hold to affine level $k=-2$ (and this is the critical level, where null states appear). Assume that we have established that
the generators at $k=-n$ consists of $T^{(n)}=\langle T_{\mu_1\ldots\mu_n}\rangle$
and $J^{(n)}=\langle J_{\mu_1\ldots\mu_n}\rangle$, where, by definition,
$J_{\mu_1\ldots\mu_n}=[T_{\mu_1},J_{\mu_2\ldots\mu_n}]$ modulo $T^{(n)}$.
Then consider $[J_\mu,J_{\mu_1\ldots\mu_n}]$. We will show that this
can be written in terms of expressions with one $J$ only, modulo an
ideal. We therefore act with a $\bar T$.
$[\bar T^\nu,[J_\mu,J_{\mu_1\ldots\mu_n}]]$ can (using the Jacobi
identity)
by the assumption be written
in terms of (some matrices acting on) $J^{(n)}$ and $T^{(n)}$.
The only observation needed is that any expression with $J^{(n)}$ can
be obtained as $[\bar T,J^{(n+1)}]$ (modulo terms with $T^{(n)}$). This is only based on the fundamental $\fg^{++}$ representation $R(-\lambda)$.
The same argument holds for $[J^{(m)},J^{(n+1-m)}]$.
In all, this shows that there is an ideal at affine level $k=-(n+1)$ that can be used to
remove anything beyond the fundamental and adjoint modules.
Since the assumption holds for $n=2$, we show by induction that it
holds for all $n\in{\mathbb N}$.

We can now continue to show that $[J_M,J_N]=g_{MN}{}^PJ_P$, \ie, that
the bracket between two $J$'s does not produce any terms with
$\fg^{++}$ generators $T_\alpha$, to all affine levels.
The idea is to consider the action of the lowest state $J$, already identified as acting as $L_1$, with some modification. Generators $J_M$ in the fundamental will be identified as precisely those transforming homogeneously under $J$.

A basis for the part of $\fg^{++}$ at negative affine level $k$ is
provided by $T_{\mu_1\ldots\mu_n}=[T_{\mu_1},\ldots,T_{\mu_n}]$.
Let us consider the transformation under $L_1=J$. We obtain
\begin{align}
&[J,T_{\mu_1\ldots\mu_n}]\nn\\
&=-(\ell_1\otimes1\otimes\ldots\otimes1
+1\otimes\ell_1\otimes\ldots\otimes1\nn\\
&\qquad\qquad+\ldots+1\otimes1\otimes\ldots\otimes\ell_1
)_{\mu_1\ldots\mu_n}{}^{\nu_1\ldots\nu_n}T_{\nu_1\ldots\nu_n}\\
&-([J_{\mu_1},T_{\mu_2},\ldots,T_{\mu_n}]
+[T_{\mu_1},J_{\mu_2},T_{\mu_3},\ldots,T_{\mu_n}]+\ldots
+[T_{\mu_1},\ldots,T_{\mu_{n-1}},J_{\mu_n}])\;.\nn
\end{align}
The terms in the last row look like the inhomogeneous part of the 
transformation. This is however only so for $n=1$, for higher $n=-k$
it will consist of a mixture of $T$ and $J$. An example of this for
$n=2$ is displayed by eq. \eqref{LDeltaLevelTwo}. Continuing to act with
$J$ one finally reaches the expression
\be
\tilde J_{\mu_1\ldots\mu_n}=[J_{\mu_1},\ldots,J_{\mu_n}]\;,\label{JTildeDef}
\ee
which transforms homogeneously,
\begin{align}
[J,\tilde J_{\mu_1\ldots\mu_n}]
&=-(\ell_1\otimes1\otimes\ldots\otimes1
+1\otimes\ell_1\otimes\ldots\otimes1\nn\\
&\qquad\qquad+\ldots+1\otimes1\otimes\ldots\otimes\ell_1
)_{\mu_1\ldots\mu_n}{}^{\nu_1\ldots\nu_n}\tilde J_{\nu_1\ldots\nu_n}\;.
\end{align}
This is a pure $J$ generator. Homogeneous transformation then might seem to be
equivalent
to a generator not containing $T$. This is indeed the case, as is argued below.
Homogeneously transforming generators form a subalgebra, but they are not spanned only by the generators $\tilde J$. This is due to the non-simplicity of the subalgebra of $J$'s. It can be noted that the 
$\tilde J$'s above are not enough to obtain the whole fundamental module. One
also needs the generators
$K_{\mu\nu}=[T_{(\mu},J_{\nu)}]$ and $L_{\mu\nu\lambda}=[T_{[\mu},J_\nu,J_{\lambda]}]$, filling in
plethysms present in the fundamental but not in the adjoint, and
transforming homogeneously, thanks to antisymmetry of the bracket and
the Jacobi identity. The inhomogeneous term received from the transformation of the single $T_\mu$ results in terms that vanish identically.
Also, in the light of the properties of the Virasoro characters mentioned earlier, there will be additional parts of the fundamental module that are not obtained in $\tilde J_{\mu_1\ldots\mu_n}$. They will also need to be added, much like the generators $K$ and $L$, although explicit expressions, containing some $T_\mu$, will be complicated, involving projections on individual coset Virasoro states.
We can however be sure that they transform homogeneously under $J$, since the inhomogeneous terms arising from the transformation of a $T_\mu$ will be precisely of the forms that do not appear in the shifted 
expressions, similarly to what happened with $K$ and $L$ above.
Note that the basis used here for the $\tilde J$'s is not the same as used earlier in the affine expansion. An expression \eqref{JTildeDef} is shifted by $n-1$ in mode number compared to $J_{\mu_1\ldots\mu_n}$, and will be obtained from it by multiplication of some invariant tensor. See for example the last
equation in \eqref{OneOneTwo}, which contains the invariant tensor $\Cn1$.

The important observation is that all elements in the fundamental are generated by the homogeneously transforming generators, that this property is preserved by the bracket, and that the fundamental module, spanned by $J_M$ therefore forms a subalgebra. 

\subsection{Some representations}
The algebra $\AA$ defined by the brackets \eqref{TJBracket} can be represented on $F=R(-\lambda)$ or $\bar F=R(\lambda)$, the fundamental modules, alone. 
This is almost obvious, since
it is generated from $L_1$, $e_{-1}$ and $T^A_m$, all of
which have well-defined actions on $F$ as a lowest weight module.

Now, consider the bracket between an element in the module
$R(\lambda)$ with highest weight state $\Phi$ and some element in the
weakly lowest weight module $R(-\lambda)$ with lowest weight state
$J$.
(We examine the representation by forming brackets for the semidirect product 
$\AA\loplus F$.)
We thus consider brackets of the type
\begin{align}
  [\Phi',J']=[[f_{i_1},\ldots,f_{i_p},\Phi],[e_{j_1},\ldots,e_{j_q},J]]\;.
  \label{PhiJBrackets}
\end{align}

Repeated use of Jacobi identities of the types
$[\Phi',[e,J']]=[[\Phi',e],J']+[e,\Phi',J']$ and
$[[f,\Phi'],J']=[f,\Phi',J']-[\Phi',f,J']$
relates expressions \eqref{PhiJBrackets} to expressions of the same
form with lower $p$ or $q$. There is no need to keep track of the weak
lowest state property of $J$, since it only produces terms without
$\epsilon$.
In the end, we only need to rely on the defining condition $[\Phi,J]=0$.
This shows that $R(\lambda)$ with highest weight state $\Phi$ is a
module of the algebra $\AA$ spanned by $T_\alpha$ and $J_M$.

One example of such a module, appearing in the tensor hierarchy algebra 
$S(\fg^{++})$, appears at level $-2$. The highest weight state then is
$\Phi=[\phi,e_0,e_\theta,\phi]$, where $\phi=\phi_{-2,1}$, see Section
\ref{SgplusplusGenRelSection}.

The brackets can be represented by representation
matrices
$t_{\alpha M}{}^N$ and $j_{MN}{}^P$ with commutators
\begin{align}
  [t_\alpha,t_\beta]&=f_{\alpha\beta}{}^\gamma t_\gamma\;,\label{FundIdI}\\
  [t_\alpha,j_M]&=-t_{\alpha M}{}^Nj_N
      +u_{M\alpha}{}^\beta t_\beta\;,\label{FundIdII}\\
  [j_M,j_N]&=g_{MN}{}^Pj_P\;.\label{FundIdIII}
 \end{align}

We can derive the beginning of the representation matrices $j_{MN}{}^P$
for the fundamental representation in the affine decomposition. We check the identities \eqref{FundIdI}-\eqref{FundIdIII}
by extending the bracket to the semidirect sum of the algebra with the
fundamental module, with basis elements
$E_M=(E,E_\mu,E_{\mu\nu},\ldots)$, and checking the Jacobi identities.
The brackets $[T,E]$ are given, and include (beyond the affine
transformations already manifest) 
\begin{align}
[\bar T^\mu,E]&=0\;,\qquad[T_\mu,E]=E_\mu\;,
\qquad[T_{\mu\nu},E]=E_{[\mu\nu]}\;,\nn\\
[\bar T^\mu,E_\nu]&=\delta^\mu_\nu E\;,
\qquad[T_\mu,E_\nu]=E_{\mu\nu}\;,\label{FundamentalModuleI}\\
[\bar T^\mu,E_{\nu\kappa}]
&=(\Cn0-1+\sigma)_{\nu\kappa}{}^{\mu\lambda}E_\lambda\;.\nn
\end{align}
The brackets $[J_M,E_N]=-j_{MN}{}^PE_P$ consistent with these are
\begin{align}
[J,E]&=0\;,\nn\\
[J,E_\mu]&=-(\ell_1)_\mu{}^\nu E_\nu\;,\nn\\
[J_\nu,E]&=0\;,\label{FundamentalModuleII}\nn\\
[J,E_{\mu\nu}]&=-((\Cn0-1+\sigma)(1\otimes\ell_1+\ell_1\otimes1)
    (\Cn0-1+\sigma)^{-1})_{\mu\nu}{}^{\kappa\lambda}E_{\kappa\lambda}\;,\\
[J_\mu,E_\nu]&=-(\Cn1(\Cn0-1+\sigma)^{-1})_{\mu\nu}{}^{\kappa\lambda}E_{\kappa\lambda}\;,\nn\\
[J_{\mu\nu},E]&=(\Cn1\Cn0{}^{-1})_{(\mu\nu)}{}^{\kappa\lambda}E_{\kappa\lambda}\;.\nn
\end{align}
We have checked all Jacobi identities $[T,J,E]$ and $[J,J,E]$
within $k=1,0,-1,-2$ involving elements at $k=1,0,-1,-2$.
 
 It may seem peculiar that there are two sets of representation matrices with the same index structure, $g_{MN}{}^P$ for the $F$ subalgebra of $\AA$ and $j_{MN}{}^P$ for the fundamental representation of $\AA$. They turn out to be related. If we list a few of them in the affine decomposition, we have
\begin{align}
j_{0,0}{}^0&=0\;,\nn\\
j_{0,\mu}{}^\nu&=(\ell_1)_\mu{}^\nu\;,\nn\\
j_{\mu,0}{}^\nu&=0\;,\nn\\
j_{0,\mu\nu}{}^{\kappa\lambda}
&=((\Cn0-1+\sigma)(1\otimes\ell_1+\ell_1\otimes1)
    (\Cn0-1+\sigma)^{-1})_{\mu\nu}{}^{\kappa\lambda}\;,\label{jcoefficients}\\
j_{\mu,\nu}{}^{\kappa\lambda}
            &=(\Cn1(\Cn0-1+\sigma)^{-1})_{\mu\nu}{}^{\kappa\lambda}\;,\nn\\
j_{\mu\nu,0}{}^{\kappa\lambda}&=-(\Cn1\Cn0{}^{-1})_{(\mu\nu)}{}^{\kappa\lambda}\;,\nn
\end{align}
and
\begin{align}
g_{0,\mu}{}^\nu&=-(\ell_1)_\mu{}^\nu\;,\nn\\
g_{0,\mu\nu}{}^{\kappa\lambda}
           &=-((\Cn0-1+\sigma)^{{1-\sigma\over2}}(1\otimes\ell_1+\ell_1\otimes1)
(\Cn0-1+\sigma)^{-{1-\sigma\over2}})_{\mu\nu}{}^{\kappa\lambda}\;,\\
g_{\mu,\nu}{}^{\kappa\lambda}&=-2(\Cn1(\Cn0-2)^{-1})_{[\mu\nu]}{}^{\kappa\lambda}\;.\nn
\end{align}
We observe that they satisfy 
\be
g_{MN}{}^P=-2j_{[MN]}{}^P\;.
\ee
We also note that the Jacobi identities \eqref{gtid} and \eqref{ggid} are
automatically satisfied when eqs. \eqref{FundIdI}-\eqref{FundIdIII} hold and
$g_{MN}{}^P=-2j_{[MN]}{}^P$. This is the only possible relation between $g$ and $j$.
We can prove that it holds in the following way.
Linearly combining eqs. \eqref{gtid} and \eqref{FundIdII} gives
\begin{align}
0&=
t_{\alpha M}{}^Q(g_{QN}{}^P+2j_{QN}{}^P)+t_{\alpha N}{}^Q(g_{MQ}{}^P+2j_{MQ}{}^P)\nn\\
&-t_{\alpha Q}{}^P(g_{MN}{}^Q+2j_{MN}{}^Q)
-t_{\beta M}{}^Pu_{N\alpha}{}^\beta-t_{\beta M}{}^Pu_{N\alpha}{}^\beta\;.
\end{align}
Antisymmetrisation in $[MN]$ states that $g_{MN}{}^P+2j_{[MN]}{}^P$ is a
$\fg^{++}$-invariant tensor. Since $\wedge^2R(-\lambda)\not\supset R(-\lambda)$ as 
$\fg^{++}$ modules, it vanishes.

Consider the tensor product $F\otimes F$. As a module of the over-extended
algebra, it will decompose in an infinite number of irreducible
modules,
\be
\label{FOtimesF}
F\otimes F=(S_2^\compl\oplus S_2)\oplus (A_2^\compl\oplus A_2)\;,
\ee
where $S_2^\compl$ and $A_2^\compl$ are the leading irreducible symmetric and antisymmetric 
$\fg^{++}$-modules, with lowest weights $-2\lambda$ and $-(2\lambda-\alpha_{-1})$, respectively.
$S_2$ and $A_2$ are infinitely reducible as $\fg^{++}$-modules.

We can perform the analysis for $k=0,-1,-2$.
Given the content of $F$ at $k=0,-1,-2$, $F=(1,f,s_2\oplus
a_2,\ldots)$, we get
\begin{align}
\vee^2F&=(1,f,s_2^\compl\oplus2s_2\oplus a_2,\ldots)\;\nn\\
\wedge^2F&=(0,f,s_2\oplus a_2^\compl\oplus2a_2,\ldots)\;.
\end{align}
On the other hand, it is straightforward to verify that
\begin{align}
S_2^\compl&=(1,f,s_2^\compl\oplus s_2\oplus a_2,\ldots)\;\nn\\
A_2^\compl&=(0,f,s_2\oplus a_2^\compl\oplus a_2,\ldots)\;.
\end{align}
The leading ``missing'' states in $s_2$ and $a_2$ are given by
\begin{align}
S_{\mu\nu}&=E\otimes E_{(\mu\nu)}+E_{(\mu\nu)}\otimes E
+\Cn0_{(\mu\nu)}{}^{\kappa\lambda}E_\kappa\otimes E_\lambda\;,\nn\\
A_{\mu\nu}&=E\otimes E_{[\mu\nu]}-E_{[\mu\nu]}\otimes E
+(\Cn0-2)_{[\mu\nu]}{}^{\kappa\lambda}E_\kappa\otimes E_\lambda\;.
\end{align}
They are annihilated by $\bar T^\lambda$, and therefore provide the lowest
affine modules in an infinite number of irreducible over-extended modules
in $S_2$ and $A_2$. From the discussion above it follows that they are
not necessarily reducible as modules of $\AA$. This of course continues at higher
negative $k$. $S_2$ and $A_2$ are $\AA$-modules, and we conjecture that they may be irreducible.

\subsection{Covariance\label{CovarianceSection}}

The structure constants and representation matrices are not quite ``tensors'' in the usual sense, since the generators in $F$ do not (only) transform as elements of the fundamental module of $\fg^{++}$. 
The structure constants $u_{M\alpha}{}^\beta$ and $j_{MN}{}^P$ are not (and should not be) invariant tensors under $\fg^{++}$.
This is quite an unfamiliar situation. 
The all-important property is of course that the Jacobi identities are satisfied.
A useful tool is to view identities for ``tensors'' as specifying the deviation from invariance.
To this end, we introduce the symbols $\Delta_\alpha$ and $\Delta_M$, measuring precisely this.
The definitions are
\begin{align}
\Delta_\alpha v_\beta&=-f_{\alpha\beta}{}^\gamma v_\gamma\;,\nn\\
\Delta_\alpha v_M&=t_{\alpha M}{}^Nv_N\;,\nn\\
\Delta_M v_\alpha&=u_{M\alpha}{}^\beta v_\beta\;,
\label{DeltaDefinitions}\\
\Delta_M v_N&=j_{MN}{}^Pv_P\;,\nn\\
\Delta_M v_{\dot N}&=-g_{MN}{}^Pv_P\;\nn
\end{align}
and, of course, distributive action on tensor products. Note that there are two distinct actions of $\Delta_M$ on a fundamental index; we distinguish the fundamental in $\AA$ by a dotted index.

All identities in eqs. \eqref{ffid}--\eqref{ggid} and \eqref{FundIds} are then equivalent to the statements
\begin{align}
\Delta_\alpha f_{\beta\gamma}{}^\delta&=0\;,\nn\\
\Delta_\alpha t_{\beta M}{}^N&=0\;,\nn\\
\Delta_\alpha j_{MN}{}^P&=u_{M\alpha}{}^\beta t_{\beta N}{}^P\;,
       \label{DeltaAlphaCoeffs}\\
2\Delta_{[\alpha}u_{|M|\beta]}{}^\gamma&=-f_{\alpha\beta}{}^\delta u_{M\delta}{}^\gamma\;,\nn
\end{align}
and
\begin{align}
\Delta_M f_{\alpha\beta}{}^\gamma&=-2u_{N[\alpha}{}^\gamma t_{\beta]M}{}^N\;,\nn\\
\Delta_M t_{\alpha N}{}^P&=t_{\alpha M}{}^Qj_{QN}{}^P\;,\nn\\
\Delta_M j_{\dot NP}{}^Q&=0\;,
         \label{DeltaMCoeffs}\\
\Delta_Mu_{\dot N\alpha}{}^\beta&=0\;.\nn
\end{align}
The last equation in \eqref{DeltaAlphaCoeffs} and the first one in \eqref{DeltaMCoeffs}
are equivalent statements of the transforming derivation property
\eqref{fuid}. These equations will be used extensively in the superalgebra $S(\fg^{++})$.

If we consider commuting the $\Delta$'s, we quite obviously find
\begin{align}
[\Delta_\alpha,\Delta_\beta]&=-f_{\alpha\beta}{}^\gamma\Delta_\gamma\;,\nn\\
[\Delta_{\dot M},\Delta_{\dot N}]&=-g_{MN}{}^P\Delta_P\;.
\end{align}
What is more surprising is that 
\be
[\Delta_M,\Delta_N]=0\;,
\ee
a kind of flatness property, which will turn of to be important.
It follows as
\begin{align}
[\Delta_M,\Delta_N]v_\alpha&=\Delta_M(u_{N\alpha}{}^\beta v_\beta)-(M\leftrightarrow N)\nn\\
&=(\Delta_M u_{N\alpha}{}^\beta)v_\beta+u_{N\alpha}{}^\beta\Delta_Mv_\beta
                        -(M\leftrightarrow N)\nn\\
&=(\Delta_M u_{\dot N\alpha}{}^\beta+(j_{MN}{}^P+g_{MN}{}^P)u_{P\alpha}{}^\beta)v_\beta
                        +u_{N\alpha}{}^\beta u_{M\beta}{}^\gamma v_\gamma\nn\\
             &\qquad-(M\leftrightarrow N)\\
&=(g_{MN}{}^Pu_{P\alpha}{}^\beta-[u_M,u_N]_\alpha{}^\beta)v_\beta\nn\\
&=0\nn
\end{align}
(since $\Delta_M u_{\dot N\alpha}{}^\beta=0$ and $g_{MN}{}^P=-2j_{[MN]}{}^P$),
and a similar calculation for $v_P$.
It is also straightforward to show that
\begin{align}
[\Delta_\alpha,\Delta_M]v_\beta&=t_{\beta M}{}^Nu_{N\alpha}{}^\gamma v_\gamma\;,\nn\\
[\Delta_\alpha,\Delta_M]v_N&=0\;,\label{DeltaAlphaDeltaMCommutator}\\
[\Delta_\alpha,\Delta_M]v_{\dot N}&=-t_{\beta M}{}^Pu_{N\alpha}{}^\beta v_P\;.\nn
\end{align}

\section{The tensor hierarchy extension of an over-extended algebra
\label{SgplusplusSection}}

\subsection{Generators and relations\label{SgplusplusGenRelSection}}

The action of $J$ as $L_1$ in the affine subalgebra follows from Section \ref{SAffineSection}, since $S(\fg^+)$ is a subalgebra of $S(\fg^{++})$.
In addition to vanishing of $[f_i,J]$, ($i=0,\ldots,r$) obtained there, we need to calculate the action of $f_{-1}$. We use the relations \eqref{oddweyl}, to express $J$ in $S(\fg^{++})$ in terms of generators in the standard basis, \ie, with only one fermionic node, obtaining
\be
J=[\epsilon_{-2},e_{-1},e_0,e_\theta,f_{-1},\phi_{-2,1}]\;.
\ee 
Then,
\begin{align}
[f_{-1},J]=[\epsilon_{-2},e_{-1},e_0,e_\theta,f_{-1},f_{-1},\phi_{-2,1}]=0\;
\end{align}
thanks to eq. \eqref{ffphiis0}.
The weakly lowest weight state $J$ is thus annihilated by all $f_i$'s
except $f_0$, under which it is annihilated modulo an element in the
adjoint. Its $\fg^{++}$ weight is $\alpha_0+\theta=-\Lambda_{-1}=-\lambda$. It
is the lowest weight state in a fundamental $\fg^{++}$ module $R(-\lambda)$, and
acts as $L_1$ on the $\fg^+$ subalgebra.

The only raising operator that can be applied to $J$, without getting
only the adjoint, is $e_{-1}$ (this also follows from Section \ref{SAffineSection}).
We can write $J$ as
$J=[f_{-1},\epsilon_{-2},e_{-1},e_0,e_\theta,\phi_{-2,1}]$. Then,
\begin{align}
  [e_{-1},J]&=[h_{-1},\epsilon_{-2},e_{-1},e_0,e_\theta,\phi_{-2,1}]
            +[f_{-1},e_{-1},\epsilon_{-2},e_{-1},e_0,e_\theta,\phi_{-2,1}]\;.
\end{align}
In the second term, we use
\begin{align}
  [\ldots,e_{-1},\epsilon_{-2},e_{-1},\ldots]
=\tfrac12([\ldots,e_{-1},e_{-1},\epsilon_{-2},\ldots]
+[\ldots,\epsilon_{-2},e_{-1},e_{-1},\ldots])\;,\label{SerreComm}
\end{align}
which holds for singly connected simple raising operators thanks to
the Serre relation
$[e_{-1},e_{-1},\epsilon_{-2}]=0$. The first of these trivially gives $0$ acting on
$[e_0,e_\theta,\phi_{-2,1}]$. In the second one, we use the same trick
to get at least one $e_{-1}$ past $e_0$, so it also vanishes.
Thus,
\begin{align}
[e_{-1},J]&=[\epsilon_{-2},e_{-1},e_0,e_\theta,\phi_{-2,1}]\;.
\end{align}
An even simpler calculation uses the ``inside-out'' form of $J$,
$J=[\phi_{-1,1},e_\theta,e_0,\epsilon_{-1}]
=-[f_{-1},\phi_{-2,1},e_\theta,e_0,e_{-1},\epsilon_{-2}]$, to obtain
\begin{align}
[e_{-1},J]&=-[\phi_{-2,1},e_\theta,e_0,e_{-1},\epsilon_{-2}]\;.
\end{align}
Here, we have used $[e_{-1},e_0,e_{-1},\epsilon_{-2}]=0$ by the same method
as in eq. \eqref{SerreComm}.

We now want to calculate $[J,[e_{-1},J]]$, as a first example of a
bracket between elements in the fundamental.
Set $\epsilon\deltaeq\epsilon_{-2}$, $\phi\deltaeq\phi_{-2,1}$.
We have
\begin{align}
  [J,[e_{-1},J]]&=[J,\epsilon,e_{-1},e_0,e_\theta,\phi]\nn\\
  &=\underset{(1)}{[\epsilon,[J,e_{-1}],e_0,e_\theta,\phi]}
  +\underset{(2)}{[\epsilon,e_{-1},[J,e_0],e_\theta,\phi]}\label{JeJBracket}\\
  &\qquad+\underset{(3)}{[\epsilon,e_{-1},e_0,e_\theta,J,\phi]}\;.\nn
\end{align}

First, determine $[J,\phi]$. We always use the identities (\ref{e-identitet}) and (\ref{f-identitet}), which in this case read
$[e_{-1},e_0,\phi]=0$ and 
$[f_{-1},e_0,\phi]=0$, 
allowing us to pass $e_{-1}$ or $f_{-1}$ through $[e_0,\phi]$.
We have
$[J,\phi]=[\epsilon,X]$,
where $X=[\phi,e_{-1},e_0,e_\theta,f_{-1},\phi]$.
Also, let $Y=[\phi,e_0,e_\theta,\phi]$. We will show that $X$ and $Y$ are
proportional. This is desirable, since we expect a single generator in
$S(\fg^{++})$ at level $-2$ and $\fg^{++}$ weight $\lambda$.

It is straightforward to show that $X$ and $Y$ are both annihilated by
$e_{-1}$ and $e_0$. They have eigenvalue $1$ under $h_{-1}$ and $0$
under $h_0$, so they are both highest weight states in a triplet
representation
of the $A_2$ subalgebra corresponding to the nodes $-1,0$.
(They are also annihilated by $e_i$, $i=1,\ldots,r$, and provide
highest weight states in the $\fg^{++}$ module $R(\lambda)$.)
We calculate the other states 
in the triplets. Acting with $f_{-1}$ gives
\begin{align}
  [f_{-1},X]&=[\phi,e_0,e_\theta,f_{-1},\phi]\;,\\
  [f_{-1},Y]&=[\phi,e_0,e_\theta,f_{-1},\phi]+[[f_{-1},\phi],e_0,e_\theta,\phi]\nn\\
  &=[\phi,e_0,e_\theta,f_{-1},\phi]+[e_0,[[f_{-1},\phi],e_\theta],\phi]\\
  &=[f_{-1},X]-[e_0,\phi,e_\theta,f_{-1},\phi]\;.\nn
\end{align}
We then act with $f_0$. In this step of the calculation, we need to use
$[e_0,f_0,\phi]=-\phi_{-2,0}\deltaeq-\phi'$ and
$[e_\theta,\phi']=\frac{(\alpha_0,\theta)}{(\alpha_1,\theta)}[e_\theta,\phi]
=-2[e_\theta,\phi]$, together with $[\phi',f_{-1},\phi]=0$, to obtain
\begin{align}
  [f_0,f_{-1},X]
  &=[[f_0,\phi],e_0,e_\theta,f_{-1},\phi]+[\phi,f_0,e_0,e_\theta,f_{-1},\phi]\nn\\
  &=[[[f_0,\phi],e_0],e_\theta,f_{-1},\phi]
  +\cancel{[e_0,[f_0,\phi],e_\theta,f_{-1},\phi]}\nn\\
  &\qquad+[\phi,e_\theta,f_{-1},\phi]\\
  &=[\phi',e_\theta,f_{-1},\phi]+[\phi,e_\theta,f_{-1},\phi]\nn\\
  &=-[\phi,e_\theta,f_{-1},\phi]\;,\nn\\
  [f_0,f_{-1},Y]
  &=[f_0,f_{-1},X]-[\phi,e_\theta,f_{-1},\phi]\nn\\
  &=-2[\phi,e_\theta,f_{-1},\phi]\;.
\end{align}
This shows that $X=\frac12Y$. 
 
Going back to the terms in eq. \eqref{JeJBracket}, they are
evaluated as
\begin{align}
  (1)&=-[\epsilon,[,e_{-1},e_0,e_\theta,\phi],e_0,e_\theta,\phi]\nn\\
  &=[[,e_{-1},e_0,e_\theta,\phi],e_0,e_\theta]
  =-[\epsilon,[e_0,e_\theta],e_{-1},e_0,e_\theta,\phi]\nn\\
  &=[\epsilon,e_\theta,e_0,e_{-1},e_0,e_\theta,\phi]\\
  &=\tfrac12[\epsilon,e_\theta,e_0,e_0,e_{-1},e_\theta,\phi]
  +\tfrac12[\epsilon,e_\theta,e_{-1},e_0,e_0,e_\theta,\phi]\nn\\
  &=\tfrac12[\epsilon,e_{-1},e_0,e_\theta,e_0,e_\theta,\phi]\;,\nn\\
  (2)&=\tfrac12[\epsilon,e_{-1},[e_0,e_0,e_\theta],e_\theta,\phi]\nn\\
  &=\tfrac12[\epsilon,e_{-1},e_0,[e_0,e_\theta],e_\theta,\phi]
  -\tfrac12[\epsilon,e_{-1},[e_0,e_\theta],e_0,e_\theta,\phi]\\
  &=-\tfrac12[\epsilon,e_{-1},e_0,e_\theta,e_0,e_\theta,\phi]\;,\nn\\
  (3)&=\tfrac12[\epsilon,e_{-1},e_0,e_\theta,\epsilon,\phi,e_0,e_\theta,\phi]\nn\\
  &=\tfrac12[\epsilon,e_{-1},e_0,e_\theta,[e_0,e_\theta],\phi]=0\;.
\end{align}
Here, we have used $[e_0,[e_0,e_\theta],\phi]=0=[e_\theta,[e_0,e_\theta],\phi]$
($[[e_0,e_\theta],\phi]$ is the central element in the $\fg^+$ adjoint),
as well as $[[e_0,e_\theta],e_0,e_\theta,\phi]=0$
and $[e_\theta,e_{-1},e_0,e_\theta,\phi]=0$.
Altogether,
\begin{align}
[J,[e_{-1},J]]=0\;,\label{JeJis0}
\end{align}
as desired. Eq. \eqref{JeJis0} can now be used as a basic building block when
considering elements in the algebra at level $0$ in $S(\fg^{++})$,
constructed sequentially by the eigenvalue of $\KK=-h_\lambda$.

At $k=-1$, we have elements in $\fg^{++}$ of the type
$[t_1,\ldots,t_p,e_{-1}]$, where each of $t_1,\ldots,t_p$ is one of the generators
$e_0,...,e_r$. We think of $e_{-1}$ as the lowest state in a
fundamental affine module, which is annihilated by $L_1$.
Thus, acting with $J$ gives the action of $L_1$ in the affine module,
plus a term $[t_1,\ldots,t_p,J,e_{-1}]$. Identifying $[J,e_{-1}]$ as
the lowest weight state in the (shifted) affine fundamental at $k=-1$
in the $\fg^{++}$ module $R(-\lambda)$, we arrive at
\begin{align}
  [J,T_\mu]&=-(\ell_1)_\mu{}^\nu T_\nu-J_\mu\;,\nn\\
  [J,J_\mu]&=-(\ell_1)_\mu{}^\nu J_\nu\;.
\end{align}

At level $-1$ in $S(\fg^{++})$, we have $\phi_i=\phi_{-2,\alpha_i}$, $i=0,\ldots,r$. 
They have $\fg^{++}$ weight $\lambda=\Lambda_1$. Notice that $\lambda=-\delta=-\alpha_0-\theta$ is light-like, and that it lies in the $\fg^+$ root space. We can choose a basis
with $\phi_\lambda$ and $\phi_j$, $j=1,\ldots,r$. 
Since $(\lambda,\alpha_i)=0$, $i=0,\ldots,r$, we find that $\phi_\lambda$ is highest weight in a $\fg^{++}$-module $R(\lambda)$. The remaining $\phi_j$, $j=1,\ldots,r$ turn out to be Cartan generators at mode number $-1$ in a $\fg^+$ adjoint part of a $\fg^{++}$ adjoint. 
We can also deduce that the transformation of this ``adjoint'' under $\fg^{++}$ however also ``leaks'' into the module $R(\lambda)$. 
This corroborates the information obtained (more easily) from duality in the following subsection, eq. \eqref{LevelMinusOneTransf}.

At level $-2$, we find a highest weight state $\Phi=[\phi,e_0,e_\theta,\phi]$ for an anti-fundamental $R(\lambda)$, paired with $E_M$ at level $1$ through the invariant bilinear form.
(There is also an infinite number of highest weight modules in $\overline{S_2}$, which have much more complicated expressions in terms of the generators.)
In addition to the identities $[e_a,\Phi]=0$, we only need to verify that 
$[J,\Phi]=0$ in order to show that this $R(\lambda)$ forms an anti-fundamental $\AA$-module.
To calculate this bracket, we write
\begin{align}
  [\Phi,J]=[[\phi,J],e_0,e_\theta,\phi]
  +[\phi,[e_0,J],e_\theta,\phi]+[\phi,e_0,e_\theta,[\phi,J]]\;,
  \label{PhiJZeroBracket}
\end{align}
where $[e_0,J]=-\frac12[e_0,e_0,e_\theta]$ and
\begin{align}
  [\phi,J]&=[\phi,\epsilon,e_{-1},e_0,e_\theta,f_{-1},\phi]\nn\\
  &=\cancel{[h_1,e_{-1},e_0,e_\theta,f_{-1},\phi]}
  -[\epsilon,\phi,e_{-1},e_0,e_\theta,f_{-1},\phi]\nn\\
  &=-\frac12[\epsilon,\phi,e_0,e_\theta,\phi]\\
  &=-\frac12\cancel{[h_1,e_0,e_\theta,\phi]}
  +\frac12[\phi,\epsilon,e_0,e_\theta,\phi]\nn\\
  &=\frac12[\phi,e_0,e_\theta,h_1]
        =\frac12[[e_0,e_\theta],\phi]\nn\;.
\end{align}
All three terms in eq. \eqref{PhiJZeroBracket} vanish.

We notice that, although the generators and defining relations provide a concise definition of the superalgebra, calculations soon become complicated. It becomes more efficient to use the 
``covariant'' tensor formalism, with its basic properties obtained from the definitions, to derive further properties (sometimes supported by the 
affine decomposition) as in Section
\ref{ExtendByFund}. In the following, we will rely on the results from that Section, in order to investigate the $\fg^{++}$-covariant decomposition of $S(\fg^{++})$ in the grading with respect to the fermionic node.

\subsection{The local superalgebra\label{LocalSubalgebraSection}}
Level $0$, in the grading with respect to the fermionic node, consists of the algebra $\AA$ described in Section \ref{ExtendByFund}, consisting of $\fg^{++}$ and the module $R(-\lambda)$, with generators $T_\alpha$ and $J_M$.
This is also consistent with the $\gl$-grading, as described in Section \ref{GLGradingSection}, applied to 
$S(\fg^{++})$.

From the existence of the non-degenerate invariant bilinear form we know that level $-1$ is dual to level $0$, \ie, it consists of the coadjoint module of $\AA$. Let the basis elements be $V^M$, $U^\alpha$. They  transform according to
\begin{align}
[T_\alpha,V^M]&=t_{\alpha N}{}^MV^N\;,\nn\\
[T_\alpha,U^\beta]&=-f_{\alpha\gamma}{}^\beta U^\gamma
               -u_{M\alpha}{}^\beta V^M\;,\nn\\
[J_M,V^N]&=-g_{MP}{}^NV^P-t_{\alpha M}{}^NU^\alpha\;,\label{LevelMinusOneTransf}\\
[J_M,U^\alpha]&=u_{M\beta}{}^\alpha U^\beta\;.\nn
\end{align}
Note in particular that $U^\alpha$, in spite of spanning a $\fg^{++}$ adjoint, in addition transforms into $V^M$.

At level $1$ there is of course $\{E_M\}$, as a basis for $\BB_1=R(-\lambda)$, with lowest weight state $\epsilon_{-2}$.
We can then use information from the $\gl$-grading, where we find, in addition to $\BB_1$, both leading and some subleading parts of $S_2$ (defined in eq. \eqref{FOtimesF}). We thus introduce the level $1$ generators $\{E_{MN}\}$ as a basis for $S_2$.
How do these transform under $\AA$?
The generators in $S_2$ will contain one uncancelled $\phi$. It must in some way be generated in the bracket $[J_M,E_N]$. An investigation of this in terms of the generators should be possible, but becomes awkward, mainly because even the lowest of the lowest weight states in $S_2$ is quite high. Instead we rely on the tensor formalism developed in Section \ref{ExtendByFund}.

We can thus postulate the transformations
\begin{align}
[T_\alpha,E_M]&=-t_{\alpha M}{}^NE_N\;,\nn\\
[T_\alpha,E_{MN}]&=-2t_{\alpha(M}{}^PE_{N)P}\;,\nn\\
[J_M,E_N]&=-j_{MN}{}^PE_P+E_{MN}\;,\\
[J_M,E_{NP}]&=-2j_{M(N}{}^QE_{P)Q}\;.\nn
\end{align}
The only ``non-covariant'' ingredient is the Jordan cell structure in the action of $J$. This describes an $\AA$-module. The only non-trivial Jacobi identity needed to show this is
\begin{align}
&[[J_M,J_N],E_P]-2[J_{[M},[J_{N]},E_P]]\nn\\
&=[g_{MN}{}^Q,E_P]-2[J_{[M},-j_{N]P}{}^QE_Q+E_{N]P}]\nn\\
&=g_{MN}{}^Q(-j_{QP}{}^RE_R+E_{QP})
     +2j_{[N|P|}{}^Q(-j_{M]Q}{}^RE_R+E_{M]Q})\nn\\
 &\qquad-2(-j_{[MN]}{}^QE_{QP}-j_{[M|P|}{}^QE_{N]Q})\\
&=([j_M,j_N]-g_{MN}{}^Qj_Q)_P{}^RE_R
     +(g_{MN}{}^Q+2j_{[MN]}{}^Q)E_{PQ}\nn\\
 &=0\;.\nn
\end{align}

It then remains to determine the brackets between level $1$ and level $-1$ generators, \ie, $[R_1,R_{-1}]$. Due to the graded antisymmetry of the structure constants, these brackets contain the same structure constants as $[R_{-1},R_{-1}]$, and are thus symmetric.
The calculation can be systematised as follows.
Introduce a sequence of basis elements
$E_\MP=E_{M_1\ldots M_p}$ for various
values of $p$, and let them transform covariantly (\ie, with
$t_\alpha$ and $j_M$) apart from an extra term:
\be
[J_M,E_{M_1\ldots M_p}]=-(j_M\cdot E)_{M_1\ldots M_P}+E_{MM_1\ldots
M_p}\;,
\label{JShift}
\ee
\ie, $[J_M,E_\NP]=-(j_M\cdot E)_\NP+E_{M\NP}$.
It is in fact possible to begin already at $p=0$.
We then write
\begin{align}
[E_\MP,U^\alpha]&=a_\MP{}^{\alpha\beta}T_\beta
+b_\MP{}^{\alpha N}J_N\;,\nn\\
[E_\MP,V^N]&=b_\MP{}^{\alpha N}T_\alpha
+c_\MP{}^{NP}J_P\;.
\label{ABCDefinitions}
\end{align}
for some coefficients $a,b,c$, which are to be determined.
The double occurrence of $b$ follows from the graded antisymmetry of the structure constants, which also implies that $a$ and $c$ are symmetric in upper indices.

The Jacobi identities with $T$ are equivalent to
\begin{align}
&\Delta_\alpha a_\MP{}^{\beta\gamma}
      +2u_{N\alpha}{}^{(\beta}b_\MP{}^{\gamma) N}=0\;,\nn\\
&\Delta_\alpha b_\MP{}^{\beta N}
      +u_{P\alpha}{}^{\beta}c_\MP{}^{NP}=0\;,\label{TJacobis}\\
&\Delta_\alpha c_\MP{}^{NP}=0\;.\nn
\end{align}
The Jacobi identities with $J$ relate coefficients at different $p$:
\begin{align}
a_{M\MP}{}^{\alpha\beta}&=\Delta_Ma_\MP{}^{\alpha\beta}\;,\nn\\
b_{M\MP}{}^{\alpha N}&=\Delta_Mb_\MP{}^{\alpha\dot N}
+t_{\beta M}{}^Na_\MP{}^{\alpha\beta}\;,\label{JJacobis}\\
c_{M\MP}{}^{NP}&=\Delta_Mc_\MP{}^{\dot N\dot P}
+2t_{\alpha M}{}^{(N}b_\MP{}^{|\alpha|P)}\;.\nn
\end{align}
Here, one must remember that lower indices belong to the fundamental and
transform with $j$, while upper ones belong to the adjoint and
transform with $g$. The latter are dotted when $\Delta_M$ acts, in accordance with
eq. \eqref{DeltaDefinitions}.

We can then use induction and show that if eq. \eqref{TJacobis} is
satisfied at $p$, then the coefficients $a,b,c$ at $p+1$, given by
eq. \eqref{JJacobis} also satisfy eq. \eqref{TJacobis}.
This is done by commuting $\Delta_\alpha$ and $\Delta_M$, as in Section
\ref{CovarianceSection}.
Checking the first equation in \eqref{TJacobis} at $p+1$:
\begin{align}
\Delta_\alpha a_{M\MP}{}^{\beta\gamma}
     &=\Delta_\alpha\Delta_Ma_\MP{}^{\beta\gamma}\nn\\
&=\Delta_M\Delta_\alpha a_\MP{}^{\beta\gamma}
     +[\Delta_\alpha,\Delta_M]a_\MP{}^{\beta\gamma}\nn\\
&=-2u_{N\alpha}{}^{(\beta}\Delta_Mb_\MP{}^{\gamma)N}
     -2t_{\delta M}{}^Nu_{N\alpha}{}^{(\beta}a_\MP{}^{\gamma)\delta}\\
&=-2u_{N\alpha}{}^{(\beta}b_{M\MP}{}^{\gamma)N}\;.\nn
\end{align}
In passing, we used $\Delta_Mu_{\dot N\alpha}{}^\beta=0$ (where the $N$
index transforms with $g$).
Analogous calculations hold for $b$ and $c$.

Starting from $p=0$ with 
\begin{align}
a^{\alpha\beta}&=\eta^{\alpha\beta}\;,\nn\\
b^{\alpha M}&=0\;,\\
c^{MN}&=0\;,\nn
\end{align}
where $\eta$ is the inverse Killing metric on $\fg^{++}$,
this gives a constructive recipe to determine the whole set of coefficients recursively.
The concrete expressions obtained (which however are less useful than the recursion relations themselves) become, after using various identities among eqs. \eqref{ffid}--\eqref{ggid}, 
\begin{align}
a_M{}^{\alpha\beta}&=-2u_M{}^{(\alpha\beta)}\;,\nn\\
b_M{}^{\alpha N}&=t^\alpha{}_M{}^N\;,\nn\\
c_M{}^{NP}&=0\;;\nn\\
a_{MN}{}^{\alpha\beta}&=2(u_{(M}{}^{\gamma(\alpha}+u_{(M}{}^{(\alpha|\gamma|})
        u_{N)\gamma}{}^{\beta)}-2j_{(MN)}{}^Pu_P{}^{(\alpha\beta)}\;,\\
b_{MN}{}^{\alpha P}&=-4t_{\beta(M}{}^Pu_{N)}{}^{(\alpha\beta)}
        +2j_{Q(M}{}^Pt^\alpha{}_{N)}{}^Q\;,\nn\\
 c_{MN}{}^{PQ}&=2t^\alpha{}_M{}^{(P}t_{\alpha N}{}^{Q)}\;.\nn
\end{align}
We always lower and raise adjoint indices with the $\fg^{++}$ Killing metric and its inverse, but need to remember its non-invariance under $\Delta_M$.

Note that $c$ appears first at $p=2$, and that
$c_{MN}{}^{PQ}$ is proportional to the symmetrised $Y$-tensor \cite{Cederwall:2017fjm}, annihilating
$S_2^\compl$. 
In order for $S_2^\compl$ to be an ideal at level $1$, we also need $a_{MN}{}^{\alpha\beta}$ and $b_{MN}{}^{\alpha P}$ to contain only $S_2$ in the lower indices. 
By using this property for $c_{MN}{}^{PQ}$, we trace it back to the coefficients $b$ and $a$ using
eq. \eqref{TJacobis}. Since there are no $\fg^{++}$-invariant tensors $a^{(0)}$ and $b^{(0)}$ with
$\Delta_\alpha a^{(0)}_{MN}{}^{\alpha\beta}=0$, $\Delta_\alpha b^{(0)}_{MN}{}^{\alpha P}=0$, the desired result follows.

The module we are considering is of course obtained by
discarding $p=0$ and starting at $p=1$.

One criterion, from the Jacobi identity $[J,J,E]$, must be
$E_{[M_1M_2]M_3\ldots M_p}=0$, which needs to be automatically
satisfied by the recursively obtained $a$, $b$ and $c$. This follows (after some calculation)
from using $[\Delta_M,\Delta_N]=0$ in the recursion.

One can conclude, just from the observation that $S_2^\compl$ is absent
at $p=2$, that the modules at higher $p$ are empty.
If we think of $\Delta_M$ as a ``commuting bosonic object'' in $R(\lambda)$, it
is constrained (at least in the recursion relations involving $a$, $b$ and $c$) 
so that $\Delta^2$ has no part in the leading symmetric
module $R(2\lambda)$. This is the ``opposite'' to an object in a
minimal orbit. Such a constraint is strong enough to kill
all degrees of freedom in $\Delta$. We would have liked to show that already $\Delta^3$
is empty by this condition, \ie, that $R(2\lambda)\otimes R(\lambda)$
contains all of $\vee^3R(\lambda)$. This is certainly true for modules
of finite-dimensional simple Lie algebras. We have checked it in examples, \eg\ $E_{10}$, where it holds (with great and increasing margin) to degree $7$ in a $\gl(10)$ grading with respect to the exceptional node, and also to degree $10$ in the $\gl(3)$ grading of $A_1^{++}$. Presumably it can be proven, but it remains an assumption.

To summarise the local superalgebra, we have at levels $0,\pm1$:
\begin{align}
R_{-1}&=\bar\AA=\langle V^M\rangle\oplus\langle U^\alpha\rangle\;,\nn\\
R_0&=\AA=\langle T_\alpha\rangle\oplus\langle J_M\rangle\;,\\
R_1&=\langle E_M\rangle\oplus\langle E_{MN}\rangle\;.\nn
\end{align}
The brackets in the local superalgebra are
\begin{align}
[T_\alpha,T_\beta]&=f_{\alpha\beta}{}^\gamma T_\gamma\;,\nn\\
[T_\alpha,J_M]&=-t_{\alpha M}{}^NJ_N+u_{M\alpha}{}^\beta T_\beta\;,\\
[J_M,J_N]&=g_{MN}{}^PJ_P\;;\nn
\end{align}
\begin{align}
[T_\alpha,E_M]&=-t_{\alpha M}{}^NE_N\;,\nn\\
[T_\alpha,E_{MN}]&=-2t_{\alpha(M}{}^PE_{N)P}\;,\nn\\
[J_M,E_N]&=-j_{MN}{}^PE_P+E_{MN}\;,\\
[J_M,E_{NP}]&=-2j_{M(N}{}^QE_{P)Q}\;;\nn
\end{align}
\begin{align}
[T_\alpha,V^M]&=t_{\alpha N}{}^MV^N\;,\nn\\
[T_\alpha,U^\beta]&=-f_{\alpha\gamma}{}^\beta U^\gamma
               -u_{M\alpha}{}^\beta V^M\;,\nn\\
[J_M,V^N]&=-g_{MP}{}^NV^P-t_{\alpha M}{}^NU^\alpha\;,\\
[J_M,U^\alpha]&=u_{M\beta}{}^\alpha U^\beta\;;\nn
\end{align}
\begin{align}
[E_M,V^N]&=t^\alpha{}_M{}^NT_\alpha\;,\nn\\
[E_M,U^\alpha]&=-2u_M{}^{(\alpha\beta)}T_\beta+t^\alpha{}_M{}^NJ_N\;,\nn\\
[E_{MN},V^P]&=b_{MN}{}^{\alpha P}T_\alpha+c_{MN}{}^{PQ}J_Q\;,\\
[E_{MN},U^\alpha]&=a_{MN}{}^{\alpha\beta}T_\beta+b_{MN}{}^{\alpha
P}J_P\;.\nn
\end{align}

In fact, due to the existence of the invariant bilinear form, we directly know also $R_{-2}=\bar R_1$, its transformations under $R_0=\AA$, and the brackets $[R_1,R_{-2}]$, which contain the same structure constants as $[R_1,R_0]$.

The tensor hierarchy algebra $S(\fg^{++})$ is the superalgebra generated by the local superalgebra, modulo the maximal ideal intersecting it trivially.

\begin{figure}
  \begin{center}
\includegraphics[scale=.75]{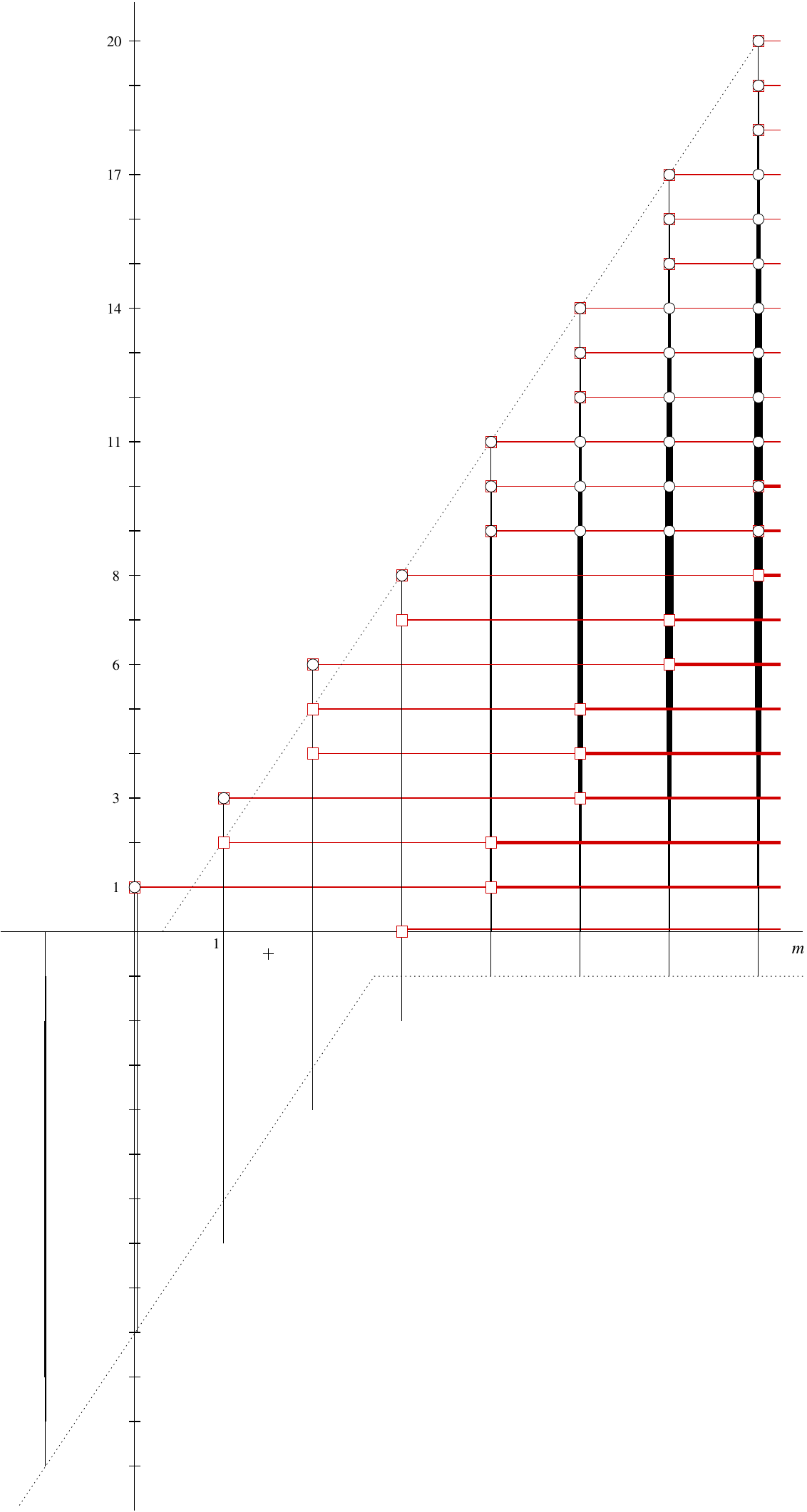}
\caption{\it The double
grading of $S(E_{10})$, with 
respect to the fermionic and the exceptional nodes. Each grade $(\ell,m)$ contains a
$\gl(10)$ module. These modules organise into ``columns'', containing
tensor products of some $\gl(10)$ module with all forms. The collection
of such columns at given degree $m$ with respect to the exceptional node forms a
module of the tensor hierarchy algebra $W(10)$ at degree $0$. Tops of
columns are marked with a black circle. Horizontally, the $\gl(10)$
modules are organised in $E_{10}$ modules. The lowest states in
$\BB_\ell$ and $\BB_{\ell+1}$ at level $\ell$ are denoted with red squares.}    
\label{SE10Figure}
\end{center}
\end{figure}

\begin{figure}
\begin{center}
\includegraphics{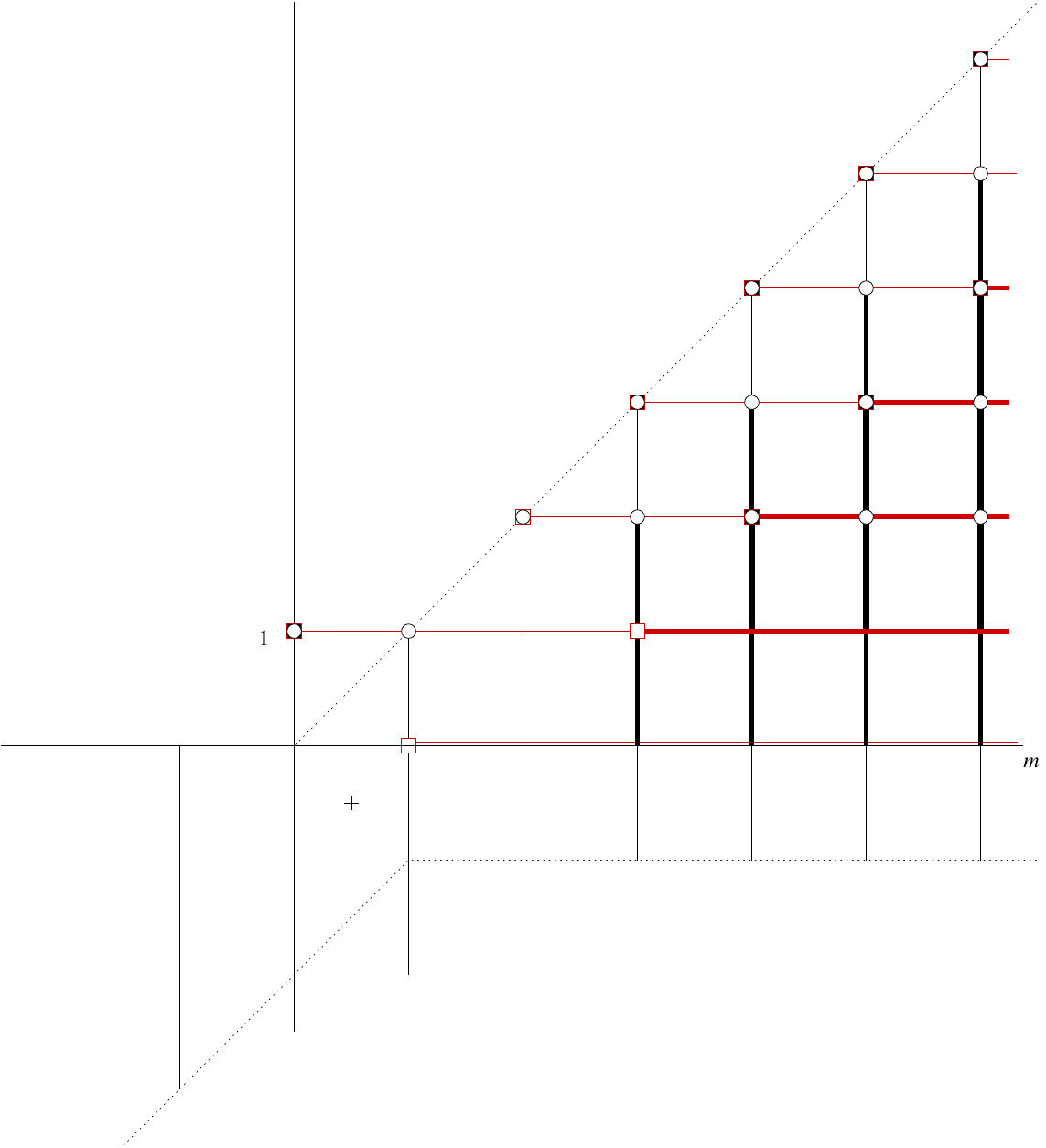}
\caption{\it The double grading of $S(A_1^{++})$, with
respect to the fermionic and the rightmost node. Each grade $(\ell,m)$ contains a
$\gl(3)$ module. These modules organise into ``columns'', containing
tensor products of some $\gl(3)$ module with all forms. The collection
of such columns at given degree $m$ with respect to node $1$ forms a
module of the tensor hierarchy algebra $W(3)$ at degree $0$. Tops of
columns are marked with a black circle. Horizontally, the $\gl(3)$
modules are organised in $A_1^{++}$ modules. The lowest states in
$\BB_\ell$ and $\BB_{\ell+1}$ at level $\ell$ are denoted with red squares.
The superalgebra on the diagonal is freely generated. Our computation goes to $m=10$.}    
\label{SA1++Figure}
\end{center}
\end{figure}

\def\sw#1{\baselineskip=6pt$\scriptscriptstyle#1$}
\def\mw#1{\baselineskip=6pt$\scriptstyle#1$}

\renewcommand{\arraystretch}{.8}	

\begin{table}
\vspace{-1.5cm}
\begin{center}
\begin{tabular}{r | m {1.4 cm}| m {1.4 cm} | m {1.4 cm} | m {1.4 cm} | m {1.4 cm} | m {1.4 cm} | m {1.4cm} | m {1.4 cm} |}
&\mw{m=0}&\mw{1}&\mw2&\mw3&\mw4&\mw5&\mw6&\mw7\\ \hline
\mw{\ell=20}&&&&&&&&\sw{(000000001)}\\ \hline
\mw{19}&&&&&&&&\sw{2(000000002)} \sw{2(000000010)}\\ \hline
\mw{18}&&&&&&&&\sw{3(000000003)} \sw{5(000000011)} \sw{3(000000100)}\\ \hline
\mw{17}&&&&&&
&\sw{(000000001)}
&\sw{2(000000004)} \sw{5(000000012)} \sw{3(000000020)}
         \sw{7(000000101)} \sw{4(000001000)}\\ \hline
\mw{16}&&&&&&
&\sw{(000000002)} \sw{2(000000010)}
&\sw{(000000005)} \sw{2(000000013)} \sw{2(000000021)}
         \sw{6(000000102)} \sw{5(000000110)} \sw{9(000001001)} \sw{4(000010000)}\\ \hline
\mw{15}&&&&&&
&\sw{2(000000003)} \sw{3(000000011)} \sw{(000000100)}
&\sw{2(000000103)} \sw{3(000000111)} \sw{2(000000200)}
         \sw{7(000001002)}  \sw{7(000001010)} \sw{8(000010001)} \sw{4(000100000)} \\ \hline
\mw{14}&&&&&
&\sw{(000000001)}
&\sw{(000000004)} \sw{(000000012)} \sw{(000000020)} \sw{3(000000101)} \sw{3(000001000)} 
&\sw{(000000201)} \sw{2(000001003)} \sw{4(000001011)}
         \sw{4(000001100)}  \sw{5(000010002)} \sw{6(000010010)} 
         \sw{8(000100001)} \sw{4(001000000)} \\ \hline
\mw{13}&&&&&
&\sw{(000000002)} \sw{(000000010)}
&\sw{2(000000102)} \sw{(000000110)} \sw{3(000001001)} \sw{	(000010000)} 
&\sw{2(000001101)} \sw{2(000002000)} \sw{(000010003)}
         \sw{3(000010011)}  \sw{3(000010100)} \sw{5(000100002)} 
         \sw{5(000100010)} \sw{7(001000001)} \sw{3(010000000)} \\ \hline
\mw{12}&&&&&
&\sw{(000000003)} \sw{(000000011)} \sw{	(000000100)}
&\sw{(000001002)} \sw{2(000001010)} \sw{2(000010001)} \sw{2(000100000)} 
&\sw{(000002001)} \sw{(000010101)} \sw{(000011000)}
         \sw{(000100003)}  \sw{2(000100011)} \sw{3(000100100)} 
         \sw{3(001000002)} \sw{5(001000010)} \sw{5(010000001)} 
         \sw{3(100000000)} \\ \hline
\mw{11}&&&&
&\sw{(000000001)}
&\sw{(000000101)} \sw{(000001000)}
&\sw{(000001100)} \sw{(000010002)} \sw{2(000100001)} \sw{(001000000} 
&\sw{(000000000)} \sw{(000100101)} \sw{(000101000)}
         \sw{(001000011)} \sw{2(001000100)} \sw{3(010000002)} 
         \sw{3(010000010)} \sw{4(100000001)}  \\ \hline
\mw{10}&&&&
&\sw{(000000002)}
&\sw{(000001001)}
&\sw{(000100010)} \sw{(001000001)} \sw{(010000000)}
&\sw{2(000000001)} \sw{(001001000)} \sw{(010000011)}
         \sw{(010000100)} \sw{(100000002)} \sw{2(100000010)}  \\ \hline
\mw9&&&&
&\sw{(000000100)}
&\sw{(000100000)}
&\sw{(010000001)}
&\sw{(000000002)} \sw{(100000100)}  \\ \hline
\mw8&&&&\sw{(000000001)}&&&& \\ \hline
\mw7&&&&&&&& \\ \hline
\mw6&&&\sw{(000000000)}&&&&& \\ \hline
\mw5&&&&&&&& \\ \hline
\mw4&&&&&&&& \\ \hline
\mw3&&\sw{(000000000)}&&&&&& \\ \hline
\mw2&&&&&&&& \\ \hline
\mw1&\sw{(100000000)}&&&&&&& \\ \hline
\end{tabular}
\captionof{table}{\it The $\sl(10)$ modules $r_{\ell,m}$ for the columns in $S(E_{10})$. The 
$\gl(10)$ weights, normalised so that a $1$-form has weight $1$, are $3m-\ell$.}
\label{SE10ColumnsTable}
\end{center}
\end{table}

\begin{table}
\vspace{-2cm}
\begin{center}
\begin{tabular}{r | m {.4cm}| m {.4cm} | m {.4cm} | m {.4cm} | m {.5cm} | m {.5cm} | m {.6cm} | m {.6cm} | m {1.8cm} | m {1.95cm} | m {2.2cm} |}
\mw{m=}&\mw{0}&\mw{1}&\mw2&\mw3&\mw4&\mw5&\mw6&\mw7&\mw8&\mw9&\mw{10}\\ \hline
\mw{\ell=10}&&&&&&&&&&&\sw{20(01)} \sw{24(04)} \sw{3(07)} \sw{44(12)} \sw{16(15)} \sw{(18)} \sw{27(20)} \sw{31(23)} \sw{4(26)} \sw{29(31)} \sw{7(34)} \sw{10(42)} \sw{3(50)}\\ \hline
\mw{9}&&&&&&&&&
&\sw{4(00)} \sw{13(03)} \sw{3(06)} \sw{19(11)} \sw{12(14)} \sw{(17)} \sw{18(22)} \sw{3(25)} \sw{9(30)} \sw{5(33)} \sw{5(41)} 
&\sw{273(02)} \sw{142(05)} \sw{11(08)} \sw{191(10)} \sw{399(13)} \sw{75(16)} \sw{(19)} \sw{413(21)} \sw{211(24)} \sw{10(27)} \sw{287(32)} \sw{37(35)} \sw{143(40)} \sw{65(43)} \sw{57(51)} \\ \hline
\mw{8}&&&&&&&&
&\sw{6(02)} \sw{3(05)} \sw{6(10)} \sw{8(13)} \sw{(16)} \sw{9(21)} \sw{2(24)} \sw{4(32)} \sw{(40)} 
&\sw{84(01)} \sw{85(04)} \sw{9(07)} \sw{180(12)} \sw{54(15)} \sw{(18)} \sw{103(20)} \sw{123(23)} \sw{9(26)} \sw{120(31)} \sw{27(34)} \sw{37(42)} \sw{18(50)} 
&\sw{341(00)} \sw{1186(03)} \sw{336(06)} \sw{10(09)} \sw{1579(11)} \sw{1285(14)} \sw{125(17)} \sw{1955(22)} \sw{510(25)} \sw{9(28)} \sw{1020(30)} \sw{982(33)} \sw{64(36)} \sw{892(41)} \sw{165(44)} \sw{225(52)} \sw{104(60)} \\ \hline
\mw{7}&&&&&&&
&\sw{3(01)} \sw{2(04)} \sw{5(12)} \sw{(15)} \sw{3(20)} \sw{2(23)} \sw{2(31)} 
&\sw{15(00)} \sw{46(03)} \sw{7(06)} \sw{65(11)} \sw{37(14)} \sw{(17)} \sw{63(22)} \sw{8(25)} \sw{35(30)} \sw{18(33)} \sw{18(41)} 
&\sw{401(02)} \sw{195(05)} \sw{8(08)} \sw{290(10)} \sw{580(13)} \sw{88(16)} \sw{623(21)} \sw{290(24)} \sw{8(27)} \sw{432(32)} \sw{44(35)} \sw{220(40)} \sw{95(43)} \hfill\sw{90(51)} 
&\sw{1978(01)} \sw{2406(04)} \sw{387(07)} \sw{2(0,10)} \sw{4787(12)} \sw{1980(15)} \sw{91(18)} \sw{2652(20)} \sw{4194(23)} \sw{578(26)} \sw{3(29)} \sw{3954(31)} \sw{1574(34)} \sw{46(37)} \sw{2188(42)} \sw{198(45)} \sw{1043(50)} \sw{383(53)} \sw{348(61)} \\ \hline
\mw{6}
&&&&&&
&\sw{(00)} \sw{(03)} \sw{3(11)} \sw{(14)} \sw{2(22)} 
&\sw{20(02)} \sw{6(05)} \sw{15(10)} \sw{23(13)} \sw{(16)} \sw{27(21)} \sw{6(24)} \sw{11(32)} \sw{6(40)} 
&\sw{102(01)} \sw{101(04)} \sw{6(07)} \sw{217(12)} \sw{57(15)} \sw{132(20)} \sw{146(23)} \sw{7(26)} \sw{152(31)} \sw{28(34)} \sw{49(42)} \hfill\sw{23(50)} 
&\sw{244(00)} \sw{815(03)} \sw{202(06)} \sw{2(09)} \sw{1118(11)} \sw{849(14)} \sw{60(17)} \sw{1370(22)} \sw{310(25)} \sw{2(28)} \sw{721(30)} \sw{667(33)} \sw{30(36)} \sw{638(41)} \sw{106(44)} \sw{158(52)} \sw{80(60)} 
&\sw{4283(02)} \sw{2576(05)} \sw{214(08)} \sw{2954(10)} \sw{7035(13)} \sw{1581(16)} \sw{27(19)} \sw{7221(21)} \sw{4594(24)} \sw{322(27)} \sw{6384(32)} \sw{1285(35)} \sw{13(38)} \sw{3199(40)} \sw{2458(43)} \sw{112(46)} \sw{2223(51)} \sw{322(54)} \sw{458(62)} \hfill \sw{220(70)}\\ \hline
\mw{5}&&&&&
&\sw{(02)} \sw{(10)} \sw{(13)} \sw{(21)} 
&\sw{6(01)} \sw{5(04)} \sw{12(12)} \sw{(15)} \sw{8(20)} \sw{4(23)} \sw{6(31)} 
&\sw{16(00)} \sw{41(03)} \sw{5(06)} \sw{63(11)} \sw{33(14)} \sw{62(22)} \sw{5(25)} \sw{34(30)} \sw{17(33)} \sw{19(41)} 
&\sw{220(02)} \sw{90(05)} \sw{2(08)} \sw{155(10)} \sw{302(13)} \sw{35(16)} \sw{336(21)} \sw{144(24)} \sw{(27)} \sw{231(32)} \sw{18(35)} \sw{124(40)} \sw{48(43)} \hfill\sw{51(51)} 
&\sw{691(01)} \sw{801(04)} \sw{97(07)} \sw{1648(12)} \sw{603(15)} \sw{15(18)} \sw{939(20)} \sw{1392(23)} \sw{151(26)} \sw{1374(31)} \sw{488(34)} \sw{8(37)} \sw{721(42)} \sw{52(45)} \sw{366(50)} \sw{122(53)} \sw{122(61)} 
&\sw{1220(00)} \sw{4456(03)} \sw{1462(06)} \sw{48(09)} \sw{5856(11)} \sw{5376(14)} \sw{633(17)} \sw{2(1,10)} \sw{8198(22)} \sw{2602(25)} \sw{79(28)} \sw{4198(30)} \sw{5117(33)} \sw{516(36)} \sw{(39)} \sw{4687(41)} \sw{1443(44)} \sw{25(47)} \sw{1974(52)} \sw{131(55)} \sw{981(60)} \sw{278(63)} \sw{265(71)} \\ \hline
\mw{4}&&&&
&\sw{(01)} \sw{(12)}
&\sw{(00)} \sw{3(03)} \sw{5(11)} \sw{(14)} \sw{3(22)} \sw{2(30)} 
&\sw{12(02)} \sw{4(05)} \sw{12(10)} \sw{16(13)} \sw{19(21)} \sw{3(24)} \sw{9(32)} \sw{4(40)} 
&\sw{39(01)} \sw{34(04)} \sw{(07)} \sw{82(12)} \sw{17(15)} \sw{49(20)} \sw{52(23)} \sw{(26)} \sw{59(31)} \sw{9(34)} \sw{18(42)} \sw{10(50)} 
&\sw{58(00)} \sw{191(03)} \sw{35(06)} \sw{268(11)} \sw{184(14)} \sw{7(17)} \sw{315(22)} \sw{57(25)} \sw{178(30)} \sw{146(33)} \sw{4(36)} \sw{152(41)} \sw{19(44)} \sw{37(52)} \sw{18(60)} 
&\sw{703(02)} \sw{371(05)} \sw{18(08)} \sw{496(10)} \sw{1105(13)} \sw{194(16)} \sw{(19)} \sw{1185(21)} \sw{663(24)} \sw{27(27)} \sw{1007(32)} \sw{157(35)} \sw{521(40)} \sw{357(43)} \sw{9(46)} \sw{350(51)} \sw{39(54)} \sw{66(62)} \sw{35(70)} 
&\sw{1793(01)} \sw{2320(04)} \sw{396(07)} \sw{3(0,10)} \sw{4575(12)} \sw{2060(15)} \sw{106(18)} \sw{	2522(20)} \sw{4417(23)} \sw{700(26)} \sw{5(29} \sw{4203(31)} \sw{2001(34)} \sw{83(37)} \sw{2770(42)} \sw{386(45)} \sw{(48)} \sw{1393(50)} \sw{799(53)} \sw{20(56)} \sw{754(61)} \sw{72(64)} \sw{118(72)} \hfill\sw{58(80)}\\ \hline
\mw{3}&&&
&\sw{(11)}
&\sw{(02)} \sw{(10)} \sw{(13)} \sw{2(21)} 
&\sw{3(01)} \sw{2(04)} \sw{5(12)} \sw{3(20)} \sw{2(23)} \sw{3(31)} 
&\sw{2(00)} \sw{9(03)} \sw{15(11)} \sw{6(14)} \sw{13(22)} \sw{(25)} \sw{8(30)} \sw{3(33)} \sw{5(41)} 
&\sw{29(02)} \sw{10(05)} \sw{23(10)} \sw{39(13)} \sw{2(16)} \sw{46(21)} \sw{15(24)} \sw{30(32)} \sw{(35)} \sw{17(40)} \sw{5(43)} \sw{7(51)} 
&\sw{67(01)} \sw{64(04)} \sw{5(07)} \sw{151(12)} \sw{42(15)} \sw{86(20)} \sw{117(23)} \sw{6(26)} \sw{125(31)} \sw{33(34)} \sw{58(42)} \sw{2(45)} \sw{33(50)} \sw{8(53)} \sw{10(61)} 
&\sw{81(00)} \sw{287(03)} \sw{70(06)} \sw{395(11)} \sw{315(14)} \sw{21(17)} \sw{520(22)} \sw{125(25)} \sw{(28)} \sw{280(30)} \sw{294(33)} \sw{15(36)} \sw{296(41)} \sw{66(44)} \sw{110(52)} \sw{3(55)} \sw{58(60)} \sw{11(63)} \sw{14(71)} 
&\sw{885(02)} \sw{546(05)} \sw{37(08)} \sw{622(10)} \sw{1516(13)} \sw{334(16)} \sw{4(19)} \sw{1581(21)} \sw{1035(24)} \sw{68(27)} \sw{1516(32)} \sw{319(35)} \sw{2(38)} \sw{763(40)} \sw{676(43)} \sw{32(46)} \sw{641(51)} \sw{122(54)} \sw{192(62)} \sw{5(65)} \sw{103(70)} \sw{16(73)} \sw{19(81)} \\ \hline
\mw{2}&&
&\sw{(02)}
&\sw{(12)}
&\sw{(11)} \sw{(22)}
&\sw{(02)} \sw{(10)} \sw{(13)} \sw{2(21)} \sw{(32)} 
&\sw{2(01)} \sw{2(04)} \sw{4(12)} \sw{3(20)} \sw{2(23)} \sw{3(31)} \sw{(42)} 
&\sw{(00)} \sw{6(03)} \sw{10(11)} \sw{5(14)} \sw{10(22)} \sw{(25)} \sw{6(30)} \sw{3(33)} \sw{5(41)} \sw{(52)} 
&\sw{18(02)} \sw{6(05)} \sw{12(10)} \sw{24(13)} \sw{2(16)} \sw{28(21)} \sw{12(24)} \sw{21(32)} \sw{(35)} \sw{13(40)} \sw{5(43)} \sw{7(51)} \sw{(62)} 
&\sw{34(01)}  \sw{35(04)}  \sw{4(07)}  \sw{81(12)}  \sw{26(15)}  \sw{45(20)}  \sw{69(23)}  \sw{5(26)}  \sw{72(31)}  \sw{24(34)}  \sw{40(42)}  \sw{2(45)}  \sw{23(50)}  \sw{8(53)}  \sw{10(61)} \hfill \sw{(72)}  
&\sw{38(00)} \sw{136(03)} \sw{40(06)} \sw{187(11)} \sw{163(14)} \sw{14(17)} \sw{265(22)} \sw{73(25)} \sw{(28)} \sw{136(30)} \sw{165(33)} \sw{11(36)} \sw{163(41)} \sw{47(44)} \sw{73(52)} \sw{3(55)} \sw{40(60)} \sw{11(63)} \sw{14(71)} \hfill\sw{(82)} \\ \hline
\mw{1}&\sw{(10)}&\sw{(01)}&&&&&&&&&\\ \hline
\end{tabular}
\captionof{table}{\it The $\sl(3)$ modules $r_{\ell,m}$ for the columns in $S(A_1^{++})$. The 
$\gl(3)$ weights, normalised so that a $1$-form has weight $1$, are $2m-\ell$.}
\label{SA1++ColumnsTable}
\end{center}
\end{table}

\renewcommand{\arraystretch}{1.3}

When $(\lambda,\lambda)=0$, the standard definition of $W(\fg^{++})$ does not produce a simple superalgebra, roughly because $\BB(\fg^{++})$ is already a subalgebra of $S(\fg^{++})$. 
A meaningful definition of $W(\fg^{++})$ would instead be as a subalgebra of $S(\fg^{+++})$. This will be dealt with in a forthcoming paper \cite{CederwallPalmkvistForth}.

\subsection{A conjecture}

As we have shown, $S(\fg^{++})$ at levels $\ell=-2,-1,0,1$ contains
\begin{align}
\ell=-2\,:&\quad \bar S_2\oplus R(\lambda)\;,\nn\\
-1\,:&\quad R(\lambda)\oplus\adj\;,\nn\\
0\,:&\quad\adj\oplus R(-\lambda)\;,\\
1\,:&\quad R(-\lambda)\oplus S_2\;,\nn
\end{align}
where $\adj$ is the adjoint of $\fg^{++}$.
This is precisely $\BB_\ell\oplus\BB_{\ell+1}$.
We would like to conjecture that this vector space (not algebra) decomposition holds to all levels, \ie, that 
\be
S(\fg^{++})=\BB(\fg^{++})\oplus\BB(\fg^{++})[1]\;.
\ee

Consider
$R_2$, obtained 
as $[R_1,R_1]$. Let
\be
Z_{\MP,\NP}=[E_\MP,E_\NP]\;.
\ee
Then,
\begin{align}
[T_\alpha,Z_{\MP,\NP}]&=-(t_\alpha\cdot Z)_{\MP,\NP}\;,\nn\\
[J_M,Z_{\NP,\PP}]&=-(j_M\cdot Z)_{\NP,\PP}+Z_{M\NP,\PP}+Z_{\NP,M\PP}\;.
\end{align}
We can investigate the appearance of ideals by checking the Jacobi identities
with $U^\alpha$ and $V^M$. One directly obtains
\begin{align}
[Z_{\MP,\NP},U^\alpha]&=A_{\MP,\NP}{}^{\alpha\PP}E_\PP\;,\nn\\
[Z_{\MP,\NP},V^P]&=B_{\MP,\NP}{}^{P\QP}E_\QP\;,
\end{align}
where
\begin{align}
A_{\MP,\NP}{}^{\alpha\PP}&=t_{\beta\NP}{}^\PP a_\MP{}^{\alpha\beta}
    +j_{Q\NP}{}^\PP b_\MP{}^{\alpha Q}\nn\\
    &\qquad-\delta_{Q\NP}^\PP b_\MP{}^{\alpha Q}+(\MP\leftrightarrow\NP)\;,\nn\\
B_{\MP,\NP}{}^{P\QP}&=t_{\alpha\NP}{}^\QP b_\MP{}^{\alpha P}
    +j_{R\NP}{}^\QP c_\MP{}^{PR}\\\
    &\qquad-\delta_{R\NP}^\QP c_\MP{}^{PR}
    +(\MP\leftrightarrow\NP)\;\nn
\end{align}
(The notation with Kronecker deltas is shorthand; projectors on $S_2$
should be inserted when appropriate).

In order to systematise the calculation, we introduce the coefficients
\begin{align}
v_{\MP,N}{}^{P\alpha}&=t_{\beta N}{}^Pa_\MP{}^{\alpha\beta}
             +j_{QN}{}^Pb_\MP{}^{\alpha Q}\;,\nn\\
w_{\MP,N}{}^{PQ}&=t_{\beta N}{}^Pb_\MP{}^{\beta Q}
             +j_{RN}{}^Pc_\MP{}^{QR}\;.\label{VWDefs}
\end{align}
It is straightforward to show that
\begin{align}
0&=\Delta_\alpha v_{\MP,N}{}^{P\beta}
        +u_{Q\alpha}{}^\beta w_{\MP N}{}^{PQ}\;,\nn\\
0&=\Delta_\alpha w_{\MP,N}{}^{PQ}\;,\nn\\
v_{M\NP,P}{}^{Q\alpha}&=\Delta_M v_{\NP,P}{}^{Q\alpha}\;,\\
w_{M\NP,P}{}^{QR}&=\Delta_M w_{\NP,P}{}^{Q\dot R}
         +t_{\alpha M}{}^Rv_{\NP,P}{}^{Q\alpha} \;.\nn
\end{align}
We also observe that symmetrisation in the part with two lower indices
yields
\begin{align}
v_{(M,N)}{}^{P\alpha}&=\fr2b_{MN}{}^{P\alpha}\;,\nn\\
w_{(M,N)}{}^{PQ}&=\fr2c_{MN}{}^{PQ}\;.\label{vbwcrel}
\end{align}
Concerning the 3-index part, a short calculation gives
\begin{align}
v_{M(N,P)}{}^{Q\alpha}&=\Delta_Mv_{(N,P)}{}^{Q\alpha}
       =\fr2\Delta_Mb_{NP}{}^{\alpha Q}\nn\\
&=\fr2(\Delta_Mb_{NP}{}^{\alpha\dot Q}
       -(g_{MR}{}^Q+j_{MR}{}^Q)b_{NP}{}^{\alpha R})\nn\\
&=\fr2(-t_{\beta M}{}^Qa_{NP}{}^{\alpha\beta}
       -j_{RM}{}^Qb_{NP}{}^{\alpha R})\label{Threeindexvw}\\
&=-\fr2v_{NP,M}{}^{Q\alpha}\;,\nn
\end{align}
and, by an analogous calculation, $w_{M(N,P)}{}^{QR}=-\fr2w_{NP,M}{}^{QR}$.
This shows that the completely symmetric parts vanish,
$v_{(MN,P)}{}^{Q\alpha}=0=w_{(MN,P)}{}^{QR}$.

Consider first $Z_{M,N}$. The relevant coefficients are
\begin{align}
A_{M,N}{}^{\alpha P}&=2(t_{\beta(M}{}^Pa_{N)}{}^{\alpha\beta}
                   +j_{Q(M}{}^Pb_{N)}{}^{\alpha Q})
                   =2v_{(M,N)}{}^{P\alpha}=b_{MN}{}^{\alpha P}\;,\nn\\
A_{M,N}{}^{\alpha,PQ}&=-2t^\alpha_{(M}{}^{(P}\delta_{N)}^{Q)}\;,\nn\\
B_{M,N}{}^{P,Q}&=2w_{(M,N)}{}^{PQ}=c_{MN}{}^{PQ}\;,\\
B_{M,N}{}^{P,QR}&=0\;,\nn
\end{align}
This shows that $S_2^\compl$ is an ideal in $Z_{M,N}=[E_M,E_N]$. This is of course completely expected, since it is an ideal for the Borcherds superalgebra.

Next, investigate $Z_{MN,P}$. The coefficients become
\begin{align}
A_{MN,P}{}^{\alpha Q}&=v_{MN,P}\;,\nn\\
A_{MN,P}{}^{\alpha QR}&=2v_{P,(M}{}^{\alpha Q}\delta_{N)}^{R)}
            -b_{MN}{}^{\alpha(Q}\delta_P^{R)}\;,\nn\\
B_{MN,P}{}^{Q,R}&=w_{MN,P}{}^{QR}\;,\\
B_{MN,P}{}^{Q,RS}&=2w_{P,(M}{}^{Q(R}\delta_{N)}^{S)}
            -c_{MN}{}^{Q(R}\delta_P^{S)}\;.\nn
\end{align}
Using eqs. \eqref{vbwcrel},\eqref{Threeindexvw}, we find that the symmetrised parts
in $(MNP)$ vanish in all four expressions. Therefore, $Z_{MN,P}$ only
contains $\BB_3$ (level 3 in the Borcherds superalgebra).

Finally, consider $Z_{MN,PQ}$. 
The same relations that were used for $Z_{MN,P}$ show that
$Z_{(MN,PQ)}=0$, which leaves the plethysm \yng(2,2) .
It does not obviously vanish.
However, we have checked to high degrees in $\gl$-gradings that an object in \yng(2,2) for which any symmetric pair is in $S_2$ vanishes identically.
 There, one observes that $\vee^4R(-\lambda)$ is large enough to
contain all of $\vee^2S_2$.
Provided this holds, $Z_{MN,PQ}=0$.

We conjecture that, as a vector space (and as a $\fg^{++}$ module for $p>0$), 
\be
R_\ell=\BB_\ell\oplus\BB_{\ell+1}
\ee
for all $\ell\in\ZZ$.
The first term is level $\ell$ in the Borcherds subalgebra $\BB(\fg^{++})\subset S(\fg^{++})$.
The invariant bilinear form pairs the Borcherds part $\BB_\ell$ at level $\ell$ with the second part in
$R_{-\ell-1}=\BB_{-\ell-1}\oplus\BB_{-\ell}$. Such a bilinear form is invariant under the brackets given for $-3\leq\ell\leq2$.
The conjecture is supported by the $\gl$ grading, depicted for $S(A_1^{++})$ in Figure \ref{SA1++Figure} and Table \ref{SA1++ColumnsTable}, and for $S(E_{10})$ in Figure \ref{SE10Figure} and Table \ref{SE10ColumnsTable}.
In both cases, the two copies of the Borcherds superalgebras produce complete columns of forms, satisfying the other restrictions explained in Section \ref{InvariantFormSection}. This property has been checked for $S(E_{10})$ to degree $7$ and for $S(A_1^{++})$ to degree $10$.

\subsection{$\fg^+$-covariant double grading}

This is the grading with respect to the two leftmost nodes. 
Its central part is shown in Table \ref{DoubleGradingTable}.
The leftmost fermionic generator $\epsilon=\epsilon_{-2}$ carries degree $(p,q)=(0,-1)$, while $\epsilon_{-1}$ is the lowest weight state at $(p,q)=(1,0)$. The diagonal $p=q$ contains level $0$ in the grading with respect to the fermionic node in Figure \ref{DynkinSE} and its analogues.
The relation is $\ell=p-q$.
We find the decomposition of the $\fg^{++}$ adjoint as
$(\ldots,\bar A'^{\flat\mu\nu},\bar J^{\flat\mu},(T^A_m,\KK,\dd),E^\sh_\mu,A^\sh_{\mu\nu},\ldots)_0$ and the decomposition of the fundamental as
$(L_1,J_\mu,(S'{}^\sh_{\mu\nu},A'_{\mu\nu}),\ldots)_1$.
Since $\epsilon$ is a $\fg^+$ scalar, there is a ``lowering'' operation $X^\flat=-[X,\epsilon]$, mapping an irreducible module at $(p,q)$ to the same module at $(p,q-1)$. 
Its ``inverse'' is the raising operator $X^\sh$. All modules come in pairs, so that 
$(X^\flat)^\sh+(X^\sh)^\flat=X$ for all $X$. We refer to refs. 
\cite{Cederwall:2018aab,Cederwall:2019qnw,Cederwall:2019bai} for details in the general setting.
The generators denoted $S^{(')}$ form a basis for $s_2$, the ones denoted $A^{(')}$ form a basis for $a_2$. The number below each set of generators is the mode number shift in the action of $\dd$, as compared to how they appear in the decomposition of the over-extended adjoint.

This is the grading relevant for extended geometry with $\fg^+$ as structure algebra.
As such, it contains information on the fields and gauge symmetries of extended geometry with affine structure group \cite{Bossard:2017aae,Bossard:2018utw}. For example, the constrained ancillary fields
\cite{Cederwall:2018aab,Cederwall:2019bai} of ref. \cite{Bossard:2018utw} correspond to the position of the generators $\bar E^\mu$.

\renewcommand{\arraystretch}{1.8}		
\begin{table}
\begin{center}
\begin{tabular}{r|c|c|c|c|c}
&$p=-2$&$-1$&$0$&$1$&$2$\\\hline
$q=3$&&&&&$\shift0{A'{}^\sh_{\mu\nu}}$\\
$2$&&&$\shift{-1}{\pi^\sh}$&$\shift0{J^\sh_\mu}$
          &$\shift0{A^\sh_{\mu\nu}},\shift1{S'^\sh_{\mu\nu}},\shift1{A'_{\mu\nu}}$\\
$1$&$\shift{-2}{\bar S^{\mu\nu}}$&$\shift{-1}{\bar E^\mu}$
          &$\shift0{T^{\sh A}_{1-m}},\shift{-1}{\KK^\sh},\shift0{L_1^\sh},\shift0\pi$
          &$\shift0{E^\sh_\mu},\shift1{J_\mu}$
          &$\shift1{S^\sh_{\mu\nu}},\shift1{A_{\mu\nu}},\shift2{S'_{\mu\nu}}$\\
$0$&$\shift{-1}{\bar S^{\flat\mu\nu}},\shift{-1}{\bar A^{\mu\nu}},\shift{-2}{\bar S'^{\mu\nu}}$
          &$\shift{0}{\bar E^{\flat\mu}},\shift{-1}{\bar J^\mu}$
          &$\shift0{T^A_m},\shift0\KK,\shift1{L_1},\shift0\dd$
          &$\shift1{E_\mu}$&$\shift2{S_{\mu\nu}}$\\
$-1$&$\shift0{\bar A^{\flat\mu\nu}},\shift{-1}{\bar S'^{\flat\mu\nu}},\shift{-1}{\bar A'^{\mu\nu}}$
          &$\shift0{\bar J^{\flat\mu}}$&$\shift1\epsilon$&&\\
$-2$&$\shift0{\bar A'^{\flat\mu\nu}}$&&&&
\end{tabular}
\captionof{table}{\it The central part of the decomposition of $S(\fg^{++})$ in $\fg^+$ modules. Note the symmetry around $(p,q)=(0,\frac12)$.}
\label{DoubleGradingTable}
\end{center}
\end{table}
\renewcommand{\arraystretch}{1.3}

\appendix

\section{Calculation of modules in Borcherds superalgebras\label{BorcherdsAppendix}}
Extra modules in tensor hierarchy algebras at non-negative levels, as
compared to Borcherds superalgebras, may be predicted \eg\ in a $\gl$
grading, as described in Section \ref{ExtraModulesSubSection}. We thus need to calculate the content in a Borcherds superalgebra to some level. To this end, we use the Koszul duality established in ref. \cite{Cederwall:2015oua}.

A Borcherds superalgebra $\BB(\fa)$ is constructed from a Kac--Moody algebra $\fa$ by attaching a single fermionic (null) node
so that level $1$ in the corresponding grading is a lowest weight module $R(-\lambda)$ of $\fa$.
Let $\mu\in R(-\lambda)$ lie in the minimal orbit, meaning that functions of $\lambda$ with degree of homogeneity $p$ contain the unique lowest weight module $R(-p\lambda)$.
Then, the partition function of the Lie superalgebra $\BB_+(\fa)$, the positive level subalgebra, twisted with fermion number, is the inverse of the partition function of functions of $\mu$. This is interpreted as a denominator formula for $\BB(\fa)$.

Concretely, the partition function of functions on the minimal orbit, taking values in the representation ring, is
\be
Z_\mu(t)=\bigoplus\limits_{p=0}^\infty R(-p\lambda)t^p\;.
\ee
The partition function for $\BB_+(\fg)$ is
\be
Z_\BB(t)=\bigotimes\limits_{q=1}^\infty z_F(\BB_{2q-1},t^{2q-1})\otimes z_B(\BB_{2q},t^{2q})\;,
\ee
where 
\begin{align}
z_F(R,t)&=\bigoplus\limits_{i=0}^\infty\wedge^i R(-t)^i\;,\nn\\
z_B(R,t)&=\bigoplus\limits_{i=0}^\infty\vee^i R\, t^i=\frac1{z_F(R,t)}
\end{align}
are the partitions for a level $1$ fermion/boson in $R$, and where $\BB_\ell$ is the module at level $\ell$ in 
$\BB(\fa)$.
Then, the statement is
\be
Z_\mu(t)\otimes Z_{\BB}(t)=1\;.
\ee

This duality is useful for calculating the content of Borcherds superalgebras in a level expansion, since the modules occurring in the partition function for $\mu$ are much more tractable.
There is a single irreducible lowest weight module at each positive level. The superalgebra
$\BB(\fa)$, on the other hand, contains an infinite number or irreducible modules already at level $2$ when $\fa$ is infinite-dimensional.
If a further grading, \eg\ a $\gl$ grading, is performed on the $\fa$ modules, it will be dealt with using the method of Appendix \ref{GradingAppendix}.

\section{Weyl--Kac character formula in gradings \label{GradingAppendix}}

There are computational tools for decomposition of representations of Kac--Moody algebras, notably 
SimpLie \cite{SimpLie}, which handles infinite-dimensional algebras. We need results that go beyond the capacity of SimpLie (on a small computer).
For example, the search for extra modules in $S(E_{10})$, illustrated in Figure \ref{SE10Figure}, uses
the decomposition into $\gl(10)$ of the $E_{10}$ modules with lowest weights $-\EWeight{\ell00000}0{00}0$, $1\leq \ell\leq 22$, to degree $7$ (where the degree of the weight 
$\EWeight{100000}0{00}0$ is shifted to $0$). The reason it is sufficient to find the $\gl$ gradings of these irreducible modules at each level, is that they are used to calculate $\BB_\ell$ using the duality described in Appendix \ref{BorcherdsAppendix}.

We use the Weyl--Kac character formula, containing a sum over the Weyl group. Since we are not interested in the character formulas for the finite-dimensional $\sl(10)$ representations, but only their highest/lowest weights, it is sufficient to consider Weyl group elements that, through the transformation shifted by the Weyl vector $\varrho$,
\be
W(\Lambda)=w(\Lambda+\varrho)-\varrho\;,
\ee
map the dominant $E_{10}$ weights to weights which are ''$\sl(10)$-dominant'', \ie, have non-negative entries for all nodes except the exceptional one.
Their restriction to $\sl(10)$ weights will then be the highest/lowest weights for the corresponding $\sl(10)$ modules appearing at some degree. 
In the example at hand, the Weyl transformations with this property are listed in Table
\ref{WeylTable}, together with the images of $\EWeight{\ell 00000}0{00}0$ and their degree.
The table is of course truncated---the branching contains an infinite number of $\sl(10)$ representations---but gives the complete result to the degree considered.

\begin{table}
\begin{center}
\begin{tabular}{ccc}\hline
Weyl group element&$\varrho$-shifted image of $\EWeight{\ell 00000}0{00}0$&degree\\ \hline
$1$&$\EEWeight{\ell 00000}0{00}{0}$&$0$\\
$w_8$&$\EEWeight{\ell 00000}1{00}{-2}$&$1$\\
$w_8w_5$&$\EEWeight{\ell 00001}0{10}{-3}$&$2$\\
$w_8w_5w_4$&$\EEWeight{\ell 00010}0{20}{-4}$&$3$\\
$w_8w_5w_6$&$\EEWeight{\ell 00002}0{01}{-4}$&$3$\\
$w_8w_5w_4w_3$&$\EEWeight{\ell 00100}0{30}{-5}$&$4$\\
$w_8w_5w_4w_6$&$\EEWeight{\ell 00011}0{11}{-5}$&$4$\\
$w_8w_5w_6w_7$&$\EEWeight{\ell 00003}0{00}{-5}$&$4$\\
$w_8w_5w_4w_3w_2$&$\EEWeight{\ell 01000}0{40}{-6}$&$5$\\
$w_8w_5w_4w_3w_6$&$\EEWeight{\ell 00101}0{21}{-6}$&$5$\\
$w_8w_5w_4w_6w_7$&$\EEWeight{\ell 00012}0{10}{-6}$&$5$\\
$w_8w_5w_4w_6w_5$&$\EEWeight{\ell 00020}1{02}{-6}$&$5$\\
$w_8w_5w_4w_3w_2w_1$&$\EEWeight{\ell 10000}0{50}{-7}$&$6$\\
$w_8w_5w_4w_3w_2w_6$&$\EEWeight{\ell 01001}0{31}{-7}$&$6$\\
$w_8w_5w_4w_3w_6w_7$&$\EEWeight{\ell 00102}0{20}{-7}$&$6$\\
$w_8w_5w_4w_3w_6w_5$&$\EEWeight{\ell 00110}1{12}{-7}$&$6$\\
$w_8w_5w_4w_6w_7w_5$&$\EEWeight{\ell 00021}1{01}{-7}$&$6$\\
$w_8w_5w_4w_6w_5w_8$&$\EEWeight{\ell 00030}0{03}{-6}$&$6$\\
$w_8w_5w_4w_3w_2w_1w_0$&$\EEWeight{\ell+1,00000}0{60}{-8}$&$7$\\
$w_8w_5w_4w_3w_2w_1w_6$&$\EEWeight{\ell 10001}0{41}{-8}$&$7$\\
$w_8w_5w_4w_3w_2w_6w_7$&$\EEWeight{\ell 01002}0{30}{-8}$&$7$\\
$w_8w_5w_4w_3w_2w_6w_5$&$\EEWeight{\ell 01010}1{22}{-8}$&$7$\\
$w_8w_5w_4w_3w_6w_7w_5$&$\EEWeight{\ell 00111}1{11}{-8}$&$7$\\
$w_8w_5w_4w_3w_6w_5w_4$&$\EEWeight{\ell 00200}2{03}{-8}$&$7$\\
$w_8w_5w_4w_6w_7w_5w_6$&$\EEWeight{\ell 00030}2{00}{-8}$&$7$\\
$w_8w_5w_4w_3w_6w_5w_8$&$\EEWeight{\ell 00120}0{13}{-7}$&$7$\\
$w_8w_5w_4w_6w_7w_5w_8$&$\EEWeight{\ell 00031}0{02}{-7}$&$7$\\
\hline
\end{tabular}
\captionof{table}{\it Weyl transformations used in the decomposition of $E_{10}$ representations into $\gl(10)$ representations. The labels of simple Weyl reflections is according to the numbering of simple roots in Figure \ref{E-figur}.}
\label{WeylTable}
\end{center}
\end{table}

\begin{table}
\begin{center}
\begin{tabular}{ccc}\hline
Weyl group element&$\varrho$-shifted image of $(\ell,0,0)$&degree\\ \hline
$1$&$(\ell,0,0)$&$0$\\
$w_1$&$(\ell,2,-2)$&$1$\\
$w_1w_0$&$(\ell+1,4,-4)$&$3$\\
$w_1w_0w_{-1}$&$(0,3\ell+7,-(2\ell+6))$&$2\ell+5$\\
$w_1w_0w_1$&$(\ell+3,6,-6)$&$6$\\
$w_1w_0w_{-1}w_1$&$(2,3\ell+9,-(2\ell+8))$&$2\ell+8$\\
$w_1w_0w_1w_0$&$(\ell+6,8,-8)$&$10$\\
$w_1w_0w_{-1}w_1w_0$&$(4,3\ell+14,-(2\ell+12))$&$2\ell+14$\\
$w_1w_0w_1w_0w_1$&$(\ell+10,10,-10)$&$15$\\
\hline
\end{tabular}
\captionof{table}{\it Weyl transformations used in the decomposition of $A_1^{++}$ representations into $\gl(3)$ representations. The labels of simple Weyl reflections is according to the numbering of simple roots in Figure \ref{A1-figur}.}
\label{WeylTableA}
\end{center}
\end{table}

The degrees appearing in Table \ref{WeylTable} are independent of $\ell$, since $w_{-1}$ is not used, 
$\Lambda_{-1}$ being the only fundamental weight entering the lowest weights in question.  At level $8$ and higher one would encounter it, for the first time through $w_8w_5w_4w_3w_2w_1w_0w_{-1}$, which maps 
$\EWeight{\ell 00000}0{00}0$ to $\EEWeight{000000}0{\hspace{-17pt},\ell+7,0}{-(\ell+9)}$, appearing at
degree $\ell+8$.

Let $\{w_I\}$, $I\in\II$, for some index set $\II$, be the set of Weyl group elements that achieve this to degree $P$ in the decomposition, and $\{W_I\}$ its $\varrho$-shifted action. Let $\Lambda'$ be the restriction of the $E_{10}$ weight $\Lambda$ to an $\sl(10)$ weight and $p(\Lambda)$ its degree.
Denote an $\sl(10)$ representation with highest weight $\Lambda'$ by $r(\Lambda')$.
Then, the Weyl--Kac character formula, on its $\gl(10)$-covariant form, states that
\begin{align}
&\chi_P(\Lambda,t)\nn\\&=\left(\bigoplus_{I\in\II}(-1)^{|w_I|}r(W_I(\Lambda)')t^{p(W_I(\Lambda))}\right)
\otimes\left(\bigoplus_{I\in\II}(-1)^{|w_I|}r(W_I(0)')t^{p(W_I(0))}\right)^{-1}\;,
\end{align}
where the power of $t$ counts the degree, reproduces the correct branching up to degree $P$.
In our example, the degree is given as
$p(W(\Lambda))=-(W(\Lambda)-\Lambda,\Lambda_8)$.

Even if the list of Weyl elements is substantial, the calculation becomes efficient. The number of irreducible $\gl(10)$ modules in the branching of an $E_{10}$ module $\EWeight{\ell00000}0{00}0$ appearing up to degree $7$ is typically $\sim1500$. 
These irreducible modules are then used in order to construct $\BB_\ell$ in this grading, using the duality of Appendix \ref{BorcherdsAppendix}.
This calculation, implemented using LiE \cite{LiE}, lies behind the double grading exhibited in Figure
\ref{SE10Figure} and Table \ref{SE10ColumnsTable}, and completely analogous calculations, using the Weyl--Kac character formula for modules in the corresponding Borcherds superalgebras, are used for  $S(E_9)$ to degree $7$ (Figure \ref{SE9Figure}) and for $S(A_1^{++})$ to degree $10$ (Figure \ref{SA1++Figure}, Table \ref{SA1++ColumnsTable}).

\newpage


\providecommand{\href}[2]{#2}\begingroup\raggedright\endgroup

\end{document}